\documentclass[11pt,a4paper,reqno]{amsart}
\usepackage{tikz}
\usetikzlibrary{arrows.meta, positioning, calc, decorations.markings}
\usepackage{tikz-cd}
\usepackage[utf8]{inputenc}
\usepackage[T1]{fontenc}
\usepackage{amsmath,amssymb,amsthm,amsfonts,mathtools}
\usepackage{enumitem}
\usepackage[numbers,sort&compress]{natbib}
\usepackage{mathrsfs}
\usepackage{booktabs}
\usepackage[colorlinks=true,linkcolor=blue,citecolor=magenta,urlcolor=blue]{hyperref}
\usepackage{geometry}
\newtheorem{theorem}{Theorem}[section]
\newtheorem{lemma}[theorem]{Lemma}
\newtheorem{proposition}[theorem]{Proposition}
\newtheorem{corollary}[theorem]{Corollary}

\theoremstyle{definition}
\newtheorem{definition}[theorem]{Definition}
\newtheorem{example}[theorem]{Example}
\newtheorem{construction}[theorem]{Construction}

\theoremstyle{remark}
\newtheorem{remark}[theorem]{Remark}

\newcommand{\K}{\mathbb{K}}
\newcommand{\Z}{\mathbb{Z}}

\newcommand{\R}{\mathbb{R}}
\newcommand{\Q}{\mathbb{Q}}

\newcommand{\Hom}{\operatorname{Hom}}
\newcommand{\rad}{\operatorname{rad}}
\newcommand{\Mod}[1]{\operatorname{Mod}(#1)}
\newcommand{\md}[1]{\operatorname{mod}(#1)}

\newcommand{\Zg}{\operatorname{Zg}}
\newcommand{\PE}{\operatorname{PE}}
\newcommand{\CB}{\operatorname{CB}}
\newcommand{\KG}{\operatorname{KG}}

\newcommand{\forb}{\mathfrak{F}}
\newcommand{\Dgeo}{D_{\mathrm{geo}}}
\newcommand{\stgeo}{\operatorname{st}_{\mathrm{geo}}}

\title[Geometric Hammocks and the Ziegler Spectrum]{Geometric Hammocks: Surface Laminations and the Ziegler Spectrum of String Algebras}

\author{Shantanu Sardar}
\address{Department of Mathematics and Statistics, Shaanxi Normal University, Xi'an, China}
\email{shantanusardar17@gmail.com}

\date{\today}

\begin{document}

\begin{abstract}
We give the first complete topological, geometric, and ergodic classification of the Ziegler spectrum $\Zg(\Lambda)$ of a non-domestic string algebra $\Lambda$ over an algebraically closed field of arbitrary cardinality, via a dictionary with Thurston's theory of geodesic laminations on a cusped hyperbolic surface. This is a model-theoretic counterpart of Haiden--Katzarkov--Kontsevich's realization (\emph{Publ.\ Math.\ IHES} 2017) of stability conditions on the bounded derived category of a gentle algebra via flat surfaces: pure-injective modules correspond to admissible geodesics, and the Ziegler topology to a corridor quotient of Chabauty space. Three consequences follow. The Krull--Gabriel dimension equals the surface's topological accumulation depth, and the string closure with all higher Cantor--Bendixson strata is globally homeomorphic to an explicit corridor space, with the pushout map to $\Zg(\Lambda)$ a homeomorphism if and only if $\Lambda$ is domestic. Prefix-free dyadic blocks in the geometric perfect core verify Puninski's cardinality-free two-chain criterion, resolving the open problem of superdecomposable pure-injectives over uncountable fields for all non-domestic string and special biserial algebras. Three purely algebraic invariants on the perfect core---localized pp-shape entropy, a $K_0$-weighted pressure function satisfying an algebraic Ruelle--Walters Variational Principle, and a pp-bifurcation index with a Gauss--Bonnet simplicity obstruction---strictly separate uncountable iso-entropic and iso-pressure families, and the inverse limit of the quotient pp-lattice recovers the transversal Cantor cross-section.
\end{abstract}

\maketitle

\vspace{-0.9em}
\section{Introduction}
\label{sec:intro}

\subsection{Two central questions}
\label{subsec:two_questions}

Let $\Lambda$ be a finite-dimensional string algebra over an algebraically closed field $\K$ of arbitrary cardinality, and let $\Zg(\Lambda)$ be the Ziegler spectrum of indecomposable pure-injective right $\Lambda$-modules, topologized with the Ziegler (equivalently, the pp-) topology \cite{Ziegler1984, Prest2009}. When $\Lambda$ is \emph{non-domestic}---equivalently, when the module category $\Mod{\Lambda}$ supports exponentially growing families of string modules---the Ziegler spectrum carries infinitely many points of transfinite Cantor--Bendixson rank, its Krull--Gabriel dimension $\KG(\Lambda)$ is infinite \cite{Schroer2000b, PrestPuninski2002}, and its perfect core $\mathcal{P}(\Zg(\Lambda)) = \bigcap_\alpha D^\alpha(\Zg(\Lambda))$ is a non-empty closed subset on which the classical Cantor--Bendixson invariant is uniformly infinite. Beyond the classification of \emph{points} of $\Zg(\Lambda)$ \cite{PrestPuninski2002} and the computation of $\KG(\Lambda)$ for the domestic case \cite{Schroer2000b, LPP2018}, the global structure of $\Zg(\Lambda)$ has, until now, resisted explicit description.

This paper answers two questions.

\begin{enumerate}[label=(\roman*), leftmargin=*]
\item \emph{What is the global topological structure of $\Zg(\Lambda)$ for a non-domestic string algebra $\Lambda$, including all higher Cantor--Bendixson strata and the perfect core?}
\item \emph{Does every non-domestic string or special biserial algebra over an \emph{uncountable} algebraically closed field $\K$ admit a superdecomposable pure-injective module?}
\end{enumerate}

Question (ii) was explicitly left open by Puninski \cite[Abstract and \S4]{Puninski2008} and flagged as unresolved in Prest's MathSciNet review MR2434248. The obstruction is genuinely model-theoretic: while Puninski showed that every non-domestic string algebra has a pp-lattice of infinite width over any field, the standard extraction of a superdecomposable module via Ziegler's width theorem requires the ring to be countable, and the case of uncountable fields has remained open since 2008. We answer both questions affirmatively; in doing so, we uncover a dictionary between the model theory of $\Mod{\Lambda}$ and the hyperbolic geometry of a bordered surface with geodesic laminations, which we now describe.

\subsection{A dictionary with hyperbolic geometry: the Haiden--Katzarkov--Kontsevich parallel}
\label{subsec:hkk_parallel}

That string and gentle algebras admit geometric surface models is by now well established: Opper--Plamondon--Schroll \cite{OPS2018} and Baur--Coelho Sim\~oes \cite{BCS2020} realize finite-dimensional string modules as homotopy classes of admissible arcs on a marked ribbon surface, and canonical bases of morphism spaces are given by oriented intersections. In the influential work of Haiden--Katzarkov--Kontsevich \cite{Haiden2017}, this dictionary is carried one step further: for gentle algebras, the bounded derived category $D^b(\Lambda)$ is realized by quadratic differentials (equivalently, flat surfaces) on the ribbon surface, and Bridgeland stability conditions correspond to finite-area flat metrics, with the stability manifold conjecturally encoding the geometry of the surface's lamination space. The parallel is striking:

\medskip
\begin{center}
\renewcommand{\arraystretch}{1.30}
\small
\begin{tabular}{@{}p{6.5cm}p{9cm}@{}}
\toprule
\textbf{Haiden--Katzarkov--Kontsevich \cite{Haiden2017}} & \textbf{The present paper} \\
\midrule
Gentle algebras & String (and special biserial) algebras \\[2pt]
Flat surfaces / quadratic differentials & Cusped hyperbolic surfaces / forbidden boundaries \\[2pt]
Stability conditions on $D^b(\Lambda)$ & Indecomposable pure-injective modules in $\Zg(\Lambda)$ \\[2pt]
Stability manifold & Perfect core $\mathcal{P}(\Zg(\Lambda))$ \\[2pt]
Geodesic laminations / measured foliations & Geodesic laminations / measured laminations \\
\bottomrule
\end{tabular}
\end{center}
\medskip

Both frameworks are built on the same underlying ribbon surface $(\mathcal{S}, \mathcal{M}, \forb)$, but with different geometric structures: flat for the bounded derived category, and cusped hyperbolic for the unbounded module category. The two structures are linked by Thurston's correspondence between measured foliations (dual to flat structures via holomorphic quadratic differentials) and measured geodesic laminations on the cusped hyperbolic surface \cite{Thurston1988}. The frameworks are also complementary in what they capture: Haiden--Katzarkov--Kontsevich treat the \emph{bounded derived} category, while we treat the \emph{unbounded} module category $\Mod{\Lambda}$ and its pure-injective spectrum. Both are governed by the same underlying surface, but the model theory of $\Mod{\Lambda}$ is fundamentally a theory of \emph{limits}: pure-injective modules arise as direct limits of finite-dimensional modules (Pr\"ufer-like), as inverse limits (adic-like), and as bi-infinite string modules ``at infinity.'' Our main observation is that these limits admit precise geometric avatars: they are the complete admissible geodesics and geodesic laminations on the cusped hyperbolic surface, with the Ziegler topology corresponding to a $T_0$ corridor topology on the space of geodesics.

This dictionary allows us to ``geometrize'' the transfinite invariants of $\Zg(\Lambda)$: the Cantor--Bendixson derivative becomes \emph{geometric pruning} of pointed geodesics, the perfect core becomes the space of \emph{minimal arational geodesic laminations}, the Krull--Gabriel dimension becomes the \emph{pruning depth}, and superdecomposable pure-injective modules become pure-injective modules supported on a non-isolated closed subset of the perfect core, realized geometrically as prefixes of arational laminations. Once this correspondence is in place, the two questions (i) and (ii) admit clean geometric solutions.

\subsection{Main results}
\label{subsec:main_results}

Throughout, $\K$ denotes an algebraically closed field of \emph{arbitrary cardinality}; $\Lambda$ is a finite-dimensional string algebra modeled on a generalized marked surface $(\mathcal{S}, \mathcal{M}, \forb)$ with cellular dissection $\Delta = \{d_i\}_{i \in Q_0}$ and forbidden boundary segments $\forb$ encoding the monomial relations of $\Lambda$ (Section 2); and $\mathcal{S}^\circ = \mathcal{S} \setminus \mathcal{P}_{\mathcal{M}}$ is the punctured surface obtained by collapsing each marked segment to an ideal boundary cusp, equipped with its canonical complete finite-area hyperbolic metric. We prove five main theorems.

\medskip
\refstepcounter{maintheorem}%
\noindent\textbf{Theorem \Alph{maintheorem} (Geometric string hammock poset).}%
\label{thm:theorem_A}
\emph{For any vertex $i \in Q_0$, Schr\"oer's combinatorial hammock poset $H_\Lambda(S(i))$ is strictly order-isomorphic to the geometric fan $\mathcal{F}(\gamma_{P(S(i))})$ of pointed admissible arcs anchored at the dissection arc $d_i$. Combinatorial hook adjunction and co-hook deletion correspond geometrically to elementary boundary pivots of arcs.}

\smallskip
Theorem~\ref{thm:theorem_A} is the bridge between finite-dimensional representation theory and surface geometry; it is proved in Section 3 by translating Butler--Ringel's combinatorial classification of almost split sequences into a topological statement about elementary pivots.

\medskip
\refstepcounter{maintheorem}%
\noindent\textbf{Theorem \Alph{maintheorem} (Pure-injective limits as geodesics and laminations).}%
\label{thm:theorem_B}
\emph{Infinite monotonic pivot sequences of pointed admissible arcs converge to complete pointed admissible geodesics (forming simple geodesic laminations whenever the limiting curves are simple on $\mathcal{S}^\circ$), and the Cantor--Bendixson limit points of $\Zg(\Lambda)$ decompose into three geometric families:
\begin{enumerate}[label=(\arabic*), leftmargin=*]
\item \textbf{Global band correspondence.} Up to band inversion, unoriented primitive admissible closed geodesics $\gamma_v$ avoiding $\forb$ correspond bijectively to band generic modules $G_v \in \Zg(\Lambda)$. Decorating an oriented closed geodesic with a holonomy scalar $\lambda \in \K^\times$ parameterizes the regular tube $F_v(\Zg(\K[t, t^{-1}]))$ of finite band, Pr\"ufer, and adic modules.
\item \textbf{Local pointed periodic rays and corays.} Pointed periodic string ray and coray modules $\PE(M(w_0 v^{\infty}))$, $\PE(\widehat{M}(w_0 v^{-\infty}))$ in $\overline{H_\Lambda(S(i))}$ correspond bijectively to pointed forward and backward spiraling geodesics anchored at $p_0 \in d_i$ and accumulating asymptotically onto $\gamma_v$.
\item \textbf{Local pointed aperiodic string pure-injectives.} Pointed aperiodic string pure-injective modules $\PE(M(\vec{w}))$ in $\overline{H_\Lambda(S(i))}$ correspond bijectively to pointed complete aperiodic admissible geodesics $(\ell, p_0)$ in $\overline{\mathcal{F}(\gamma_{P(S(i))})}$, forming pointed cusp-spiraling leaves or pointed leaves of minimal arational/Cantor geodesic laminations whenever $\ell$ is simple on $\mathcal{S}^\circ$.
\end{enumerate}}

\smallskip
Theorem~\ref{thm:theorem_B} is proved in Section 5; it establishes the geometric meaning of the classification of pure-injective points of $\Zg(\Lambda)$ and provides the geometric dictionary that drives everything else.

\medskip
\refstepcounter{maintheorem}%
\noindent\textbf{Theorem \Alph{maintheorem} (Global topological classification and Krull--Gabriel dimension).}%
\label{thm:theorem_C}
\emph{The subspace Ziegler topology on the string closure $C = \bigcup_{i \in Q_0} \overline{H_\Lambda(S(i))}$ is globally homeomorphic to the $T_0$ corridor quotient space $\mathcal{C}_{\mathrm{corr}}$ of pointed complete geodesics and laminations. For every ordinal $\alpha \geq 2$ (including the perfect core), restriction induces a global homeomorphism $D^\alpha(\Zg(\Lambda)) \cong \Dgeo^\alpha(\mathcal{C}_{\mathrm{corr}})$. Globally, the closed-cover pushout map
\[
\Upsilon \colon \mathcal{C}_{\mathrm{corr}} \;\cup_{\{G_v\}}\; \bigsqcup_{v \in \operatorname{Ba}(\Lambda)} F_v\bigl(\Zg(\K[t, t^{-1}])\bigr) \;\longrightarrow\; \Zg(\Lambda)
\]
is a continuous bijection that is a homeomorphism on $C$ and on each band tube, and is a global homeomorphism if and only if $\Lambda$ is domestic. In either case, the Krull--Gabriel dimension equals the geometric pruning depth:
\[
\KG(\Lambda) \;=\; \stgeo(\mathcal{S}, \mathcal{M}, \forb) \;\coloneqq\; \max_{i \in Q_0} \CB\bigl(\overline{\mathcal{F}(\gamma_{P(S(i))})}\bigr).
\]
}

\smallskip
Theorem~\ref{thm:theorem_C} is the global topological statement answering question (i): $\Zg(\Lambda)$ is realized as an explicit quotient space of the surface's lamination space. It is proved in Section 6 by gluing the local geometric fan homeomorphisms across the finitely many vertices $Q_0$ and verifying that the finite closed cover detects openness.

\medskip
\refstepcounter{maintheorem}%
\noindent\textbf{Theorem \Alph{maintheorem} (Superdecomposable pure-injectives over arbitrary fields).}%
\label{thm:theorem_D}
\emph{Every non-domestic string or special biserial algebra $\Lambda$ over an algebraically closed field $\K$ of \emph{arbitrary} cardinality admits a superdecomposable pure-injective module. Specifically, the localized pp-lattice $\mathcal{L}_i(\Lambda)$ contains an interval in which two bounded chains of pointed string pp-formulas, each containing a densely ordered subchain order-isomorphic to $(\mathbb{Q}, <)$, freely generate a distributive sublattice.}

\smallskip
Theorem~\ref{thm:theorem_D} answers question (ii). It is proved in Section~\ref{sec:superdecomposables} by constructing prefix-free dyadic block interleavings in the geometric perfect core $\mathcal{P}_{\mathrm{geo}}$, verifying that they generate a free distributive product $\widehat{\mathcal{C}}^- \otimes \widehat{\mathcal{C}}^+ \hookrightarrow \mathcal{L}_i(\Lambda)$, and applying Puninski's cardinality-free two-chain criterion \cite[Theorem~2.3]{Puninski2008}. Corollary~\ref{cor:special_biserial_extension} extends the result to all special biserial and Brauer graph algebras.

\medskip
\refstepcounter{maintheorem}%
\noindent\textbf{Theorem \Alph{maintheorem} (Ergodic geometry and algebraic entropy on the perfect core).}%
\label{thm:theorem_E}
\emph{Let $N = \Psi(\mathcal{L}) \in \mathcal{P}_i$ be the arational pure-injective module associated with a minimal recurrent geodesic lamination $\mathcal{L} \subset \mathcal{S}^\circ \setminus \forb$. Then:
\begin{enumerate}[label=(\arabic*), leftmargin=*]
\item \textbf{Shape-Injectivity and Algebraic pp-Entropy.} The number $S_n(N)$ of distinct pairs of pp-definable subgroups of $N e_i$ cut out by active radius-$n$ oriented string templates satisfies $S_n(N) = \frac{1}{2}(N_n(N) + |P_n(N)|)$, and the \emph{localized pp-shape entropy}
\[
h^{\mathrm{shape}}(N) \;\coloneqq\; \lim_{n \to \infty} \frac{1}{2n}\log S_n(N) \;=\; h_{\mathrm{cod}}(\mathcal{L}) \;=\; \sup_{[\mu] \in \mathcal{M}_{\mathrm{erg}}(\mathcal{L})} h_{\mathbb{P}_\mu}(T)
\]
equals both the topological coding entropy of the first-return subshift and the supremum of Kolmogorov--Sinai entropies across the Choquet simplex of ergodic transverse measures.
\item \textbf{Ergodic Decomposition of Exhaustions.} Each ergodic transverse measure class $[\mu] \in \mathcal{M}_{\mathrm{erg}}(\mathcal{L})$ determines a full-$\mu$-measure union of generic leaves along which finite-prefix exhaustions converge to a constant asymptotic dimension vector $\underline{d}(\mu) = \Phi_\Delta(\mu)$ and pp-growth profile $\rho_\bullet(\mu) = \widehat{\Phi}_\Delta(\mu)$, injectively separating distinct ergodic classes.
\item \textbf{Localized pp-Shape Pressure and Algebraic Ruelle--Walters Variational Principle.} Weighting each $N$-shape class by the Grothendieck class $[M(\mathfrak{s})] \in K_0(\md{\Lambda}) \cong \Z^{Q_0}$ of its canonical witness module defines a convex pressure function $P^{\mathrm{shape}}_N \colon \R^{Q_0} \to \R$ satisfying
\[
P^{\mathrm{shape}}_N(\theta) \;=\; \sup_{[\mu] \in \mathcal{M}_{\mathrm{erg}}(\mathcal{L})} \bigl(\, h_{\mathbb{P}_\mu}(T) \;+\; \langle \theta, \underline{d}(\mu) \rangle \,\bigr).
\]
Its Legendre--Fenchel dual reconstructs the ergodic dimension simplex $\mathcal{D}(N) = \operatorname{conv}(\operatorname{Supp}_{K_0}^\infty(N))$ and the microcanonical $K_0$-shape entropy $\mathfrak{h}_{K_0}(N; d)$ purely algebraically, and strictly separates uncountable iso-entropic families in $\mathcal{P}(\Zg(\Lambda))$.
\item \textbf{Bifurcation Index and Stone--Ziegler Reconstruction.} The asymptotic rate $\kappa_{\mathrm{pp}}(N)$ of one-step interval forks satisfies the Gauss--Bonnet bound $\kappa_{\mathrm{pp}}(N) \le 6(|\mathcal{M}| - 2\chi(\mathcal{S}))$ whenever $\Phi(N)$ is simple (and $\le 4(|\mathcal{M}| - 2\chi(\mathcal{S}))$ in the triangle-free case), strictly separates infinite iso-pressure hierarchies, and yields a purely algebraic obstruction to non-self-intersection of $\Phi(N)$. Moreover, $N \in \mathcal{P}(\Zg(\Lambda))$ is a closed point if and only if $\mathcal{L}$ is minimal arational, in which case the inverse limit $\mathfrak{X}_i(N) = \varprojlim_n \mathcal{S}_n(N)$ of the shape-quotient tower of the quotient pp-lattice $\mathcal{L}_i(N)$ is a compact Cantor space canonically homeomorphic to the transversal cross-section $\mathcal{L} \cap d_i$.
\end{enumerate}}

\smallskip
Theorem~\ref{thm:theorem_E} is proved in Sections 11 and 12. It constitutes the ergodic and entropy-theoretic refinement of the geometric dictionary: on the perfect core, where the classical Cantor--Bendixson invariant is useless, we introduce three new algebraic invariants (pp-shape entropy, pp-shape pressure, and pp-bifurcation index) that are (i) purely expressible in terms of the localized pp-lattice and the Grothendieck group, (ii) strictly separating uncountable families of arational pure-injective modules, and (iii) geometrically meaningful via Thurston's measured lamination space and the Gauss--Bonnet theorem. These three invariants are \emph{strictly nested} in discriminative power: $\kappa_{\mathrm{pp}}(N)$ is not a function of $P^{\mathrm{shape}}_N$, and $P^{\mathrm{shape}}_N$ is not a function of $h^{\mathrm{shape}}(N)$. They have no direct counterpart in the classical theory of $\Zg(\Lambda)$, and we view them as the main \emph{new mathematics} arising from the dictionary.

\subsection{Techniques}
\label{subsec:techniques}

Three techniques combine to yield the results above.

\smallskip
\noindent\textbf{(i) The corridor topology and geometric evaluation of pp-formulas.} The Ziegler topology on $\Zg(\Lambda)$ has a basis of quasi-compact open sets $(\phi/\psi)$ indexed by pp-pairs. We show (Section 4) that for any string algebra, the restriction of this topology to the string closure $C$ admits a basis of \emph{pointed string test corridors}: geometric open neighborhoods $V_{\tau_\phi} \setminus K_{\tau_\psi}$ in the space of admissible geodesics, encoding the finite positive and negative constraints of a pp-pair $\phi = \sum_j \phi_{u_j}$, $\psi = \sum_k \phi_{w_k}$. This identification passes through Crawley--Boevey's monomial basis theorem for pp-definable subspaces of string modules \cite{CrawleyBoevey1989} and allows us to translate pp-evaluation into a geometric intersection problem.

\smallskip
\noindent\textbf{(ii) Prefix-free dyadic blocks on the geometric perfect core.} On a non-domestic surface, the geometric perfect core $\mathcal{P}_{\mathrm{geo}}$ is a compact perfect space of minimal arational geodesic laminations, and in particular is uncountable. We exploit this by constructing, in Section 8, an explicit binary tree of prefix-free dyadic blocks (iterated concatenations of pairs of distinct primitive band words bridged by admissible strings) inside the geometric perfect core; these blocks are uniformly bounded away from the monomial relations and hence define admissible bi-infinite strings. We then verify Birkhoff's criterion for free products of two chains in a distributive lattice, using a monomial basis argument for pointed string modules to avoid any countability hypothesis on $\K$. This is the technical heart of Theorem D.

\smallskip
\noindent\textbf{(iii) Corridor measure embedding and the algebraic Variational Principle.} For Theorems~\ref{thm:theorem_C} and~\ref{thm:theorem_E}, we combine Thurston's theory of measured geodesic laminations \cite{Thurston1988, CassonBleiler1988} with the finite-level identity (Proposition~\ref{prop:finite_level_identity}) relating pp-subquotient dimensions to geometric intersection numbers. The corridor measure map $\mathbb{P}\widehat{\Phi}_\Delta \colon \mathcal{PML}(\mathcal{S}^\circ, \forb) \hookrightarrow \mathbb{P}(\prod_{i \in Q_0} \R_{\ge 0}^{\mathcal{L}_i(\Lambda)})$ embeds the projectivized measured lamination space into an infinite product indexed by pp-pairs, providing the precise identification between ergodic measures on the lamination and asymptotic pp-growth profiles. The algebraic Ruelle--Walters Variational Principle (Theorem~\ref{thm:algebraic_variational_principle}) then follows from the classical Ruelle--Walters theorem applied to the symbolic coding of the first-return map, together with the exact identification of the topological pressure of a locally constant potential with the exponential growth rate of $K_0$-weighted shape-class partition sums.

\subsection{Organization}
\label{subsec:organization}

Section~\ref{sec:preliminaries} recalls string algebras and the generalized surface model. Section~\ref{sec:geometric_fan} proves Theorem~\ref{thm:theorem_A}. Section~\ref{sec:ziegler_convergence} introduces the pp-evaluation dictionary via test corridors. Section~\ref{sec:laminations} proves Theorem~\ref{thm:theorem_B}. Section~\ref{sec:pruning_kg} proves the Global Corridor Homeomorphism and Theorem~\ref{thm:theorem_C}, including the domestic/non-domestic pushout dichotomy. Section~\ref{sec:classification} translates $1$-domestic, $n$-domestic, and non-domestic algebra types into geometric properties of closed geodesics. Section~\ref{sec:superdecomposables} proves Theorem~\ref{thm:theorem_D} and extends it to special biserial and Brauer graph algebras. Sections~\ref{sec:stable_rank} and~\ref{sec:mapping_class_group} establish applications to the stable rank of the radical filtration and to bordered mapping class group symmetries. Sections~\ref{sec:measured} and~\ref{sec:entropy} prove Theorem~\ref{thm:theorem_E}, including the Basis Identity Lemma, the algebraic Ruelle--Walters Variational Principle, the Legendre--Fenchel reconstruction of the ergodic dimension simplex, the Gauss--Bonnet bound on the pp-bifurcation index, and the Stone--Ziegler reconstruction of the transversal Cantor cross-section. Section~\ref{sec:conclusion} summarizes the paper's contributions and discusses future directions.

\section{Preliminaries: String Algebras and Surface Models}
\label{sec:preliminaries}

In this section, we recall string algebras, their modules, and the generalized surface models of \cite{OPS2018, BCS2020, Brustle2020}. Throughout, let $\K$ be an algebraically closed field of arbitrary cardinality, and let $\Lambda$ be a finite-dimensional basic $\K$-algebra.

\subsection{String Algebras and String Modules}
\label{subsec:string_algebras}

\begin{definition}[String Algebra]
\label{def:string_algebra}
Let $\mathcal{Q} = (Q_0, Q_1, s, t)$ be a finite quiver and let $\rho$ be a finite set of monomial relations, each a path of length at least $2$. The algebra $\Lambda = \K\mathcal{Q}/\langle\rho\rangle$ is a \emph{string algebra} if:
\begin{enumerate}
    \item every vertex of $Q_0$ is the source of at most two arrows and the target of at most two arrows;
    \item for every arrow $\beta \in Q_1$, there is at most one arrow $\gamma \in Q_1$ with $s(\gamma)=t(\beta)$ and $\gamma\beta \notin \langle\rho\rangle$, and at most one arrow $\alpha \in Q_1$ with $t(\alpha)=s(\beta)$ and $\beta\alpha \notin \langle\rho\rangle$.
\end{enumerate}
\end{definition}

Let $Q_1^{-1}$ be the set of formal inverses of arrows, with $s(\alpha^{-1}) = t(\alpha)$ and $t(\alpha^{-1}) = s(\alpha)$ for each $\alpha \in Q_1$. A \emph{word} of length $n \ge 1$ is a sequence of syllables $w = l_n \dots l_1$ with $l_k \in Q_1 \cup Q_1^{-1}$ such that $t(l_k) = s(l_{k+1})$ for $1 \le k < n$ (together with the trivial words $e_i$ of length $0$ at each vertex $i$). A word is \emph{reduced} if it contains no subword of the form $\alpha\alpha^{-1}$ or $\alpha^{-1}\alpha$.

\begin{definition}[Strings and Bands]
\label{def:strings_and_bands}
A reduced word is a \emph{string} if it contains no subword in $\rho \cup \rho^{-1}$. A string $v$ of length $n \ge 1$ is \emph{cyclic} if $s(v) = t(v)$ and $v^2$ is also a string. A cyclic string $v$ is a \emph{primitive band} if it is not a proper power of a shorter string and contains both direct and inverse syllables. We denote the set of strings by $\operatorname{Str}(\Lambda)$ and the set of primitive bands up to cyclic permutation and inversion $v \sim v^{-1}$ by $\operatorname{Ba}(\Lambda)$. The algebra $\Lambda$ is \emph{domestic} if $\operatorname{Ba}(\Lambda)$ is finite, and \emph{non-domestic} otherwise.
\end{definition}

By the foundational work of Butler and Ringel \cite{ButlerRingel1987}, the indecomposable finite-dimensional $\Lambda$-modules are classified as follows.

\begin{definition}[String and Band Modules, \cite{ButlerRingel1987}]
\label{def:string_band_modules}
For $w = l_n \dots l_1 \in \operatorname{Str}(\Lambda)$, the \emph{string module} $M(w)$ is the $\K$-vector space with canonical basis $v_0, \dots, v_n$ corresponding to the vertices traversed by $w$. The right $\Lambda$-action is $v_{k-1} \cdot \alpha = v_k$ if $l_k = \alpha \in Q_1$, and $v_k \cdot \alpha = v_{k-1}$ if $l_k = \alpha^{-1} \in Q_1^{-1}$; all other arrow actions vanish. We have $M(w) \cong M(w^{-1})$.

For a primitive band $v \in \operatorname{Ba}(\Lambda)$, a scalar $\lambda \in \K^\times$, and an integer $n \ge 1$, the \emph{band module} $M(v, \lambda, n)$ assigns $\K^n$ to each vertex occurrence in $v$, where every syllable acts as the identity matrix $I_n$ except the final syllable closing the cycle, which acts via the Jordan block $J_n(\lambda)$. We have $M(v, \lambda, n) \cong M(v^{-1}, \lambda^{-1}, n)$.
\end{definition}

\subsection{String Hammocks and the Ziegler Spectrum}
\label{subsec:string_hammocks}

For a vertex $i \in Q_0$, let $S(i)$ be the simple right $\Lambda$-module at $i$, with projective cover $P(S(i)) = e_i\Lambda$. Following Schr\"oer \cite{Schroer1999}, we define the discrete hammock poset strictly on pointed finite-dimensional string modules.

\begin{definition}[String Hammock Poset, \cite{Schroer1999}]
\label{def:hammock_poset}
The \emph{string hammock poset} $H_\Lambda(S(i))$ is the locally finite partially ordered set whose elements are pointed string modules $(M(w), v_k)$, where $v_k \in M(w)e_i$ is a canonical basis vector corresponding to a non-zero graph map $P(S(i)) \to M(w)$. The partial order is generated by irreducible morphisms between string modules: $(M(w), v_k) \le (M(w'), v_{k'})$ if there is a directed path of irreducible morphisms between string modules in $\md{\Lambda}$ from $M(w)$ to $M(w')$ mapping $v_k$ to $v_{k'}$.
\end{definition}

By Prest and Puninski \cite{PrestPuninski2002}, finite-dimensional string and band modules form the set of isolated points and constitute a topologically dense subset of the Ziegler spectrum $\Zg(\Lambda)$. While finite band modules $M(v, \lambda, n)$ support non-zero morphisms from $P(S(i))$ whenever $v$ traverses $i$, any such morphism lies in the infinite radical $\rad_\Lambda^\omega$ \cite{Schroer2000a}, and hence $H_\Lambda(S(i))$ is a discrete structure consisting strictly of pointed string modules. For each vertex $i$ traversed by $v$ (or reaching $v$ via a bridge string), the band generic module $G_v$ and the infinite-dimensional string pure-injectives supported at $i$ emerge as limit points in the Ziegler topological closure $\overline{H_\Lambda(S(i))} \subset \Zg(\Lambda)$, so that globally every $G_v$ lies in $C = \bigcup_{i \in Q_0} \overline{H_\Lambda(S(i))}$.

\subsection{The Generalized Surface Model for String Algebras}
\label{subsec:surface_model}

For gentle algebras, Opper, Plamondon, and Schroll \cite{OPS2018} and Baur and Coelho Sim\~oes \cite{BCS2020} constructed a marked ribbon surface. To accommodate monomial relations of arbitrary length $k \ge 2$ that may share overlapping subpaths, we partition the boundary into marked and forbidden segments \cite{Brustle2020}.

\begin{definition}[Generalized Surface Model for a String Algebra, \cite{OPS2018, BCS2020, Brustle2020}]
\label{def:generalized_surface}
Let $\Lambda = \K\mathcal{Q}/\langle\rho\rangle$ be a string algebra. A \emph{generalized surface model} for $\Lambda$ is a compact, oriented surface $\mathcal{S}$ with non-empty boundary partitioned into three types of disjoint intervals: \emph{marked segments} $\mathcal{M}$ (admissible boundary intervals where arcs are permitted to terminate), \emph{forbidden segments} $\forb$ (spikes that admissible curves are strictly prohibited from crossing or terminating upon), and \emph{free segments} $\mathcal{E}$ (boundary intervals where arcs may pass freely during a pivot but may not anchor). The surface is equipped with a cellular dissection $\Delta = \{d_i\}_{i \in Q_0}$:
\begin{enumerate}
    \item \textbf{Ribbon Strips:} For each vertex $i \in Q_0$, take a compact oriented ribbon strip $R_i$ containing a central transversal dissection arc $d_i$ whose two endpoints lie on $\partial R_i \cap \mathcal{M}$. The \emph{syllables at $i$} --- namely the arrows $\alpha \in Q_1$ with $s(\alpha)=i$, together with the inverse arrows $\beta^{-1}$ of arrows $\beta \in Q_1$ with $t(\beta)=i$ --- partition into two disjoint sets $Q_{1,\mathrm{walk}}^+(i) \sqcup Q_{1,\mathrm{walk}}^-(i)$, according to which side of the cut $d_i$ they traverse through $R_i$.
    \item \textbf{Gluing:} For each arrow $\alpha \colon i \to j$ in $Q_1$, glue the outgoing side of $R_i$ to the incoming side of $R_j$ at a shared marked segment $M \in \mathcal{M}$ so that traversing from $d_i$ to $d_j$ around $M$ in the counter-clockwise direction corresponds to $\alpha$. Any unglued boundary intervals between marked segments that do not carry a relation become free segments $E \in \mathcal{E}$.
    \item \textbf{Forbidden Boundaries:} For each minimal monomial relation $r = \alpha_k \dots \alpha_1 \in \rho$ ($k \ge 2$), insert a forbidden segment $F_r \in \forb$ on the outer boundary of the polygon sector spanned by the sequence of dissection arcs $d_{s(\alpha_1)}, d_{t(\alpha_1)}, \dots, d_{t(\alpha_k)}$, obstructing any curve that enters across $d_{s(\alpha_1)}$ and traverses the full relation path to exit across $d_{t(\alpha_k)}$ (with local forbidden boundary slits inserted in the presence of overlapping relations to separate the respective turns).
\end{enumerate}
\end{definition}

\subsection{Geometric Realization of String and Band Modules}
\label{subsec:geometric_realization}

\begin{theorem}[Geometric Realization of String and Band Modules, \cite{OPS2018, BCS2020, Brustle2020}]
\label{thm:bls_realization}
Let $\Lambda$ be a string algebra with generalized surface model $(\mathcal{S}, \mathcal{M}, \forb)$ and cellular dissection $\Delta = \{d_i\}_{i \in Q_0}$.
\begin{enumerate}
    \item There is a canonical bijection between isomorphism classes of string modules $M(w)$ (with $w \in \operatorname{Str}(\Lambda)$ up to $w \sim w^{-1}$) and homotopy classes (relative to $\mathcal{M}$ and avoiding $\forb$) of \emph{admissible arcs} $\gamma_w$, where an admissible arc is an immersed curve $\gamma \colon [0,1] \to \mathcal{S} \setminus \forb$ with endpoints on $\mathcal{M}$ that intersects $\Delta$ minimally (without bigons) and transversally.
    \item There is a canonical bijection between primitive bands $v \in \operatorname{Ba}(\Lambda)$ (up to $v \sim v^{-1}$) and free homotopy classes of unoriented primitive, non-contractible closed curves $\gamma_v$ in $\mathcal{S} \setminus (\partial\mathcal{S} \cup \forb)$ that intersect $\Delta$ minimally. Equipping an oriented representative $\gamma_v$ with a rank-$n$ local system of monodromy $J_n(\lambda)$ for $\lambda \in \K^\times$ realizes the band module $M(v, \lambda, n)$.
\end{enumerate}
\end{theorem}

\begin{proof}
(1) Let $\gamma \colon [0,1] \to \mathcal{S} \setminus \forb$ be an admissible arc intersecting $\Delta$ minimally and transversally at points $0 < t_0 < t_1 < \dots < t_n < 1$ with $\gamma(t_m) \in d_{i_m}$. On each subinterval $[t_{m-1}, t_m]$, $\gamma$ travels inside the glued ribbon polygon between $d_{i_{m-1}}$ and $d_{i_m}$ around a shared marked segment $M \in \mathcal{M}$, recording a unique syllable $l_m \in Q_1 \cup Q_1^{-1}$ with $s(l_m) = i_{m-1}$ and $t(l_m) = i_m$. Minimal position ensures $\gamma$ forms no bigons, so the word $w(\gamma) = l_n \dots l_1$ contains no $\alpha\alpha^{-1}$ or $\alpha^{-1}\alpha$ and is reduced. Avoidance of $\forb$ ensures the word avoids $\rho \cup \rho^{-1}$, so $w(\gamma) \in \operatorname{Str}(\Lambda)$. Reversing the parameterization replaces $w(\gamma)$ by $w(\gamma)^{-1}$, giving $M(w(\gamma)) \cong M(w(\gamma)^{-1})$. Conversely, any string $w \in \operatorname{Str}(\Lambda)$ avoids $\rho \cup \rho^{-1}$ and therefore specifies a sequence of polygon crossings avoiding $\forb$, defining an admissible arc $\gamma_w$ unique up to homotopy relative to $\mathcal{M}$.

(2) A primitive band $v \in \operatorname{Ba}(\Lambda)$ is a cyclically reduced string with $s(v) = t(v)$ such that $v^2$ avoids $\rho \cup \rho^{-1}$ and $v$ is not a proper power. Identifying the initial and terminal dissection crossings in $R_{s(v)}$ closes the ribbon path into a primitive, non-contractible closed immersion $\gamma_v \colon S^1 \to \mathcal{S} \setminus (\partial\mathcal{S} \cup \forb)$ in minimal position with $\Delta$, unique up to free homotopy and orientation reversal $v \leftrightarrow v^{-1}$. Equipping $\gamma_v$ with the monodromy $J_n(\lambda)$ realizes the band module $M(v, \lambda, n)$.
\end{proof}

\subsection{Morphisms and the Auslander--Reiten Translation}
\label{subsec:morphisms_ar_translation}

\begin{definition}[Admissible Oriented Intersections]
\label{def:oriented_intersections}
For two admissible arcs $\gamma_M$ and $\gamma_N$, a transversal intersection $p \in \gamma_M \cap \gamma_N$ is an \emph{admissible oriented intersection} from $\gamma_M$ to $\gamma_N$ if the local rotation from $\gamma_M$ to $\gamma_N$ around $p$ within their shared polygon corridor follows the counter-clockwise orientation of $\mathcal{S}$ without sweeping across a forbidden segment in $\forb$.
\end{definition}

By Crawley--Boevey \cite{CrawleyBoevey1989}, the canonical basis of $\Hom_\Lambda(M(w), M(w'))$ corresponds bijectively to admissible oriented intersections from $\gamma_w$ to $\gamma_{w'}$. Since $\Hom_\Lambda(P(S(i)), M(w)) = M(w)e_i$, and the canonical basis of $M(w)e_i$ corresponds to the transversal intersections of $\gamma_w$ with $d_i$ by Definition~\ref{def:string_band_modules}, we have $\dim_\K \Hom_\Lambda(P(S(i)), M(w)) = |\gamma_w \cap d_i|$.

\begin{lemma}[The Topological Auslander--Reiten Translation, \cite{OPS2018, BCS2020}]
\label{lem:topological_ar}
Let $M(w)$ be a non-projective string module represented by the admissible arc $\gamma_w$. The Auslander--Reiten translation $\tau M(w)$ is represented by the arc $\tau(\gamma_w)$ obtained by shifting both endpoints of $\gamma_w$ by one admissible polygon boundary step in the positive (counter-clockwise) direction, avoiding $\forb$. Dually, $\tau^{-1}(\gamma_w)$ shifts both endpoints by one step in the negative (clockwise) direction.
\end{lemma}

\begin{proof}
By Butler--Ringel \cite{ButlerRingel1987}, $\tau M(w)$ is obtained combinatorially at each of the two ends of $w$ by adjoining a maximal hook $\alpha\beta_1^{-1}\dots\beta_r^{-1}$ (or deleting a co-hook). In the generalized surface model (Definition~\ref{def:generalized_surface}; compare \cite[Lemma~3.6]{BCS2020}), the boundary sector of a polygon between two consecutive marked segments, traversed counter-clockwise without crossing $\forb$, sweeps across the dissection arcs corresponding to exactly one direct arrow $\alpha$ followed by the entire maximal inverse string $\beta_1^{-1}\dots\beta_r^{-1}$. Sliding the endpoint of $\gamma_w$ to the adjacent marked segment on that polygon boundary performs this exact adjunction at both extremities.
\end{proof}

\begin{proposition}[Irreducible Morphisms as Elementary Boundary Pivots]
\label{prop:irreducible_pivot}
Let $M(w)$ be a string module represented by the admissible arc $\gamma_w$. An irreducible morphism $M(w) \to M(w')$ between string modules exists in $\md{\Lambda}$ if and only if $\gamma_{w'}$ is obtained from $\gamma_w$ by an \emph{elementary boundary pivot} (defined in \S\ref{sec:geometric_fan}): fixing one endpoint of $\gamma_w$ and sliding the other by one admissible step along the boundary of the incident polygon in $\Delta$ to the adjacent marked segment, avoiding $\forb$.
\end{proposition}

\begin{proof}
By the Butler--Ringel classification \cite{ButlerRingel1987}, irreducible morphisms between string modules are obtained by hook adjunction or co-hook deletion at one extremity of $w$ while leaving the other extremity fixed. We split into the two non-degenerate cases.

For a non-projective string module $M(w)$, the almost split sequence $0 \to M(w) \to E \to \tau^{-1} M(w) \to 0$ has middle term $E$ with at most two indecomposable summands \cite{ButlerRingel1987}, each corresponding to an irreducible morphism generated by adjoining a hook or deleting a co-hook at one extremity of $w$. By Lemma~\ref{lem:topological_ar}, $\tau^{-1}(\gamma_w)$ is obtained by shifting both endpoints of $\gamma_w$ by one admissible polygon boundary step; smoothing the crossing between $\gamma_w$ and $\tau^{-1}(\gamma_w)$ near one extremity pairs one original endpoint of $\gamma_w$ with one shifted endpoint of $\tau^{-1}(\gamma_w)$, corresponding bijectively to a pivot at the active extremity. Dually, when $M(w')$ is non-injective, the almost split sequence $0 \to \tau M(w') \to E' \to M(w') \to 0$ has middle term $E'$ with at most two indecomposable summands corresponding to pivots at the two extremities of $\gamma_{w'}$, giving the same statement for irreducible morphisms into $M(w')$.

When $M(w)$ is projective and $M(w')$ is injective simultaneously, the Butler--Ringel classification still applies: the irreducible morphisms out of the projective string module $M(w)$ correspond to the indecomposable summands of the left almost split morphism $M(w) \to \mathrm{rad}(M(w))$, which decomposes into at most two irreducible morphisms each realized by a pivot at an extremity of $\gamma_w$; dually, the irreducible morphisms into the injective string module $M(w')$ correspond to the indecomposable summands of the right almost split morphism $M(w')/\mathrm{soc}(M(w')) \to M(w')$, each realized by a pivot at an extremity of $\gamma_{w'}$. Since pivoting is invertible, in all cases we obtain the stated equivalence.
\end{proof}

\vspace{-0.2em}
\section{The Geometric Fan of Admissible Arcs and Theorem A}
\label{sec:geometric_fan}

In this section, we construct the geometric counterpart to the combinatorial string hammock posets defined by Schr\"oer \cite{Schroer1999}. By formalizing the action of elementary boundary pivots on pointed admissible arcs, we prove Theorem~\ref{thm:theorem_A}: the hammock poset $H_\Lambda(S)$ is strictly order-isomorphic to the geometric fan $\mathcal{F}(\gamma_{P(S)})$. Throughout this section, $\Lambda$ is a string algebra modeled on $(\mathcal{S}, \mathcal{M}, \forb)$ with a fixed simple module $S = S(i)$ corresponding to the dissection arc $d_i \in \Delta$.

\subsection{Pointed Admissible Arcs and Elementary Pivots}
\label{subsec:pointed_arcs}

Recall from Definition~\ref{def:hammock_poset} that the elements of $H_\Lambda(S(i))$ are pointed string modules $(M(w), v_k)$, where $v_k \in M(w)e_i$ tracks a specific occurrence of vertex $i$ in $w$ corresponding to a non-zero graph map $P(S(i)) \to M(w)$.

\begin{definition}[Pointed Admissible Arcs]
\label{def:pointed_arcs}
Let $i \in Q_0$ be a vertex and $d_i \in \Delta$ be its corresponding dissection arc. An ordered pair $(\gamma, p)$ is a \emph{pointed admissible arc} at vertex $i$ if $\gamma \colon [0,1] \to \mathcal{S} \setminus \forb$ is an admissible arc intersecting $\Delta$ minimally and transversally, and $p \in \gamma \cap d_i$ is a chosen transversal intersection point. 

Two pointed admissible arcs $(\gamma, p)$ and $(\gamma', p')$ are equivalent if there exists a continuous family of admissible arcs $\gamma_t$ ($t \in [0,1]$) such that:
\begin{enumerate}
    \item $\gamma_0 = \gamma$ and $\gamma_1 = \gamma'$;
    \item For all $t \in [0,1]$, the endpoints $\gamma_t(0)$ and $\gamma_t(1)$ remain in $\mathcal{M}$;
    \item For all $t \in [0,1]$, the image of $\gamma_t$ is strictly disjoint from $\forb$;
    \item For all $t \in [0,1]$, $\gamma_t$ intersects $\Delta$ transversally, and the specific intersection point $p_t \in \gamma_t \cap d_i$ varies continuously along $d_i$ from $p$ to $p'$.
\end{enumerate}
\end{definition}

By Theorem~\ref{thm:bls_realization}, the canonical basis vectors $v_k \in M(w)e_i$ are in strict bijection with the transversal intersections $p_k \in \gamma_w \cap d_i$. Thus, the map $(M(w), v_k) \mapsto (\gamma_w, p_k)$ is a bijection between the elements of $H_\Lambda(S(i))$ and the equivalence classes of pointed admissible arcs at $d_i$.

\begin{definition}[Elementary Boundary Pivots]
\label{def:elementary_pivots}
Let $(\gamma, p)$ be a pointed admissible arc representing $(M(w), v_k)$ and terminating at marked segments $M_a, M_b \in \mathcal{M}$. An elementary boundary pivot shifts one active endpoint of $\gamma$ in the $\tau^{-1}$-direction (Proposition~\ref{prop:irreducible_pivot}) to the adjacent marked segment on its incident polygon boundary without crossing $\forb$, while fixing the other endpoint and the anchor intersection $p \in d_i$:
\begin{enumerate}
    \item An \emph{elementary hook pivot}, denoted $P^+(\gamma, p)$, slides the active endpoint outward along the incident polygon boundary to adjoin a maximal hook $h = \alpha\beta_1^{-1}\dots\beta_r^{-1}$ at that extremity of $w$ (lengthening the dissection crossing sequence away from $p$ and inducing an irreducible monomorphism $M(w) \hookrightarrow M(wh)$).
    \item An \emph{elementary co-hook pivot}, denoted $P^-(\gamma, p)$, slides the active endpoint inward along the incident polygon boundary to delete a maximal co-hook $c = \gamma_s\dots\gamma_1\delta^{-1}$ at that extremity of $w$ (shortening the dissection crossing sequence toward $p$ and inducing an irreducible epimorphism $M(w) \twoheadrightarrow M(w')$).
\end{enumerate}
\end{definition}

\begin{definition}[The Geometric Fan]
\label{def:geometric_fan}
Let $P(S(i)) = e_i\Lambda$ be the projective cover of $S(i)$, and let $(\gamma_{P(S(i))}, p_0)$ be the pointed admissible arc representing the canonical identity evaluation map $P(S(i)) \to P(S(i))$. The \emph{geometric fan} $\mathcal{F}(\gamma_{P(S(i))})$ is the partially ordered set of homotopy classes of pointed admissible arcs generated from $(\gamma_{P(S(i))}, p_0)$ by finite sequences of elementary boundary pivots. We write $(\gamma, p) \le (\gamma', p')$ in $\mathcal{F}(\gamma_{P(S(i))})$ if $(\gamma', p')$ is obtained from $(\gamma, p)$ by a finite sequence of elementary hook pivots $P^+$ or elementary co-hook pivots $P^-$, both corresponding to directed paths of irreducible morphisms in $H_\Lambda(S(i))$.
\end{definition}

\subsection{Proof of Theorem A}
\label{subsec:proof_thm_A}

\begin{theorem}[Theorem A: The Geometric Hammock Isomorphism]
\label{thm:geometric_hammock_iso}
Let $\Lambda = \K\mathcal{Q}/\langle\rho\rangle$ be a string algebra modeled on $(\mathcal{S}, \mathcal{M}, \forb)$, and let $S = S(i)$ be a simple $\Lambda$-module. The combinatorial string hammock poset $H_\Lambda(S)$ is strictly order-isomorphic to the geometric fan $\mathcal{F}(\gamma_{P(S)})$:
\[
H_\Lambda(S) \;\cong\; \bigl(\mathcal{F}(\gamma_{P(S)}), \le\bigr).
\]
\end{theorem}

\begin{proof}
We construct the map $\Phi_i \colon H_\Lambda(S(i)) \to \mathcal{F}(\gamma_{P(S(i))})$ by $\Phi_i(M(w), v_k) = (\gamma_w, p_k)$ and verify that it is an order-isomorphism.

\emph{Step 1 (Bijection on Vertices):} By Definition~\ref{def:hammock_poset} and Theorem~\ref{thm:bls_realization}, $\Phi_i(M(w), v_k) = (\gamma_w, p_k)$ is a bijection between the pointed string modules in $H_\Lambda(S(i))$ and the pointed admissible arcs in $\mathcal{F}(\gamma_{P(S(i))})$.

\emph{Step 2 (Preservation of Covering Relations):} Let $(M(w), v_k) \lessdot (M(w'), v_{k'})$ be a covering relation in $H_\Lambda(S(i))$. By definition, this corresponds to an irreducible morphism $f \colon M(w) \to M(w')$ between string modules mapping $v_k$ to $v_{k'}$ such that the evaluation map from $P(S(i))$ is preserved. By the Butler-Ringel classification \cite{ButlerRingel1987}, $f$ is uniquely a canonical hook adjunction (inclusion) or co-hook deletion (projection) at one extremity of $w$, leaving the remainder of $w$ containing $v_k$ fixed. By Proposition~\ref{prop:irreducible_pivot}, this translates bijectively to an elementary boundary pivot $P^{\pm}(\gamma_w, p_k) = (\gamma_{w'}, p_{k'})$ fixing the interior segment containing the anchor $p_k \in d_i$. Thus, $\Phi_i$ preserves and reflects covering relations.

\emph{Step 3 (Preservation of Poset Order, Source, and Sink):} By definition, $X \le Y$ in $H_\Lambda(S(i))$ if and only if there is a finite directed path of covering relations from $X$ to $Y$. By Step 2, this corresponds bijectively to a finite sequence of elementary boundary pivots connecting $\Phi_i(X)$ to $\Phi_i(Y)$ in $\mathcal{F}(\gamma_{P(S(i))})$. The source $(P(S(i)), v_0)$ maps to $(\gamma_{P(S(i))}, p_0)$, and any maximal element of $H_\Lambda(S(i))$ (in particular, the pointed injective hull $(I(S(i)), v_{\mathrm{inj}})$ whenever $I(S(i))$ is a string module) maps to a pointed admissible arc from which no further elementary pivot is admissible without crossing $\forb$. Hence $\Phi_i$ is a strict order-isomorphism.
\end{proof}

\subsection{An Explicit Quiver Example}
\label{subsec:quiver_example}

The truncation of the geometric fan $\mathcal{F}(\gamma_{P(S(i))})$ by the forbidden segments $\forb$ directly mirrors the algebraic annihilation of paths by the monomial relations $\rho$.

\begin{figure}[htbp]
\centering
\begin{tikzpicture}[scale=1.15]
    \draw[thick] (0,0) circle (2.2cm);
    
    \draw[blue, very thick] (20:2.2) arc (20:70:2.2) node[midway, above right, black] {$M_1$};
    \draw[blue, very thick] (110:2.2) arc (110:160:2.2) node[midway, above left, xshift=-2mm, yshift=1mm, black] {$M_2$};
    \draw[blue, very thick] (200:2.2) arc (200:250:2.2) node[midway, left, black] {$M_3$};
    \draw[blue, very thick] (290:2.2) arc (290:340:2.2) node[midway, below right, black] {$M_4$};

    \draw[gray!80, very thick, dashed] (70:2.2) arc (70:110:2.2) node[midway, above, yshift=2mm, black] {$F_1$};
    \draw[gray!80, very thick, dashed] (160:2.2) arc (160:200:2.2) node[midway, left, black] {$F_2$};
    \draw[gray!80, very thick, dashed] (340:2.2) arc (340:380:2.2) node[midway, right, black] {$F_4$};

    \draw[red!80!black, very thick, dashed] (250:2.2) arc (250:290:2.2) node[midway, below=0.2cm] {$F_3 \in \forb$};
    \draw[red!80!black, very thick] (270:2.2) -- (270:0.8); 

    \draw[black, thick] (45:2.2) .. controls (90:0.9) .. (135:2.2) node[pos=0.45, above] {$\gamma_\alpha$};
    \draw[green!50!black, thick, dashed, -{Stealth[length=1.5mm]}] (45:2.2) .. controls (0,0.8) and (0,-0.8) .. (225:2.2) node[pos=0.65, right=0.2cm] {$\gamma_{\beta\alpha}$};
    \draw[green!50!black, thick, {Stealth[length=1.8mm]}-{Stealth[length=1.8mm]}] (142:1.75) arc (142:215:1.75) node[midway, above right] {\small \(P^-\)};
\end{tikzpicture}
\caption{Local polygon view of the surface model for Example~\ref{ex:quiver_pivot}. Pivoting between $M_2$ and $M_3$ across the free segment $F_2 \in \mathcal{E}$ relates $\gamma_{\beta\alpha}$ ($P(S(1)) = M(\beta\alpha)$) and $\gamma_\alpha$ ($M(\alpha)$), whereas pivoting from $M_3$ to $M_4$ is obstructed by the forbidden spike $F_3 = F_{\gamma\beta\alpha} \in \forb$.}
\label{fig:quiver_pivot}
\end{figure}
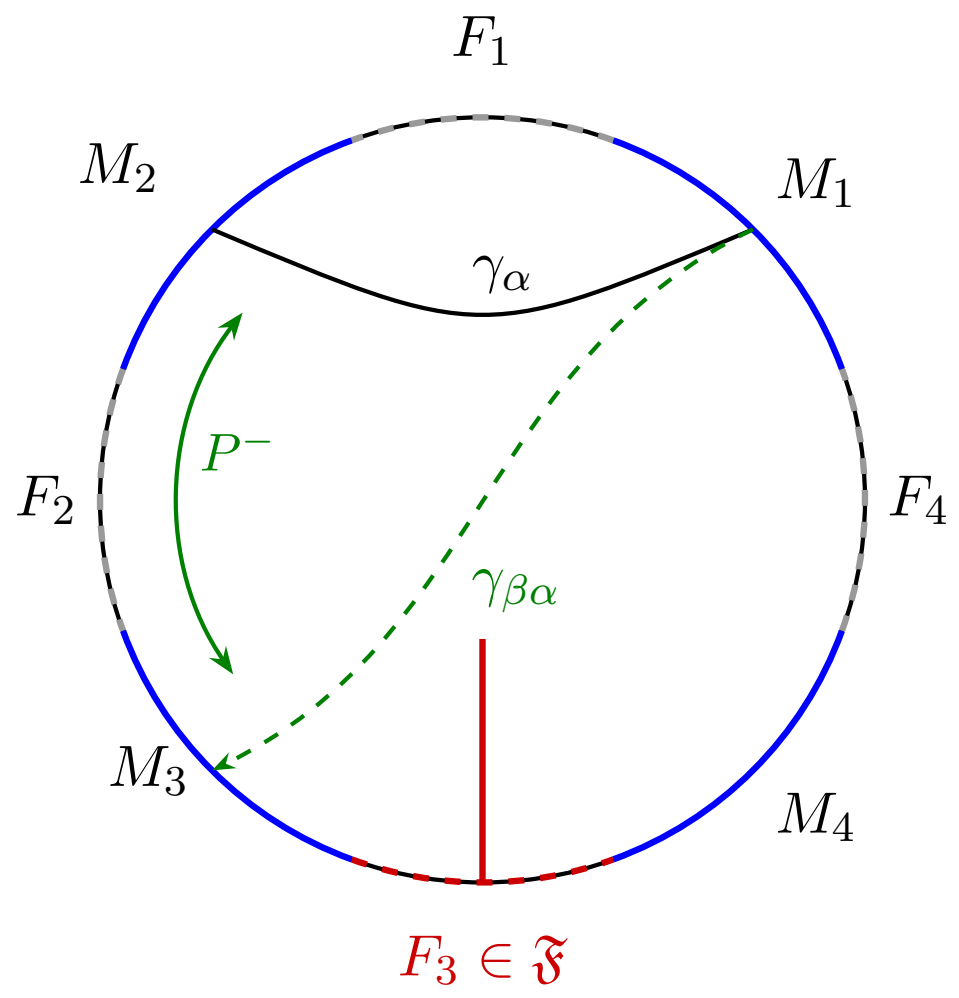

\begin{example}[Pivots and Forbidden Segments on an Ideal Quadrilateral]
\label{ex:quiver_pivot}
Consider the finite-dimensional string algebra $\Lambda = \K\mathcal{Q}/I$ given by the linear quiver
\[
1 \xrightarrow{\alpha} 2 \xrightarrow{\beta} 3 \xrightarrow{\gamma} 4
\]
with the single relation $I = \langle \gamma\beta\alpha \rangle$ (following the right-to-left composition convention). Since the relation has length $k=3$, this is a string algebra that is not gentle. Let $S = S(1)$ be the simple module at vertex $1$, whose projective cover is $P(S(1)) = e_1\Lambda = M(\beta\alpha)$ (with canonical basis $\{e_1, \alpha, \beta\alpha\}$ at vertices $1, 2, 3$, since $\gamma\beta\alpha = 0$).

The local polygon sector of the surface model $(\mathcal{S}, \mathcal{M}, \forb)$ is a topological disc ($\chi(\mathcal{S}) = 1$) whose boundary alternates between four marked segments $M_1, M_2, M_3, M_4 \in \mathcal{M}$ and four boundary intervals $F_1, F_2, F_3, F_4$ (Figure~\ref{fig:quiver_pivot}). Here, $F_3 = F_{\gamma\beta\alpha} \in \forb$ is the forbidden boundary segment encoding the relation $\gamma\beta\alpha = 0$, while $F_1, F_2, F_4 \in \mathcal{E}$ are free segments adjacent to admissible turns (e.g., $F_2$ bounds the admissible turn $\beta\alpha \notin I$). Collapsing the four marked segments $M_1, M_2, M_3, M_4$ to ideal boundary cusps turns the polygon into a hyperbolic ideal quadrilateral with positive Gauss--Bonnet area $\operatorname{Area}(\mathcal{S}^\circ) = \pi(|\mathcal{M}| - 2\chi(\mathcal{S})) = \pi(4 - 2(1)) = 2\pi > 0$.

Starting from the source arc $\gamma_{\beta\alpha}$ (representing $P(S(1)) = M(\beta\alpha)$, connecting $M_1$ to $M_3$ across the three dissection arcs $d_1, d_2, d_3$), the geometric fan $\mathcal{F}(\gamma_{P(S(1))})$ at vertex $1$ evolves as follows:
\begin{enumerate}
    \item \textbf{Admissible Co-Hook Deletion Adjacent to $F_2 \in \mathcal{E}$:} Sliding the active endpoint of $\gamma_{\beta\alpha}$ from $M_3$ to $M_2$ through the polygon interior adjacent to the free segment $F_2 \in \mathcal{E}$ is unobstructed and produces the arc $\gamma_\alpha$ (connecting $M_1$ to $M_2$ across $d_1, d_2$), which represents the two-dimensional quotient module $M(\alpha)$ under the irreducible co-hook projection $P(S(1)) = M(\beta\alpha) \twoheadrightarrow M(\alpha)$, establishing $M(\beta\alpha) \lessdot M(\alpha)$ in $H_\Lambda(S(1))$ via the co-hook pivot $P^-(\gamma_{\beta\alpha}) = \gamma_\alpha$.
    \item \textbf{Obstructed Extension at $F_3 \in \forb$:} Attempting to pivot the active endpoint of $\gamma_{\beta\alpha}$ counter-clockwise from $M_3$ toward $M_4$ through the polygon sector bounded by $F_3$ (to cross $d_4$ and form the arc that would represent $M(\gamma\beta\alpha)$ in the absence of the relation $\gamma\beta\alpha = 0$) is blocked by the forbidden boundary segment $F_3 = F_{\gamma\beta\alpha} \in \forb$, which extends into the polygon as an impenetrable spike obstructing every admissible path from $M_1$ to $M_4$.
\end{enumerate}
This geometric obstruction replicates the algebraic annihilation $\gamma\beta\alpha \in I$, preventing any extension of $P(S(1)) = M(\beta\alpha)$ across $d_4$.
\end{example}

\vspace{-0.9em}

\vspace{-0.9em}
\section{The Unbounded Category and the Ziegler Spectrum}
\label{sec:ziegler_convergence}

In Section~\ref{sec:geometric_fan}, we proved that the combinatorial hammock poset \(H_\Lambda(S)\) is strictly order-isomorphic to the geometric fan \(\mathcal{F}(\gamma_{P(S)})\). In this section, we pass to the unbounded module category \(\Mod{\Lambda}\). We review the Ziegler spectrum \(\Zg(\Lambda)\), prove that pointed string test corridors form a topological basis for the Ziegler topology on the string closure, and show that infinite monotonic pivot sequences converge topologically to pure-injective limit modules.

\subsection{Pure-Injective Modules and the Ziegler Spectrum}
\label{subsec:pure_injective}

Let \(\Mod{\Lambda}\) denote the category of all right \(\Lambda\)-modules. Recall (see \cite{Prest2009, Ziegler1984}) that a short exact sequence in \(\Mod{\Lambda}\) is \emph{pure exact} if every finitely presented module \(M \in \md{\Lambda}\) sees it as exact, and a module \(N\) is \emph{pure-injective} if every pure exact sequence starting at \(N\) splits. Every module embeds purely into a minimal pure-injective module, denoted its \emph{pure-injective envelope} \(\PE(-)\). The \emph{Ziegler spectrum} \(\Zg(\Lambda)\) is the topological space whose points are the isomorphism classes of indecomposable pure-injective modules. It is a quasi-compact \(T_0\) space with a basis of quasi-compact open sets given by
\[
\left(\frac{\phi}{\psi}\right) \;\coloneqq\; \Bigl\{\, N \in \Zg(\Lambda) \;\Big|\; \frac{\phi(N)}{\psi(N)} \neq 0 \,\Bigr\},
\]
where \(\phi/\psi\) is a pair of positive primitive (pp) formulas at some sort \(i \in Q_0\) such that \(\psi(M) \subseteq \phi(M)\) for all \(M \in \Mod{\Lambda}\). The set of pp-formulas at sort \(i\), modulo equivalence across all \(\Lambda\)-modules, forms a lattice \(\mathcal{L}_i(\Lambda)\) under conjunction \(\wedge\) (intersection of solution sets) and sum \(+\) (sum of solution sets); this lattice is modular in general. For any \(M \in \Mod{\Lambda}\), the solution set \(\phi(M)\) is a \(\K\)-subspace of \(Me_i\). By Prest and Puninski \cite{PrestPuninski2002}, the finite-dimensional string and band modules form the set of isolated points and constitute a topologically dense subset of \(\Zg(\Lambda)\).

\subsection{Geometric Evaluation of String pp-Formulas and Test Corridors}
\label{subsec:geometric_pp}

We now translate the algebraic evaluation of pp-formulas into the geometric language of the surface model \((\mathcal{S}, \mathcal{M}, \forb)\).

For a simple module \(S = S(i)\) corresponding to the dissection arc \(d_i \in \Delta\), a pointed string pp-formula \(\phi_u(x) \in \mathcal{L}_i(\Lambda)\) is determined by a finite pointed template string \((u, v_0)\) passing through vertex \(i\), where \(v_0 \in M(u)e_i\) is the distinguished canonical basis vector. Functorially, \(\phi_u(M)\) is the image of \(v_0\) under all module homomorphisms from \(M(u)\) to \(M\). If \(u\) factors into two one-sided arms \((u^-, u^+)\) exiting vertex \(i\) along the two opposite walk directions across \(d_i\), then \(\phi_u(x) = \phi_{u^-}(x) \wedge \phi_{u^+}(x)\) asserts that an element \(x \in Me_i\) extends along the finite string walk \(u\).

\begin{definition}[Geometric Test Corridors]
\label{def:test_corridor}
Let \(\phi_u(x)\) be a pointed string formula at sort \(i\), where \(u = (u^-, u^+)\) represents an admissible finite arc segment \(\tau_\phi\) anchored at \(d_i\), and let \(\Psi = \{\phi_{w_1}(x), \dots, \phi_{w_m}(x)\}\) be a finite family of pointed string formulas where each \(w_k\) strictly extends \(u\) (corresponding to compact extension arc segments \(\tau_{\psi, k}\)). The \emph{geometric test corridor} \(\mathcal{C}_{u, \Psi}\) (also denoted \(V_{\tau_\phi} \setminus K_{\tau_\psi}\)) is the localized topological condition on admissible curves \(\gamma\) crossing \(d_i\) requiring that at some transversal intersection \(p \in \gamma \cap d_i\):
\begin{enumerate}
    \item \textbf{Inclusion:} \(\gamma\) traverses the exact sequence of polygons specified by \(u\) around \(p\) (so \(\gamma\) enters \(V_{\tau_\phi}\));
    \item \textbf{Exclusion:} \(\gamma\) diverges from every extension template \(w_k \in \Psi\) strictly before completing the polygon crossings of \(w_k\) (so \(\gamma\) avoids the extension locus \(K_{\tau_\psi}\)).
\end{enumerate}
\end{definition}

\begin{lemma}[Geometric Evaluation via Test Corridors]
\label{lem:geometric_pp_evaluation}
Let \(\Lambda\) be a string algebra modeled on \((\mathcal{S}, \mathcal{M}, \forb)\), and let \(S = S(i)\) correspond to the dissection arc \(d_i \in \Delta\). Let \(\phi/\psi = \phi_u / \sum_{k=1}^m \phi_{w_k}\) be a basic open set defined by pointed string formulas at vertex \(i\), where each \(w_k\) extends \(u\). For any pointed string module \((M(w), v_r) \in H_\Lambda(S(i))\), or more generally any string pure-injective module \(N = \PE(M(\vec{w})) \in \overline{H_\Lambda(S(i))}\) with geometric realization \(\gamma_{\vec{w}}\), we have
\[
N \in \left(\frac{\phi_u}{\sum_{k=1}^m \phi_{w_k}}\right) \quad \Longleftrightarrow \quad \gamma_{\vec{w}} \in \mathcal{C}_{u, \Psi}.
\]
\end{lemma}

\begin{proof}
Because \(M(\vec{w}) \hookrightarrow \PE(M(\vec{w})) = N\) is a pure embedding, \(\phi(N)/\psi(N) \neq 0\) if and only if \(\phi(M(\vec{w}))/\psi(M(\vec{w})) \neq 0\) (pure embeddings preserve and reflect the non-triviality of pp-pairs; see \cite[\S 4.1]{Prest2009}). By Crawley-Boevey \cite{CrawleyBoevey1989} and Prest–Puninski \cite{PrestPuninski2002}, pp-definable subspaces of \(M(\vec{w})e_i\) are monomial with respect to the canonical basis. Thus \(\phi_u(M(\vec{w})) \supsetneq \sum_{k=1}^m \phi_{w_k}(M(\vec{w}))\) if and only if there exists a canonical basis vector \(v_r \in M(\vec{w})e_i\) such that \(v_r \in \phi_u(M(\vec{w}))\) and \(v_r \notin \phi_{w_k}(M(\vec{w}))\) for all \(k\).

Evaluating \(\phi_u(x)\) at \(v_r\) is equivalent to the existence of a pointed graph map \(f \colon M(u) \to M(\vec{w})\) mapping \(v_0\) to \(v_r\). In the surface model, this holds if and only if the pointed curve \((\gamma_{\vec{w}}, p_r)\) tracks the polygon corridor \(V_{\tau_\phi}\) of \(u\) around \(p_r \in d_i\). The condition \(v_r \notin \phi_{w_k}(M(\vec{w}))\) holds if and only if this homomorphism cannot factor through any of the extension modules \(M(w_k)\), which geometrically means that \(\gamma_{\vec{w}}\) diverges from each extension corridor \(K_{\tau_{\psi, k}}\) before completion. Thus \(N \in (\phi_u / \sum_{k=1}^m \phi_{w_k})\) if and only if \(\gamma_{\vec{w}} \in \mathcal{C}_{u, \Psi}\).
\end{proof}

We now show that these pointed string test corridors form a topological basis for the subspace Ziegler topology on the global string closure \(C \coloneqq \bigcup_{i \in Q_0} \overline{H_\Lambda(S(i))}\). Recall that while \(\mathcal{L}_i(\Lambda)\) is modular in general, its evaluation on string and band modules is distributive \cite[Theorem~2.3]{PrestPuninski2002}, a fact used in the proof below.

\begin{lemma}[String Corridors Form a Basis]
\label{lem:corridors_form_basis}
The basic open sets of the form \((\phi_u / \sum_{k=1}^m \phi_{w_k})\), corresponding to the geometric test corridors \(\mathcal{C}_{u, \Psi}\), form a topological basis for the subspace Ziegler topology on \(C = \bigcup_{i \in Q_0} \overline{H_\Lambda(S(i))} \subseteq \Zg(\Lambda)\).
\end{lemma}

\begin{proof}
Let \(U = (\phi / \psi)\) be an arbitrary basic Ziegler open set at sort \(i\), and let \(N \in U \cap C\). By Prest and Puninski \cite[Theorem~2.3]{PrestPuninski2002}, the localized pp-lattice \(\mathcal{L}_i(\Lambda)\) evaluated on string and band modules is distributive and generated by one-sided pointed string formulas. Every pp-formula can therefore be written as a finite sum of meets of one-sided string formulas. Because a meet of two one-sided string formulas on opposite walk arms is a two-sided pointed string formula \(\phi_u\) (while meets of one-sided string formulas on the same arm reduce to a single one-sided string formula by \cite[Lemma~2.3]{PrestPuninski2002}), we may express \(\phi = \sum_{j=1}^{r} \phi_{u_j}\) and \(\psi = \sum_{l=1}^{s} \phi_{\tilde{w}_l}\), where \(u_j\) and \(\tilde{w}_l\) are two-sided finite string templates anchored at \(i\).

Since \(N \in (\phi/\psi) \cap C\), we have \(\phi(N) \supsetneq \psi(N)\). Writing \(N = \PE(M(\vec{w}))\), purity and monomiality of string pp-evaluation imply that there exists a canonical basis element \(e_x\) in \(M(\vec{w})e_i\) such that \(e_x \in \phi(M(\vec{w})) \setminus \psi(M(\vec{w}))\). Consequently, \(e_x \in \phi_{u_j}(M(\vec{w}))\) for some \(j\), while \(e_x \notin \phi_{\tilde{w}_l}(M(\vec{w}))\) for all \(l\). Hence \(N\) belongs to the basic open set \(B_j \coloneqq (\phi_{u_j} / (\psi \wedge \phi_{u_j})) \subseteq (\phi/\psi) = U\).

Because every pointed string formula \(\phi_u(M(\vec{w})) \subseteq M(\vec{w})e_i\) is spanned by a subset of the canonical basis of \(M(\vec{w})e_i\) \cite{CrawleyBoevey1989}, and the coordinate subspaces of a vector space spanned by subsets of a fixed canonical basis form a distributive lattice, intersection distributes over finite sums of string pp-formulas on every module in \(C\):
\[
(\psi \wedge \phi_{u_j})\bigl(M(\vec{w})\bigr) \;=\; \sum_{l=1}^s \bigl(\phi_{\tilde{w}_l} \wedge \phi_{u_j}\bigr)\bigl(M(\vec{w})\bigr) \;=\; \sum_{k=1}^m \phi_{w_k}\bigl(M(\vec{w})\bigr),
\]
where \(\{w_k\}_{k=1}^m\) is the finite set of minimal common string extensions of \(u_j\) and \(\tilde{w}_1, \dots, \tilde{w}_s\) strictly extending \(u_j\). Thus \(B_j \cap C = (\phi_{u_j} / \sum_{k=1}^m \phi_{w_k}) \cap C\) is a geometric test corridor neighborhood of \(N\) contained in \(U \cap C\).
\end{proof}

\subsection{Infinite Pivots and Topological Convergence}
\label{subsec:infinite_pivots}

\begin{definition}[Infinite Pivot Sequences]
\label{def:infinite_pivots}
An \emph{infinite pivot sequence} in \(\mathcal{F}(\gamma_{P(S(i))})\) is an infinite sequence of pointed admissible arcs \((\gamma_0, p_0) \xrightarrow{P^\pm} (\gamma_1, p_1) \xrightarrow{P^\pm} \dots\) such that the anchor intersections \(p_n \in d_i\) remain fixed and the combinatorial word length of the arcs \(\gamma_n\) (i.e., the number of syllables in the corresponding words) tends to infinity. We say the sequence is \emph{monotonic} if at each active extremity the pivots proceed monotonically (either via forward hook adjunctions \(P^+\) or backward co-hook deletions \(P^-\)).
\end{definition}

\begin{theorem}[Topological Convergence of Infinite Pivots]
\label{thm:topological_convergence}
Let \(\Gamma = ((\gamma_n, p_n))_{n \ge 0}\) be an infinite monotonic pivot sequence consisting of forward hook adjunctions in \(\mathcal{F}(\gamma_{P(S(i))})\), corresponding to pointed string modules \((M(w_n), v_0)\). Then:
\begin{enumerate}
    \item The direct limit \(M_\infty = \varinjlim M(w_n)\) is an infinite-dimensional string module \(M(\vec{w})\) (and dually, a monotonic sequence of backward co-hook projections yields an inverse limit coray module \(\widehat{M}(\vec{w})\)).
    \item In the Ziegler topology on \(\Zg(\Lambda)\), the sequence of finite-dimensional string modules \(M(w_n)\) converges to the pure-injective envelope \(\PE(M_\infty)\).
\end{enumerate}
\end{theorem}

\begin{proof}
(1) Since the pivot sequence consists of monotonic forward hook adjunctions fixing the anchor \(v_0\), the induced morphisms \(M(w_n) \hookrightarrow M(w_{n+1})\) are canonical monomorphisms in \(\md{\Lambda}\), whose direct limit \(M_\infty = \varinjlim M(w_n) = M(\vec{w})\) is an infinite-dimensional string module.

(2) Let \((\phi/\psi)\) be a basic Ziegler open set containing \(\PE(M_\infty)\). By Lemma~\ref{lem:corridors_form_basis}, we may assume \((\phi/\psi) = (\phi_u / \sum_{k=1}^m \phi_{w_k})\) is a pointed string test corridor. Since \(M_\infty \hookrightarrow \PE(M_\infty)\) is pure, \(\phi(M_\infty)/\psi(M_\infty) \neq 0\). Because pp-formulas commute with directed colimits and the templates \(u, w_1, \dots, w_m\) have bounded finite radius \(R\), there exists an integer \(N > 0\) such that for all \(n \ge N\), the initial segment of \(w_n\) within radius \(R\) of the anchor stabilizes permanently to that of \(\vec{w}\). By Lemma~\ref{lem:geometric_pp_evaluation}, for all \(n \ge N\), the pointed arc \((\gamma_n, p_n)\) permanently traverses the test corridor \(V_{\tau_\phi} \setminus K_{\tau_\psi}\), so \(M(w_n) \in (\phi/\psi)\) for all \(n \ge N\).
\end{proof}

\vspace{-0.2em}
\section{Pure-Injective Limits and Generic Modules}
\label{sec:laminations}

In Section~\ref{sec:ziegler_convergence}, we established that infinite monotonic pivot sequences in the geometric fan converge topologically to pure-injective modules in the Ziegler spectrum. In this section, we equip the punctured surface \(\mathcal{S}^\circ\) with a complete finite-area hyperbolic metric and invoke Thurston's theory of geodesic laminations \cite{Thurston1988, CassonBleiler1988}. We prove Theorem~\ref{thm:theorem_B_sec5} (Main Theorem~\ref{thm:theorem_B}): the Cantor-Bendixson limit points of the hammock posets in \(\Zg(\Lambda)\) correspond bijectively to complete pointed geodesics and geodesic laminations on \((\mathcal{S}, \mathcal{M}, \forb)\).

\subsection{Algebraic Classification of Pure-Injective Limit Modules}
\label{subsec:algebraic_limits}

\begin{definition}[Infinite Hammock Paths and Limit Modules]
\label{def:infinite_hammock_paths}
An \emph{infinite directed ray} in the hammock poset \(H_\Lambda(S(i))\) is an infinite sequence of covering relations generated by elementary pivots at the extremities of the string:
\[
\Gamma \colon (M(w_0), v_0) \;\xrightarrow{f_1}\; (M(w_1), v_0) \;\xrightarrow{f_2}\; (M(w_2), v_0) \;\xrightarrow{f_3}\; \dots
\]
The limit along a one-sided infinite sequence defines a one-sided infinite word \(\vec{w}\), while extending at both extremities defines a bi-infinite admissible word \(\overleftrightarrow{w}\).
\end{definition}

By Ziegler's completeness theorem \cite{Ziegler1984} and the classification of Prest and Puninski \cite{PrestPuninski2002, Prest2009}, the infinite-dimensional indecomposable pure-injective modules in the Ziegler closure \(\overline{H_\Lambda(S(i))} \subset \Zg(\Lambda)\) fall into three structural families:
\begin{enumerate}
    \item \textbf{Periodic String Ray and Coray Modules:} Let \(v \in \operatorname{Ba}(\Lambda)\) be a primitive band and \(w_0\) a finite bridge string such that \(w_0 v^\infty\) is admissible. The forward periodic rays (Pr\"ufer-like) are direct limits \(M(w_0 v^\infty) = \varinjlim M(w_0 v^n)\) of canonical hook inclusions, yielding \(\PE(M(w_0 v^\infty))\). The backward periodic corays (Adic-like) are inverse limits \(\widehat{M}(w_0 v^{-\infty}) = \varprojlim M(w_0 v^{-n})\) of canonical co-hook projections, yielding \(\PE(\widehat{M}(w_0 v^{-\infty}))\). In \(\mathcal{L}_i(\Lambda)\), these define distinct points in \(\overline{H_\Lambda(S(i))}\).
    \item \textbf{Band Generic Modules:} For each primitive band \(v\) reachable from vertex \(i\), the Ziegler closure of the periodic family over \(v\) contains a unique endofinite generic module \(G_v\) \cite{Krause1998}, sitting at \(\overline{H_\Lambda(S(i))} \cap F_v(\Zg(\K[t, t^{-1}]))\).
    \item \textbf{Aperiodic String Pure-Injective Modules:} When a one-sided or bi-infinite word \(\vec{w}\) is aperiodic, \(\PE(M(\vec{w}))\) yields an aperiodic infinite-dimensional string pure-injective module that does not specialize to any band generic module.
\end{enumerate}

\subsection{Cusped Hyperbolic Geometry, Admissible Geodesics, and Laminations}
\label{subsec:hyperbolic_geometry}

To translate infinite words into continuous surface geometry, we endow the generalized surface model \((\mathcal{S}, \mathcal{M}, \forb)\) with a canonical hyperbolic metric \cite{Thurston1988, CassonBleiler1988}. Let \(\mathcal{P}_{\mathcal{M}}\) denote the finite union of marked boundary segments \(M \in \mathcal{M}\). Collapsing each marked segment to an ideal boundary cusp \(c_M\) produces the punctured surface with boundary
\[
\mathcal{S}^\circ \;\coloneqq\; \mathcal{S} \setminus \mathcal{P}_{\mathcal{M}},
\]
whose boundary \(\partial\mathcal{S}^\circ = \partial\mathcal{S} \setminus \mathcal{P}_{\mathcal{M}} = \forb \cup \mathcal{E}\) consists of \(|\mathcal{M}|\) open intervals running between the cusps. Equipping \(\mathcal{S}^\circ\) with a complete hyperbolic metric of curvature \(-1\) in which every unmarked boundary interval in \(\forb \cup \mathcal{E}\) is totally geodesic and each of the \(|\mathcal{M}|\) cusps has interior angle \(0\) (equivalently exterior angle $\pi$), the Gauss--Bonnet identity (see \cite[Chapter~1]{CassonBleiler1988}) gives
\[
2\pi\chi(\mathcal{S}) \;=\; -\operatorname{Area}(\mathcal{S}^\circ) \;+\; \sum_{M \in \mathcal{M}}\pi \;+\; 0,
\]
where the last $0$ is the total geodesic curvature contribution along the totally geodesic boundary, and the $\pi$ terms arise from the exterior angles of the cusps. Solving gives
\[
\operatorname{Area}(\mathcal{S}^\circ) \;=\; \pi\bigl(|\mathcal{M}| - 2\chi(\mathcal{S})\bigr).
\]
Equivalently, the Schwarz reflection double of \(\mathcal{S}^\circ\) across its totally geodesic unmarked boundary \(\partial\mathcal{S} \setminus \mathcal{P}_{\mathcal{M}} = \forb \cup \mathcal{E}\) is a punctured surface without boundary of Euler characteristic \(2\chi(\mathcal{S}) - |\mathcal{M}|\) and hyperbolic area \(2\pi(|\mathcal{M}| - 2\chi(\mathcal{S}))\). Consequently, \(\mathcal{S}^\circ\) admits a complete, finite-area hyperbolic metric of constant curvature \(-1\) with totally geodesic boundary \(\forb \cup \mathcal{E}\) if and only if \(|\mathcal{M}| - 2\chi(\mathcal{S}) > 0\). For every representation-infinite string algebra \(\Lambda\), the quiver contains a cycle, so \(\chi(\mathcal{S}) \le 0\) on the cycle-supporting component and \(|\mathcal{M}| \ge 2\); hence \(|\mathcal{M}| - 2\chi(\mathcal{S}) \ge 2 > 0\) and \(\operatorname{Area}(\mathcal{S}^\circ) \ge 2\pi > 0\) hold unconditionally (even when \(\chi(\mathcal{S}) = 0\), as for the Kronecker cylinder with \(|\mathcal{M}| = 2\) boundary cusps, whose Schwarz reflection double is a twice-punctured torus of Euler characteristic \(-2\)).

Under this metric, every homotopy class of admissible arcs avoiding \(\forb\) contains a unique geodesic representative terminating at the cusps, and every free homotopy class of primitive closed curves avoiding \(\forb\) contains a unique closed geodesic \(\gamma_v\).

\begin{definition}[Complete Admissible Geodesics and Geodesic Laminations, \cite{Thurston1988, CassonBleiler1988}]
\label{def:geodesic_laminations}
An \emph{admissible complete geodesic} on \(\mathcal{S}^\circ \setminus \forb\) is a complete immersed geodesic \(\ell \colon \R \to \mathcal{S}^\circ \setminus \forb\) (bi-infinite in the interior, half-infinite anchored at a cusp, or closed) in minimal position with the cellular dissection \(\Delta\).
\begin{enumerate}
    \item \textbf{Simple Geodesics and Geodesic Laminations:} A complete geodesic \(\ell\) is \emph{simple} if it has no self-intersections. A \emph{geodesic lamination} \(\mathcal{L} \subset \mathcal{S}^\circ \setminus \forb\) is a closed subset that is a disjoint union of simple complete geodesics (its \emph{leaves}). It is \emph{minimal} if every leaf is dense in \(\mathcal{L}\).
    \item \textbf{General Self-Intersecting Geodesics:} Complete admissible geodesics with transversal self-intersections correspond to pure-injective limits of strings with overlapping subwords; while not simple laminations, they are uniquely encoded by their bi-infinite sequence of crossings with \(\Delta\).
\end{enumerate}
\end{definition}

By the Casson--Bleiler structure theorem \cite{CassonBleiler1988}, any geodesic lamination \(\mathcal{L}\) decomposes into a finite union of minimal sub-laminations \(\mathcal{L}_{\min}\) together with finitely many isolated non-minimal spiraling leaves \(\ell_{\mathrm{spiral}}\).

\begin{lemma}[Hausdorff and Corridor Limits of Infinite Pivoting Arcs]
\label{lem:hausdorff_limits}
Let \(\Gamma = ((\gamma_n, p_0))_{n \ge 0}\) be an infinite sequence of pointed geodesic admissible arcs generated by monotonic pivots on \((\mathcal{S}, \mathcal{M}, \forb)\). As \(n \to \infty\), the pointed arcs \((\gamma_n, p_0)\) converge to a complete pointed admissible geodesic \((\ell_\infty, p_0)\) (and when the arcs are simple, \(\{\gamma_n\}\) converges in the Hausdorff metric to a geodesic lamination \(\mathcal{L}_\infty = \mathcal{L}_{\mathrm{spiral}} \cup \mathcal{L}_{\min}\)). The minimal asymptotic component belongs to one of three geometric classes:
\begin{enumerate}
    \item \textbf{Closed Geodesic:} A primitive closed geodesic \(\gamma_v\) onto which the spiraling leaf accumulates.
    \item \textbf{Cusp Spiral:} A geodesic ray tending asymptotically into an ideal boundary cusp.
    \item \textbf{Arational / Cantor Lamination:} An uncountable minimal geodesic lamination with no closed leaves (or an aperiodic complete geodesic in a Cantor subshift of communicating closed geodesics).
\end{enumerate}
\end{lemma}

\begin{proof}
Because each \(\gamma_n\) represents an admissible string avoiding \(\rho\), the geodesic representatives remain strictly disjoint from the totally geodesic forbidden boundary \(\forb\). Lift to the universal cover \(\widetilde{\mathcal{S}^\circ} \subset \mathbb{H}^2\), anchored at a lift \(\tilde{p}_0 \in \tilde{d}_i\); the cellular dissection \(\Delta\) lifts to a tessellation of \(\mathbb{H}^2\) by ideal polygons. Since the combinatorial word length of \(\gamma_n\) tends to infinity as \(n \to \infty\), the lifted arcs \(\tilde{\gamma}_n\) are \emph{proper} geodesic rays in \(\mathbb{H}^2\) emanating from \(\tilde{p}_0\). Their terminal ideal endpoints on the compact space \(\partial_\infty\mathbb{H}^2\) therefore form a sequence that, by monotonicity of the pivot sequence and the standard contraction lemma for convex-cocompact tessellations (see \cite[Chapter~1]{CassonBleiler1988}), is Cauchy in the visual boundary. Passing to the limit gives a unique pair of ideal endpoints, defining a unique complete pointed geodesic \((\ell_\infty, p_0)\). When the arcs are simple, Thurston's compactness theorem \cite{Thurston1988, CassonBleiler1988} ensures Hausdorff convergence to a geodesic lamination \(\mathcal{L}_\infty = \mathcal{L}_{\mathrm{spiral}} \cup \mathcal{L}_{\min}\), and the classification of \(\mathcal{L}_{\min}\) follows from \cite{CassonBleiler1988}.
\end{proof}

\subsection{The Correspondence: Proof of Theorem B}
\label{subsec:proof_theorem_b}

We now prove Main Theorem~\ref{thm:theorem_B}.

\begin{theorem}[Theorem B: Pure-Injective Limits and Generic Modules as Geodesics and Laminations]
\label{thm:theorem_B_sec5}
Let \(\Lambda\) be a representation-infinite string algebra modeled on \((\mathcal{S}, \mathcal{M}, \forb)\) with cusped hyperbolic surface \(\mathcal{S}^\circ = \mathcal{S} \setminus \mathcal{P}_{\mathcal{M}}\). There is a canonical correspondence between the Cantor-Bendixson limit points in \(\Zg(\Lambda)\) and complete admissible geodesics on \(\mathcal{S}^\circ \setminus \forb\):
\begin{enumerate}
    \item \textbf{Global Band Correspondence:} Up to band inversion \(v \sim v^{-1}\), unoriented primitive admissible closed geodesics \(\gamma_v\) avoiding \(\forb\) correspond bijectively to the band generic modules \(G_v \in \Zg(\Lambda)\). Decorating an oriented closed geodesic \(\gamma_v\) with a holonomy scalar \(\lambda \in \K^\times\) parameterizes the regular tube of finite band modules \(M(v, \lambda, n)\), band Pr\"ufer modules \(M(v, \lambda, \infty)\), and band adic modules \(\widehat{M}(v, \lambda)\).
    \item \textbf{Local Pointed Periodic Rays and Corays:} For each vertex \(i \in Q_0\), pointed periodic string ray modules \(\PE(M(w_0 v^\infty))\) and coray modules \(\PE(\widehat{M}(w_0 v^{-\infty}))\) in \(\overline{H_\Lambda(S(i))}\) correspond bijectively to pointed forward and backward spiraling geodesics \((\ell^\pm_{\mathrm{spiral}}, p_0)\) anchored at \(p_0 \in d_i\) and accumulating asymptotically onto \(\gamma_v\).
    \item \textbf{Local Pointed Aperiodic String Pure-Injectives:} Pointed aperiodic string pure-injective modules \(\PE(M(\vec{w}))\) in \(\overline{H_\Lambda(S(i))}\) correspond bijectively to pointed complete aperiodic admissible geodesics \((\ell, p_0)\) in \(\overline{\mathcal{F}(\gamma_{P(S(i))})}\) (forming pointed cusp-spiraling leaves or pointed leaves of minimal arational/Cantor geodesic laminations whenever \(\ell\) is simple on \(\mathcal{S}^\circ\)).
\end{enumerate}
\end{theorem}

\begin{proof}
\emph{(1) Global Band Correspondence:} By Theorem~\ref{thm:bls_realization}(2), \(v \in \operatorname{Ba}(\Lambda)\) corresponds to a unique free homotopy class of unoriented primitive closed curves, containing a unique unoriented closed geodesic \(\gamma_v\). Because \(G_v \cong G_{v^{-1}}\) is independent of the continuous parameter \(\lambda \in \K^\times\), unoriented closed geodesics \(\gamma_v\) correspond bijectively to the band generic modules \(G_v\). Choosing an orientation on \(\gamma_v\) and a scalar \(\lambda \in \K^\times\) specifies the regular tube over \((v, \lambda)\) containing \(M(v, \lambda, n)\), \(M(v, \lambda, \infty)\), and \(\widehat{M}(v, \lambda)\), with \(\gamma_{v^{-1}}\) matching the algebraic isomorphism \(M(v, \lambda, n) \cong M(v^{-1}, \lambda^{-1}, n)\).

\emph{(2) Local Pointed Periodic Rays and Corays:} Let \(\PE(M(w_0 v^\infty)) \in \overline{H_\Lambda(S(i))}\) be anchored at \(v_0 \in M(w_0 v^\infty)e_i\). By Theorem~\ref{thm:topological_convergence} and Lemma~\ref{lem:hausdorff_limits}, the finite prefixes \(M(w_0 v^n)\) correspond to pointed arcs \((\gamma_n, p_0)\) that converge to a unique pointed geodesic ray \((\ell^+_{\mathrm{spiral}}, p_0)\) starting at an ideal cusp, crossing \(d_i\) at \(p_0\), and spiraling asymptotically onto \(\gamma_v\) in the positive orientation of \(v\). Dually, the backward periodic coray \(\PE(\widehat{M}(w_0 v^{-\infty}))\) is the limit of co-hook projections along \(v^{-1}\), corresponding to the pointed geodesic \((\ell^{-}_{\mathrm{spiral}}, p_0)\) spiraling onto \(\gamma_v\) in the opposite orientation \(v^{-1}\).

\emph{(3) Local Pointed Aperiodic String Pure-Injectives:} Let \(\PE(M(\vec{w})) \in \overline{H_\Lambda(S(i))}\) be anchored at \(d_i\). By Lemma~\ref{lem:hausdorff_limits}, the bi-infinite sequence of polygon crossings of \(\vec{w}\) lifts to \(\mathbb{H}^2\) to determine a unique complete geodesic lift \(\tilde{\ell}\) through \(\tilde{p}_0 \in \tilde{d}_i\), projecting to a unique pointed complete aperiodic admissible geodesic \((\ell, p_0) \in \overline{\mathcal{F}(\gamma_{P(S(i))})}\). Conversely, any such \((\ell, p_0)\) intersects \(\Delta\) in an admissible non-periodic syllable sequence \(\vec{w}_{(\ell, p_0)}\) avoiding \(\rho \cup \rho^{-1}\), whose pure-injective hull \(\PE(M(\vec{w}_{(\ell, p_0)}))\) recovers the pointed pure-injective module.
\end{proof}

\subsection{Example: The Kronecker Lamination on a Two-Cusped Annulus}
\label{subsec:kronecker_example}

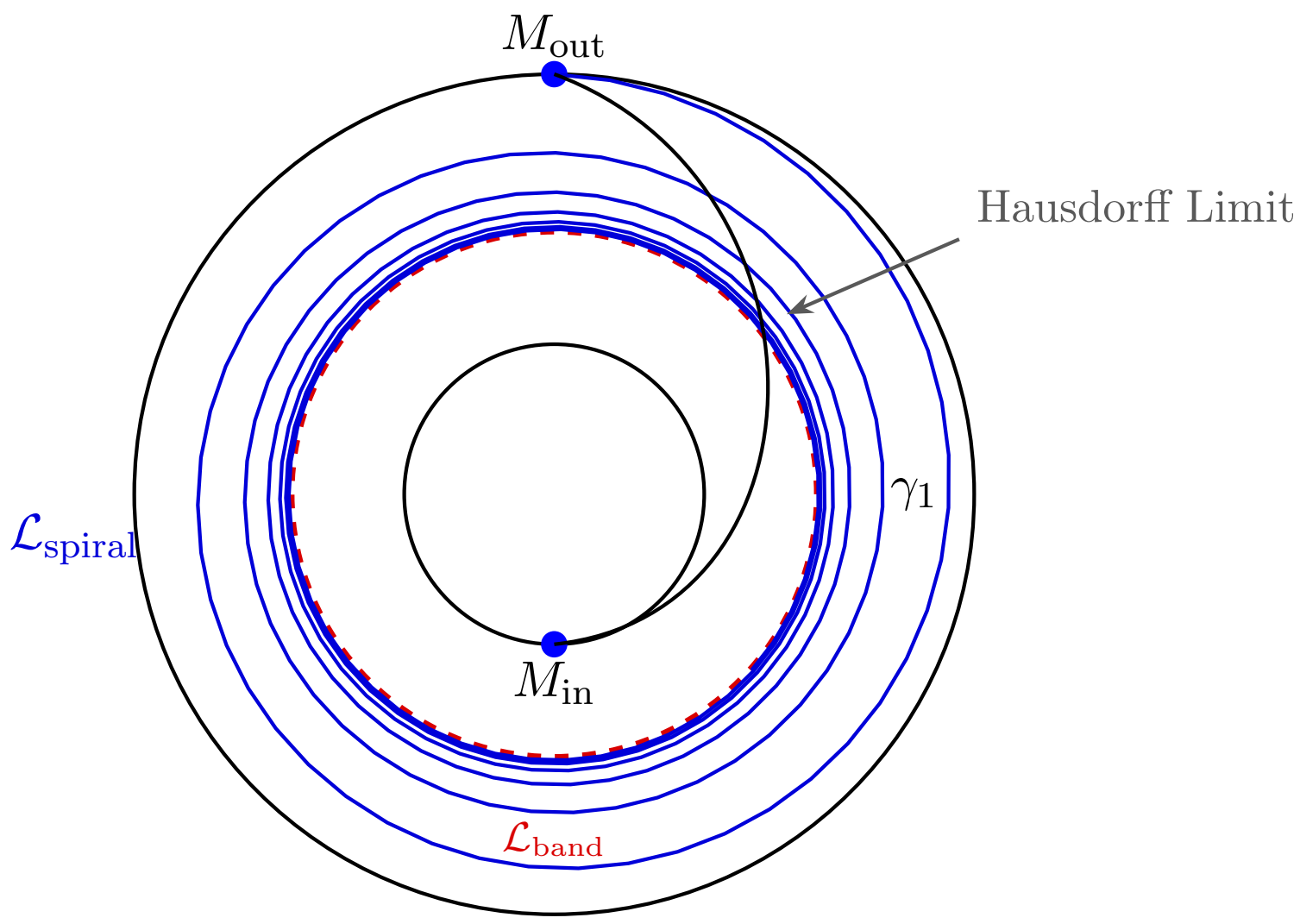
\begin{figure}[htbp]
\centering
\begin{tikzpicture}[scale=1.15]
    \draw[thick] (0,0) circle (2.8cm);
    \draw[thick] (0,0) circle (1.0cm);

    \fill[blue] (0,2.8) circle (2.5pt) node[above, black] {$M_{\mathrm{out}}$};
    \fill[blue] (0,-1.0) circle (2.5pt) node[below, black] {$M_{\mathrm{in}}$};

    \draw[red!85!black, dashed, very thick] (0,0) circle (1.75cm);
    \node[red!85!black, font=\footnotesize] at (0,-2.30) {$\mathcal{L}_{\mathrm{band}}$};

    \draw[blue!85!black, thick, domain=0:40, samples=300, variable=\x]
        plot ({ (1.75 + 1.05*exp(-0.11*\x)) * cos(90 - \x*(180/pi)) },
              { (1.75 + 1.05*exp(-0.11*\x)) * sin(90 - \x*(180/pi)) });

    \node[blue!85!black] at (-3.2, -0.3) {$\mathcal{L}_{\mathrm{spiral}}$};

    \draw[black, thick] (0,2.8) .. controls (1.9,2.1) and (1.9,-0.8) .. (0,-1.0);
    \node[black] at (2.4, 0) {$\gamma_1$};

    \draw[-{Stealth[length=2mm]}, gray!70!black, thick] (2.7,1.7) -- (1.55,1.2);
    \node[gray!70!black, above right] at (2.7,1.7) {\small Hausdorff Limit};
\end{tikzpicture}
\caption[Visualization of Theorem~\ref{thm:theorem_B_sec5}]{Visualization of Theorem~\ref{thm:theorem_B_sec5} on the annular neighborhood of the Kronecker surface model. The finite arcs $\gamma_n$ converge to the non-minimal spiraling geodesic ray $\mathcal{L}_{\mathrm{spiral}}$ (encoding the periodic ray limit module), which accumulates asymptotically onto its minimal component: the closed central geodesic $\mathcal{L}_{\mathrm{band}}$ (encoding the band generic module $G_v$).}
\label{fig:kronecker_torus}
\end{figure}

\begin{example}[The Kronecker Lamination]
\label{ex:kronecker_torus}
Consider the Kronecker algebra \(\Lambda = \K(1 \rightrightarrows 2)\) with arrows \(\alpha, \beta \colon 1 \to 2\) and \(\forb = \emptyset\). Its surface \(\mathcal{S}\) is a cylinder (\(\chi(\mathcal{S}) = 0\)) with two marked segments \(M_{\mathrm{out}}, M_{\mathrm{in}} \in \mathcal{M}\) on the two boundary components. Collapsing these to ideal boundary cusps (while keeping the two complementary unmarked boundary intervals in \(\mathcal{E}\) totally geodesic) yields the two-cusped hyperbolic crown annulus \(\mathcal{S}^\circ\) with area
\[
\operatorname{Area}(\mathcal{S}^\circ) \;=\; \pi\bigl(|\mathcal{M}| - 2\chi(\mathcal{S})\bigr) \;=\; \pi(2 - 0) \;=\; 2\pi \;>\; 0.
\]
The fundamental group \(\pi_1(\mathcal{S}^\circ) \cong \Z\) is generated by a single non-contractible loop, so there is a unique unoriented closed geodesic \(\mathcal{L}_{\mathrm{band}} = \gamma_v\) corresponding to the unique primitive band \(v = \beta^{-1}\alpha\) and the band generic module \(G_v\). The winding string modules \(W_n \coloneqq M((\beta^{-1}\alpha)^n)\) (along with the preprojective string modules \(M(\alpha(\beta^{-1}\alpha)^n)\)) correspond to finite admissible arcs \(\gamma_n\) winding \(n\) times around the cylinder. As \(n \to \infty\), \(W_n\) converges in \(\Zg(\Lambda)\) to the Pr\"ufer-like periodic string ray module, while geometrically \(\gamma_n\) converges to the non-minimal spiraling geodesic ray \(\mathcal{L}_{\mathrm{spiral}}\) accumulating onto \(\mathcal{L}_{\mathrm{band}}\) (Figure~\ref{fig:kronecker_torus}).

The absence of aperiodic laminations on the two-cusped annulus confirms that superdecomposable pure-injective modules cannot exist over the Kronecker algebra, driving the necessity for surfaces with fundamental group containing a free subgroup of rank \(\ge 2\) (such as the four-punctured sphere) to model non-domestic string algebras.
\end{example}

\vspace{-0.9em}

\vspace{-0.3em}
\section{Topological Pruning and Global Classification of the Ziegler Spectrum}
\label{sec:pruning_kg}

In Section~\ref{sec:laminations}, we established that pure-injective limit modules correspond to complete geodesics and geodesic laminations. In this section, we compare the Chabauty topology with the \(T_0\) corridor topology, glue the local geometric fans across all vertices to prove the Global Corridor Homeomorphism on \(C\), and prove Theorem~\ref{thm:theorem_C} (Theorem~\ref{thm:classification}). This establishes the global topological structure of \(\Zg(\Lambda)\) and shows that the Krull-Gabriel dimension \(\KG(\Lambda)\) equals the geometric pruning depth \(\stgeo(\mathcal{S}, \mathcal{M}, \forb)\).

\subsection{The Cantor-Bendixson Derivative and Krull-Gabriel Dimension}
\label{subsec:cb_derivative}

\begin{definition}[Cantor-Bendixson Derivative and Rank]
\label{def:cb_derivative}
Let \(X\) be a topological space. A point \(x \in X\) is \emph{isolated} if \(\{x\}\) is an open subset of \(X\). The \emph{Cantor-Bendixson (CB) derivative} \(D(X)\) is the closed subspace obtained by removing all isolated points of \(X\). For an ordinal \(\alpha\), the transfinite derivatives are \(D^0(X) \coloneqq X\), \(D^{\alpha+1}(X) \coloneqq D(D^\alpha(X))\), and \(D^\beta(X) \coloneqq \bigcap_{\alpha < \beta} D^\alpha(X)\) for a limit ordinal \(\beta\). The \emph{CB-rank} of a point \(x \in X\) is \(\CB(x) = \alpha\) for the unique ordinal \(\alpha\) with \(x \in D^\alpha(X) \setminus D^{\alpha+1}(X)\). The \emph{CB-rank of \(X\)} is \(\CB(X) \coloneqq \sup\{\CB(x) : x \in X\} = \sup\{\alpha : D^\alpha(X) \neq \emptyset\}\), with \(\CB(X) = \infty\) if the \emph{perfect core} \(\mathcal{P}(X) \coloneqq \bigcap_\alpha D^\alpha(X)\) is non-empty.
\end{definition}

For a finite-dimensional string algebra \(\Lambda\), Prest, Schr\"oer, and Ziegler established that the Krull-Gabriel dimension \(\KG(\Lambda)\) coincides with the Cantor-Bendixson rank of the Ziegler spectrum \cite{Ziegler1984, Schroer2000b, Prest2009}:\footnote{We adopt the convention that \(\Zg(\Lambda)\) is the spectrum of right modules, under which \(\KG(\Lambda) = \CB(\Zg(\Lambda))\) holds with \(\KG(\Lambda) = 0\) for representation-finite algebras and \(\KG(\Lambda) = n+1\) for \(n\)-domestic string algebras \cite{Schroer2000b, LPP2018}.}
\[
\KG(\Lambda) \;=\; \CB\bigl(\Zg(\Lambda)\bigr).
\]

\subsection{The Chabauty Space and the Corridor Topology}
\label{subsec:chabauty_corridor}

To translate the CB-rank into geometry, we endow the closure of each geometric fan \(\overline{\mathcal{F}_i} \coloneqq \overline{\mathcal{F}(\gamma_{P(S(i))})}\) with a topology compatible with the Ziegler topology.

\begin{definition}[Chabauty and Corridor Topologies]
\label{def:chabauty_corridor}
Let \(\overline{\mathcal{F}_i}\) be the closure of the geometric fan anchored at \(d_i\). On the space of closed subsets of \(\mathcal{S}\), the classical \emph{Chabauty topology} is generated by the subbasis \(\mathcal{U}_V = \{\mathcal{L} \mid \mathcal{L} \cap V \neq \emptyset\}\) for open \(V \subset \mathcal{S}\), and \(\mathcal{O}_K = \{\mathcal{L} \mid \mathcal{L} \cap K = \emptyset\}\) for compact \(K \subset \mathcal{S}\). The \emph{local corridor topology} \(\tau_{\mathrm{corr}}^i\) on \(\overline{\mathcal{F}_i}\) is the topology generated by the subbasis of geometric test corridors \(\mathcal{C}_{u, \Psi} = \mathcal{U}_{V_{\tau_\phi}} \cap \mathcal{O}_{K_{\tau_\psi}}\) (Definition~\ref{def:test_corridor}), where \(u\) is a finite pointed admissible string crossing \(d_i\) and \(\Psi\) is a finite set of templates strictly extending \(u\).
\end{definition}

\begin{remark}[Reconciling \(T_0\) Corridor and Hausdorff Chabauty Topologies]
\label{rem:t0_hausdorff}
Because a non-minimal spiraling leaf \(\mathcal{L}_{\mathrm{spiral}}\) winds asymptotically onto a minimal closed geodesic \(\gamma_v\), we have \(\gamma_v \subset \overline{\mathcal{L}}_{\mathrm{spiral}}\). In the Hausdorff Chabauty topology, they can be separated by disjoint open neighborhoods. However, in the \(T_0\) Ziegler topology they cannot be separated whenever \(G_v\) lies in the specialization closure of the spiraling ray module, since any basic corridor neighborhood of \(\gamma_v\) is traversed by the asymptotic tail of \(\mathcal{L}_{\mathrm{spiral}}\). The corridor topology \(\tau_{\mathrm{corr}}^i\) is precisely the \(T_0\) specialization coarsening of the Chabauty topology. Crucially, the transfinite Cantor-Bendixson derivatives coincide in both descriptions: a point is isolated in \(\Dgeo^\alpha(\mathcal{C}_{\mathrm{corr}})\) under the corridor topology if and only if it is isolated there under the Chabauty topology, since both notions of isolation are determined by the lamination's dissection crossing sequence.
\end{remark}

\subsection{Geometric Pruning and the Global Topological Classification of \(\Zg(\Lambda)\)}
\label{subsec:global_classification}

\begin{definition}[Geometric Pruning]
\label{def:geometric_pruning}
The \emph{geometric pruning derivative} \(\Dgeo\) is the Cantor-Bendixson derivative acting on \((\overline{\mathcal{F}_i}, \tau_{\mathrm{corr}}^i)\) and globally on the quotient corridor space \(\mathcal{C}_{\mathrm{corr}}\). At each stage \(\alpha\), \(\Dgeo^{\alpha+1}\) removes all pointed geodesics/laminations that are isolated by a test corridor in \(\Dgeo^\alpha\).
\end{definition}

\begin{lemma}[Local Equivalence of Isolation]
\label{lem:equivalence_isolation}
Let \(\Lambda\) be a string algebra, and let \(\overline{H_\Lambda(S(i))}\) be the closure of the string hammock poset in \(\Zg(\Lambda)\). The local geometric realization map \(\Phi_i \colon \overline{H_\Lambda(S(i))} \xrightarrow{\sim} (\overline{\mathcal{F}_i}, \tau_{\mathrm{corr}}^i)\) is a homeomorphism. Consequently, a pure-injective module \(N \in \overline{H_\Lambda(S(i))}\) is isolated at stage \(\alpha\) (i.e., \(\CB(N) = \alpha\)) if and only if \(\Phi_i(N)\) is isolated at stage \(\alpha\) in \(\Dgeo^\alpha(\overline{\mathcal{F}_i})\).
\end{lemma}

\begin{proof}
By Lemmas~\ref{lem:geometric_pp_evaluation} and \ref{lem:corridors_form_basis}, the basic open sets of the subspace Ziegler topology on \(\overline{H_\Lambda(S(i))}\) are given by the pointed string test corridors, which correspond bijectively under \(\Phi_i\) to the basic open test corridors generating \(\tau_{\mathrm{corr}}^i\). Hence \(\Phi_i\) is a homeomorphism and intertwines \(D^\alpha\) with \(\Dgeo^\alpha\) for every ordinal \(\alpha\).
\end{proof}

We now glue the local corridor spaces \(\overline{\mathcal{F}_i}\) across all vertices \(i \in Q_0\) to obtain a global coordinate system on the string closure \(C \coloneqq \bigcup_{i \in Q_0} \overline{H_\Lambda(S(i))} \subseteq \Zg(\Lambda)\).

\begin{lemma}[Global Corridor Homeomorphism on \(C\)]
\label{lem:global_corridor_homeomorphism}
Let \(\Lambda\) be a string algebra modeled on \((\mathcal{S}, \mathcal{M}, \forb)\). Define the anchor-forgetting equivalence relation on \(\bigsqcup_{i \in Q_0} \overline{\mathcal{F}_i}\) by \((\ell_i, p_i) \sim (\ell_j, p_j)\) if and only if they are pointed versions of the same underlying unpointed geodesic/lamination \(\ell\). Then the local homeomorphisms \(\Phi_i\) descend to a global homeomorphism
\[
\Phi \colon C \;\xrightarrow{\;\sim\;} \;\mathcal{C}_{\mathrm{corr}} \;\coloneqq\; \Bigl(\bigsqcup_{i \in Q_0} \overline{\mathcal{F}_i}\Bigr) \Big/\!\sim
\]
from \(C\) (with the subspace Ziegler topology) onto the quotient space equipped with the global corridor quotient topology.
\end{lemma}

\begin{proof}
\emph{Step 1 (Well-Definedness and Bijectivity):} Let \(N \in \overline{H_\Lambda(S(i))} \cap \overline{H_\Lambda(S(j))}\). By Theorem~\ref{thm:theorem_B_sec5}, the isomorphism class of \(N \in C\) is uniquely determined by its underlying admissible word \(\vec{w}_N\) (or lamination \(\mathcal{L}\)). Thus \(\Phi_i(N) = (\ell_N, p_i)\) and \(\Phi_j(N) = (\ell_N, p_j)\) denote the same underlying geodesic equipped with different anchors. Hence \(\Phi_i(N) \sim \Phi_j(N)\), and \(\Phi(N) \coloneqq [\Phi_i(N)]\) is a well-defined bijection.

\emph{Step 2 (Anchor-Shift Compatibility):} Suppose \(N \in \overline{H_\Lambda(S(i))} \cap \overline{H_\Lambda(S(j))}\), so \(\ell_N\) crosses \(d_i\) at \(p_i\) and \(d_j\) at \(p_j\). The segment from \(p_i\) to \(p_j\) traverses a finite bridge string \(t_{ij}\) of length \(d\). Let \(\mathcal{C}_{\phi_j/\psi_j}\) be a basic \(j\)-corridor neighborhood of \(\Phi_j(N)\) determined by a radius-\(R\) numerator template \(u_{\phi_j}\) and denominator extensions of maximal radius \(R+R''\). Set \(R' \coloneqq R + d + R''\), and let \(\mathcal{C}_{\phi_i/\psi_i}\) be the \(i\)-corridor neighborhood at \(p_i\) whose numerator records the radius-\((R+d)\) prefix around \(p_i\) and whose denominator excludes the radius-\(R'\) preimages at \(p_i\) of the extension templates. Any \(N'\) whose \(\Phi_i(N')\) lies in \(\mathcal{C}_{\phi_i/\psi_i}\) agrees with \(\ell_N\) along \(t_{ij}\) and the radius-\(R\) corridor at \(p_j\) while avoiding the extension locus, forcing \(\Phi_j(N') \in \mathcal{C}_{\phi_j/\psi_j}\). Symmetrically, every \(i\)-corridor neighborhood of \(N\) contains a \(j\)-corridor neighborhood. Thus, \(\tau_{\mathrm{corr}}^i\) and \(\tau_{\mathrm{corr}}^j\) induce identical neighborhood filters on the overlap.

\emph{Step 3 (Topological Homeomorphism):} Because \(Q_0\) is finite, \(\{\overline{H_\Lambda(S(i))}\}_{i \in Q_0}\) is a finite closed cover of \(C\). A subset \(U \subseteq C\) is open if and only if \(U \cap \overline{H_\Lambda(S(i))}\) is open for every \(i\). By Lemma~\ref{lem:equivalence_isolation}, each \(\Phi_i\) is a homeomorphism, so \(U\) is open in \(C\) if and only if \(\Phi(U)\) is open in the quotient space \(\mathcal{C}_{\mathrm{corr}}\). Hence \(\Phi\) is a global homeomorphism.
\end{proof}

\begin{theorem}[Theorem C: Global Topological Structure of the Ziegler Spectrum and Krull-Gabriel Dimension]
\label{thm:classification}
Let \(\Lambda = \K\mathcal{Q}/\langle\rho\rangle\) be a string algebra modeled on \((\mathcal{S}, \mathcal{M}, \forb)\), and let \(\mathcal{C}_{\mathrm{corr}}\) be the global corridor quotient space of Lemma~\ref{lem:global_corridor_homeomorphism}. 
\begin{enumerate}
    \item \textbf{Unconditional Homeomorphism:} For every string algebra \(\Lambda\), \(\Phi \colon C \xrightarrow{\sim} \mathcal{C}_{\mathrm{corr}}\) is a global homeomorphism. For every ordinal \(\alpha \ge 2\) (including the perfect core \(\mathcal{P}(\Zg(\Lambda))\)), restriction induces a global homeomorphism \(D^\alpha(\Zg(\Lambda)) = D^\alpha(C) \xrightarrow{\sim} \Dgeo^\alpha(\mathcal{C}_{\mathrm{corr}})\).
    \item \textbf{Closed-Cover Pushout Dichotomy:} The closed embeddings \(C \hookrightarrow \Zg(\Lambda)\) and, for each \(v \in \operatorname{Ba}(\Lambda)\), \(F_v \colon \Zg(\K[t, t^{-1}]) \hookrightarrow \Zg(\Lambda)\) induce a canonical continuous bijection from the topological pushout
    \[
    \Upsilon \colon \mathcal{C}_{\mathrm{corr}} \;\cup_{\{G_v\}}\; \bigsqcup_{v \in \operatorname{Ba}(\Lambda)} F_v\bigl(\Zg(\K[t, t^{-1}])\bigr) \;\longrightarrow\; \Zg(\Lambda)
    \]
    that restricts to a homeomorphism on \(C\) and on each band tube \(F_v\). The map \(\Upsilon\) is a global homeomorphism if and only if \(\Lambda\) is domestic. If \(\Lambda\) is non-domestic, the pushout topology is strictly finer than the Ziegler topology.
\end{enumerate}
In either case, the Krull-Gabriel dimension operates on \(\mathcal{C}_{\mathrm{corr}}\) via geometric pruning, yielding
\[
\KG(\Lambda) \;=\; \stgeo(\mathcal{S}, \mathcal{M}, \forb) \;\coloneqq\; \max_{i \in Q_0} \CB\big(\overline{\mathcal{F}(\gamma_{P(S(i))})}\big).
\]
\end{theorem}

\begin{proof}
(1) By Lemma~\ref{lem:global_corridor_homeomorphism}, \(\Phi \colon C \xrightarrow{\sim} \mathcal{C}_{\mathrm{corr}}\) is an unconditional global homeomorphism. To prove \(D^\alpha(\Zg(\Lambda)) = D^\alpha(C)\) for all \(\alpha \ge 2\), note that the isolated points of \(\Zg(\Lambda)\) are the finite-dimensional string and band modules \cite{PrestPuninski2002}, while the isolated points of \(C\) are the finite-dimensional string modules. Hence \(D^1(\Zg(\Lambda)) \cap C = D^1(C)\), and by the Prest--Puninski classification, the complement \(E \coloneqq D^1(\Zg(\Lambda)) \setminus D^1(C)\) consists precisely of the band Pr\"ufer and adic modules. Each point of \(E\) has CB-rank \(1\) and is isolated in \(D^1(\Zg(\Lambda))\), so no point of \(E\) survives into \(D^2(\Zg(\Lambda))\). Furthermore, since the only accumulation points of \(E\) in \(C\) are the band generic modules \(G_v\), which are already accumulation points of periodic string rays in \(D^1(C)\) (so \(G_v \in D^2(C)\) is non-isolated in \(D^1(C)\)), no point \(x \in D^1(C)\) that is isolated in \(D^1(C)\) lies in the closure \(\overline{E}^{\Zg(\Lambda)}\). Consequently, a point \(x \in D^1(C)\) is isolated in \(D^1(\Zg(\Lambda)) = D^1(C) \cup E\) if and only if it is isolated in \(D^1(C)\). Therefore \(D^2(\Zg(\Lambda)) = D(D^1(\Zg(\Lambda))) = D(D^1(C)) = D^2(C)\). Since all subsequent derivatives take place entirely inside \(C\), transfinite induction on \(\alpha \ge 2\) yields \(D^\alpha(\Zg(\Lambda)) = D^\alpha(C) \cong \Dgeo^\alpha(\mathcal{C}_{\mathrm{corr}})\).

(2) By \cite[\S 3--\S 4]{PrestPuninski2002}, \(\Zg(\Lambda) = C \cup \bigcup_v F_v(\Zg(\K[t, t^{-1}]))\) with \(C \cap F_v(\Zg) = \{G_v\}\) and \(F_v(\Zg) \cap F_{v'}(\Zg) = \emptyset\) for distinct bands. Since \(C\) and each band tube are closed subspaces of \(\Zg(\Lambda)\) \cite[\S 5.5]{Prest2009}, the universal property of the topological pushout yields a continuous bijection \(\Upsilon\) that is a homeomorphism on each closed piece.

If \(\Lambda\) is domestic, \(\operatorname{Ba}(\Lambda)\) is finite, making \(\{C\} \cup \{F_v(\Zg)\}_{v \in \operatorname{Ba}(\Lambda)}\) a finite closed cover of \(\Zg(\Lambda)\). Because a subset of a space with a finite closed cover is closed if and only if its intersection with each cover member is closed, \(\Upsilon\) is closed and thus a global homeomorphism.

Conversely, suppose \(\Lambda\) is non-domestic. There exist two communicating primitive bands \(v_1, v_2\) with bridge strings \(u_{12}, u_{21}\). Fix an integer \(r \ge \lceil 2|v_1|/|v_2| \rceil + 2\) (so that \(r|v_2| \ge \max\{|v_1| + |v_2| - \gcd(|v_1|, |v_2|),\; 2|v_1|\}\)), define \(v_m \coloneqq v_1^m u_{12} v_2^r u_{21}\) (\(m \ge 1\)), and consider the subset of finite band modules \(Z \coloneqq \{M(v_m, 1, 1) \mid m \ge 1\} \subset \Zg(\Lambda)\). We have \(Z \cap C = \emptyset\) (closed in \(C\)). For each \(v \in \operatorname{Ba}(\Lambda)\), \(Z \cap F_v(\Zg)\) is either \(\emptyset\) or the singleton \(\{M(v_m, 1, 1)\}\), which is a closed point in \(F_{v_m}(\Zg)\) because finite-dimensional modules are endofinite. Thus, \(\Upsilon^{-1}(Z)\) is closed in the pushout topology. 

However, \(Z\) is not closed in \(\Zg(\Lambda)\). Fix an arbitrary basic Ziegler open set \((\phi/\psi)\) containing \(G_{v_1} \in C\). Writing \(\phi = \sum_j \phi_{u_j}\) and \(\psi = \sum_k \phi_{w_k}\) as finite sums of pointed string formulas, the condition \(G_{v_1} \in (\phi/\psi)\) implies there exists a witness evaluation of some template \(u_j\) on \(G_{v_1}\) that is not covered by any extension template \(w_k\). Let \(L \coloneqq \max_j |u_j| + \max_k |w_k| + 1\) be the maximal radius of these formulas. 

Because \(v_1\) and \(v_2\) are distinct primitive bands and \(r\) satisfies \(r|v_2| \ge \max\{|v_1| + |v_2| - \gcd(|v_1|, |v_2|),\; 2|v_1|\}\), the concatenated string \(v_m = v_1^m u_{12} v_2^r u_{21}\) is primitive for all sufficiently large \(m\). Suppose for contradiction that \(v_m = w^k\) for some \(k \ge 2\).

\emph{Claim 1: \(v_m = v_1^{jk}\) for some \(j \ge 1\).} Since \(v_m = w^k\) is \(w\)-periodic and \(v_1^m\) is a prefix, \(v_1^m\) inherits the period \(|w|\). Since \(|v_1^m| \ge |v_1| + |w| - \gcd(|v_1|, |w|)\) for \(m\) large, the Fine--Wilf theorem forces \(\gcd(|v_1|, |w|)\) to be a period of \(v_1^m\). By primitivity of \(v_1\), this \(\gcd\) equals \(|v_1|\), so \(|v_1|\) divides \(|w|\). Writing \(|w| = j|v_1|\), we obtain \(w = v_1^j\) (as \(w\) is the length-\(|w|\) prefix of \(v_m\)) and hence \(v_m = v_1^{jk}\).  

\emph{Claim 2: The assumption \(v_m = v_1^{jk}\) contradicts the distinctness of \(v_1\) and \(v_2\).} By Claim~1, \(v_2^r\) is a subword of \(v_1^{jk}\) starting at position \(s = m|v_1| + |u_{12}|\) and of length \(r|v_2|\).   

\emph{(a) Periodicity of \(v_2^r\).} The word \(v_1^{jk}\) is \(|v_1|\)-periodic, so its subword \(v_2^r\) inherits period \(|v_1|\). The word \(v_2^r\) also has period \(|v_2|\) trivially. By the Fine--Wilf theorem, \(\gcd(|v_1|, |v_2|)\) is a period of \(v_2^r\). Since \(v_2\) is primitive, this forces \(\gcd(|v_1|, |v_2|) = |v_2|\), so \(|v_2|\) divides \(|v_1|\); write \(|v_1| = K|v_2|\).  

\emph{(b) Local period propagation.} Since \(v_2^r\) has period \(|v_2|\) and is a subword of \(v_1^{jk}\) at position \(s\), we have \(v_1^{jk}[s + i] = v_1^{jk}[s + i + |v_2|]\) for all \(0 \le i < (r-1)|v_2|\). Thus \(v_1^{jk}\) has period \(|v_2|\) on the interval \([s, s + r|v_2|)\), which by our choice of \(r\) has length \(\ge 2|v_1|\) and so contains at least one full \(|v_1|\)-block of \(v_1^{jk}\).  

\emph{(c) A cyclic shift of \(v_1\) becomes a proper power.} Let \(w_1\) be any length-\(|v_1|\) factor of \(v_1^{jk}\) contained in this interval. Since \(v_1^{jk}\) is \(|v_1|\)-periodic, \(w_1\) is a cyclic shift of \(v_1\). Since \(v_1^{jk}\) has period \(|v_2|\) on the interval, \(w_1\) has period \(|v_2|\), so \(w_1 = u^K\) for the length-\(|v_2|\) prefix \(u\). If \(K > 1\), this expresses the cyclic shift \(w_1\) of the primitive band \(v_1\) as a proper power, contradicting primitivity of cyclic shifts of \(v_1\). Hence \(K = 1\), i.e., \(|v_1| = |v_2|\).   

\emph{(d) Contradiction.} Since \(v_2\) is a length-\(|v_1|\) factor of \(v_1^{jk}\), the equality \(|v_1| = |v_2|\) forces \(v_2\) to be a cyclic shift of \(v_1\), contradicting that \(v_1, v_2\) are distinct primitive bands in \(\operatorname{Ba}(\Lambda)\). Hence no such \(k \ge 2\) exists, and \(v_m\) is primitive.

Fix the holonomy scalar \(\lambda = 1\) without loss of generality. For \(m \ge 2L + 1 + \lceil 2L/|v_1| \rceil\), the string \(v_m\) contains an interior copy of \(v_1^{2L+1}\) of length \((2L+1)|v_1|\) syllables, positioned at distance at least \(L\) syllables from both junctions with the bridge strings \(u_{12}\) and \(u_{21}\). The pp-type of an element is determined by the set of pp-formulas it satisfies, and a pp-formula of template-radius \(\le L\) depends only on the element's local neighborhood of radius \(L\). Since this interior \(v_1^{2L+1}\)-segment in \(M(v_m, 1, 1)\) matches the bi-infinite periodic neighborhood in \(G_{v_1}\) strictly up to radius \(L\) in both directions, the anchor at the center of this segment has the exact same pp-type as the witness occurrence in \(G_{v_1}\) for all templates of radius \(\le L\). Consequently, the \(u_j\)-template is satisfied at this interior anchor, while no \(w_k\)-template covers it. Thus \(\phi(M(v_m, 1, 1)) \not\subseteq \psi(M(v_m, 1, 1))\), implying \(M(v_m, 1, 1) \in (\phi/\psi)\). This proves \(G_{v_1} \in \overline{Z}^{\Zg(\Lambda)}\). Since \(G_{v_1} \in C\) and \(Z \cap C = \emptyset\), \(G_{v_1} \notin Z\), proving that the pushout topology is strictly finer than the Ziegler topology when \(\operatorname{Ba}(\Lambda)\) is infinite.

Finally, since \(\KG(\Lambda) = \CB(\Zg(\Lambda)) = \CB(C) = \max_{i \in Q_0} \CB(\overline{H_\Lambda(S(i))})\), Lemma~\ref{lem:equivalence_isolation} yields \(\KG(\Lambda) = \stgeo(\mathcal{S}, \mathcal{M}, \forb)\).
\end{proof}

\subsection{Example: Geometric Pruning of the Kronecker Algebra}
\label{subsec:pruning_kronecker}

\begin{example}[Three-Step Geometric Pruning of the Kronecker Algebra]
\label{ex:pruning_kronecker}
Let \(\Lambda = \K(1 \rightrightarrows 2)\) with its two-cusped annulus surface model (Example~\ref{ex:kronecker_torus}). The space \(\overline{\mathcal{F}}\) contains the finite winding arcs \(\gamma_n\), the spiraling rays \(\mathcal{L}_{\mathrm{spiral}}\), and the central closed geodesic \(\mathcal{L}_{\mathrm{band}}\):
\begin{enumerate}
    \item \textbf{Level 0 (\(\Dgeo^1\)):} Every finite arc \(\gamma_n\) has both endpoints on the boundary marked segments. A test corridor fixing the endpoints isolates \(\gamma_n\). Thus, every finite arc has CB-rank 0 and is pruned by \(\Dgeo^1\).
    \item \textbf{Level 1 (\(\Dgeo^2\)):} In \(\Dgeo^1(\overline{\mathcal{F}})\), the points are the spiraling rays and the central closed geodesic. Because the finite arcs are removed, each spiraling ray is isolated by a corridor near its boundary anchor. Hence, \(\CB(\mathcal{L}_{\mathrm{spiral}}) = 1\) and \(\Dgeo^2\) prunes all spiraling rays.
    \item \textbf{Level 2 (\(\Dgeo^3\)):} After pruning the spiraling rays, \(\Dgeo^2(\overline{\mathcal{F}}) = \{\mathcal{L}_{\mathrm{band}}\}\) consists of a single isolated point. Thus, \(\CB(\mathcal{L}_{\mathrm{band}}) = 2\) and \(\Dgeo^3(\overline{\mathcal{F}}) = \emptyset\).
\end{enumerate}
By Definition~\ref{def:cb_derivative} and Theorem~\ref{thm:classification}, \(\KG(\Lambda) = \stgeo = \sup\{0, 1, 2\} = 2\), recovering the classical Krull-Gabriel dimension purely from the three-step topological pruning of the surface model.
\end{example}

\vspace{-0.4em}
\section{Classification of Representation Type}
\label{sec:classification}

In Section~\ref{sec:pruning_kg}, we proved that the Krull-Gabriel dimension \(\KG(\Lambda)\) equals the topological accumulation depth \(\stgeo(\mathcal{S}, \mathcal{M}, \forb)\) of the corridor space \(\mathcal{C}_{\mathrm{corr}}\). In this section, we use the global corridor homeomorphism (Theorem~\ref{thm:classification}) to translate the representation-theoretic classes of string algebras (\(1\)-domestic, \(n\)-domestic, and non-domestic) into geometric properties of surface geodesics and laminations. We conclude with a geometric perspective on how the forbidden-boundary surface topology aligns with algebraic tameness.

\subsection{\(1\)-Domestic Algebras and Unnested Closed Geodesics}
\label{subsec:1_domestic}

Recall that a representation-infinite string algebra \(\Lambda\) is \emph{domestic} if the set of primitive bands \(\operatorname{Ba}(\Lambda)\) is finite, and \emph{\(1\)-domestic} if \(\Lambda\) is domestic and no admissible bridge string connects two distinct primitive bands \cite{Schroer2000b}. By Schr\"oer \cite{Schroer2000b}, representation-infinite \(1\)-domestic string algebras satisfy \(\KG(\Lambda) = 2\).

\begin{proposition}[Topological Characterization of \(1\)-Domesticity]
\label{prop:1_domestic}
Let \(\Lambda\) be a representation-infinite string algebra modeled on \((\mathcal{S}, \mathcal{M}, \forb)\). Then \(\Lambda\) is \(1\)-domestic (\(\KG(\Lambda) = 2\)) if and only if the set of primitive closed geodesics on \(\mathcal{S}^\circ \setminus \forb\) is a non-empty finite discrete set, and no heteroclinic geodesic spirals asymptotically from one closed geodesic to another.
\end{proposition}

\begin{proof}
By Theorem~\ref{thm:theorem_B_sec5}, primitive bands correspond bijectively to unoriented primitive closed geodesics \(\gamma_v\). Suppose they form a finite set with no heteroclinic spiraling. Then no bridge string connects distinct bands, so \(\Lambda\) is \(1\)-domestic. In \(\mathcal{C}_{\mathrm{corr}}\), every non-minimal infinite geodesic is anchored at an ideal cusp and spirals into a single closed geodesic \(\gamma_{v_j}\). 

The first pruning derivative \(\Dgeo^1\) removes all finite cusp-to-cusp arcs (CB-rank \(0\)). In \(\Dgeo^1(\mathcal{C}_{\mathrm{corr}})\), every pointed spiraling geodesic is isolated by a test corridor fixing its finite prefix from the anchor, so \(\Dgeo^2\) prunes all spiraling leaves (CB-rank \(1\)). The remaining space \(\Dgeo^2(\mathcal{C}_{\mathrm{corr}})\) is a non-empty finite discrete set of closed geodesics (CB-rank \(2\)), so \(\Dgeo^3(\mathcal{C}_{\mathrm{corr}}) = \emptyset\) and \(\KG(\Lambda) = \stgeo = 2\).

Conversely, suppose there exists a heteroclinic bridge \(u_{12}\) from \(\gamma_{v_1}\) to \(\gamma_{v_2}\). For each fixed \(m \ge 1\), as \(k \to \infty\) the finite pointed arcs \(\gamma_{v_1^m u_{12} v_2^k}\) (CB-rank \(0\)) converge in \(\mathcal{C}_{\mathrm{corr}}\) to the periodic spiraling ray \(\ell_m = \gamma_{v_1^m u_{12} v_2^\infty} \in \Dgeo^1(\mathcal{C}_{\mathrm{corr}})\). As \(m \to \infty\), these rays \(\ell_m \in \Dgeo^1(\mathcal{C}_{\mathrm{corr}})\) accumulate onto the heteroclinic geodesic \(\ell_{1,2} = \gamma_{v_1^{-\infty} u_{12} v_2^\infty}\), which therefore lies in \(\Dgeo^2(\mathcal{C}_{\mathrm{corr}})\). In turn, shifting the anchor along \(\ell_{1,2}\) into the \(v_2^\infty\) tail shows that \(\ell_{1,2} \in \Dgeo^2(\mathcal{C}_{\mathrm{corr}})\) accumulates onto the closed geodesic \(\gamma_{v_2}\), forcing \(\gamma_{v_2} \in \Dgeo^3(\mathcal{C}_{\mathrm{corr}}) \neq \emptyset\) and hence \(\KG(\Lambda) = \stgeo \ge 3\); the general \(n\)-domestic case (\(n \ge 2\)) is treated in Proposition~\ref{prop:n_domestic}.
\end{proof}

\subsection{\(n\)-Domestic Algebras and Heteroclinic Spiraling Chains}
\label{subsec:n_domestic}

General domestic string algebras can possess finitely many directed bridges between distinct primitive bands, forming a directed acyclic communication graph. A domestic string algebra \(\Lambda\) is \emph{\(n\)-domestic} (\(n \ge 1\)) if the maximal number of primitive bands along a directed chain of communicating bands is \(n\); in this case, \(\KG(\Lambda) = n + 1\) \cite{Schroer2000b, LPP2018}.

\begin{proposition}[Heteroclinic Spiraling Chains]
\label{prop:n_domestic}
If a string algebra \(\Lambda\) modeled on \((\mathcal{S}, \mathcal{M}, \forb)\) is \(n\)-domestic (\(n \ge 1\)), then the primitive closed geodesics on \(\mathcal{S}^\circ \setminus \forb\) form a finite set whose heteroclinic spiraling bridges form a directed acyclic graph of maximal chain length \(n-1\), and \(\stgeo(\mathcal{S}, \mathcal{M}, \forb) = \KG(\Lambda) = n + 1\).
\end{proposition}

\begin{proof}
By Theorem~\ref{thm:theorem_B_sec5}, an admissible bridge string from \(v_j\) to \(v_k\) corresponds to a heteroclinic geodesic spiraling out of \(\gamma_{v_j}\) and into \(\gamma_{v_k}\). Because \(\Lambda\) is domestic, no two bands communicate in both directions (otherwise \(v_1^m u_{12} v_2 u_{21}\) would generate infinitely many primitive bands). Thus the communication graph on \(\{\gamma_v\}\) is a finite directed acyclic graph. By definition of \(n\)-domesticity, the maximal number of vertices in a directed path is \(n\), which has \(n-1\) directed bridge edges. By \cite{Schroer2000b, LPP2018} and Theorem~\ref{thm:classification}, \(\stgeo(\mathcal{S}, \mathcal{M}, \forb) = \KG(\Lambda) = n + 1\).
\end{proof}

\subsection{Non-Domestic Algebras and Perfect Geodesic Cores}
\label{subsec:non_domestic}

A string algebra \(\Lambda\) is \emph{non-domestic} if \(\operatorname{Ba}(\Lambda)\) is infinite, which occurs if and only if there exist two distinct primitive bands that communicate in both directions \cite{Schroer2000b, PrestPuninski2002}.

\begin{theorem}[Non-Domesticity and Perfect Geodesic Cores]
\label{thm:exponential_growth}
Let \(\Lambda\) be a representation-infinite string algebra modeled on \((\mathcal{S}, \mathcal{M}, \forb)\). The following are equivalent:
\begin{enumerate}
    \item \(\Lambda\) is non-domestic.
    \item There exist two primitive closed geodesics \(\gamma_{v_1}, \gamma_{v_2}\) on \(\mathcal{S}^\circ \setminus \forb\) that communicate in both directions.
    \item The Krull-Gabriel dimension is infinite: \(\KG(\Lambda) = \stgeo(\mathcal{S}, \mathcal{M}, \forb) = \infty\).
    \item The geometric perfect core \(\mathcal{P}_{\mathrm{geo}} \coloneqq \bigcap_{\alpha} \Dgeo^\alpha(\mathcal{C}_{\mathrm{corr}})\) is non-empty and contains an uncountable compact perfect subset \(W = \{w_1, w_2\}^{\Z}\) (with \(w_1 = v_1 u_{12} v_2 u_{21}\) and \(w_2 = v_1^2 u_{12} v_2 u_{21}\)); every point of \(W\) corresponds to a pure-injective module at CB-rank \(\infty\) in \(\Zg(\Lambda)\), and every minimal aperiodic subshift \(\Omega \subset W\) (such as the uncountable family of Sturmian subshifts \(\Omega_\alpha \subset W\)) is a compact perfect Cantor set inside \(\mathcal{P}_{\mathrm{geo}}\) consisting entirely of minimal recurrent aperiodic geodesics.
\end{enumerate}
\end{theorem}

\begin{proof}
(1 \(\iff\) 2) \(\Lambda\) is non-domestic if and only if \(\operatorname{Ba}(\Lambda)\) is infinite, which holds if and only if there exist two primitive bands communicating in both directions \cite[\S 3]{PrestPuninski2002}. By Theorem~\ref{thm:theorem_B_sec5}, this is equivalent to \(\gamma_{v_1}\) and \(\gamma_{v_2}\) communicating in both directions on \(\mathcal{S}^\circ \setminus \forb\).

(2 \(\implies\) 4) Suppose \(\gamma_{v_1}\) and \(\gamma_{v_2}\) communicate in both directions via bridge strings \(u_{12}, u_{21}\) crossing \(d_i\). Set \(w_1 = v_1 u_{12} v_2 u_{21}\) and \(w_2 = v_1^2 u_{12} v_2 u_{21}\). The subshift of bi-infinite words \(W \coloneqq \{w_1, w_2\}^{\Z}\) is an uncountable Cantor space. By Theorem~\ref{thm:theorem_B_sec5}, each \(w \in W\) defines a pointed complete geodesic \((\ell_w, p_0) \in \overline{\mathcal{F}_i}\). Because any basic test corridor constrains only a finite prefix around \(p_0\), every neighborhood of \((\ell_w, p_0)\) contains uncountably many other elements of \(W\) sharing that finite prefix. Thus \(W\) is a non-empty perfect subset of \(\mathcal{C}_{\mathrm{corr}}\) devoid of isolated points, so \(W \subseteq \Dgeo^1(\mathcal{C}_{\mathrm{corr}})\). By transfinite induction, \(W \subseteq \bigcap_\alpha \Dgeo^\alpha(\mathcal{C}_{\mathrm{corr}}) = \mathcal{P}_{\mathrm{geo}}\), so \(\mathcal{P}_{\mathrm{geo}}\) contains the uncountable perfect subset \(W\). The periodic bi-infinite words in \(W\) (such as \(w_1^{\Z}\) and \(w_2^{\Z}\)) correspond to the periodic bi-infinite string pure-injective modules \(\PE(M(w_1^{\Z}))\) and \(\PE(M(w_2^{\Z}))\), which lie in \(\mathcal{C}_{\mathrm{corr}}\) at CB-rank \(\infty\) and are distinct from both the band Pr\"ufer/adic modules and the band generic modules \(G_v\). Furthermore, every minimal aperiodic subshift \(\Omega \subset W\) (such as the uncountably many Sturmian subshifts \(\Omega_\alpha \subset \{w_1, w_2\}^{\Z}\) for irrational \(\alpha \in (0,1) \setminus \Q\)) is a compact perfect Cantor set inside \(\mathcal{P}_{\mathrm{geo}}\) consisting entirely of minimal recurrent aperiodic geodesics at CB-rank \(\infty\).

(4 \(\implies\) 3) If \(\mathcal{P}_{\mathrm{geo}} \neq \emptyset\), the transfinite pruning process never terminates, so \(\stgeo(\mathcal{S}, \mathcal{M}, \forb) = \infty\). By Theorem~\ref{thm:classification}, \(\KG(\Lambda) = \stgeo = \infty\).

(3 \(\implies\) 1) If \(\Lambda\) were domestic, \(\stgeo = \KG(\Lambda) = n + 1 < \infty\) by Proposition~\ref{prop:n_domestic}, contradicting \(\stgeo = \infty\).
\end{proof}

\subsection{A Geometric Perspective on Tameness}
\label{subsec:tameness_heuristic}

\begin{remark}[Topological Obstruction Heuristic for Tameness]
\label{rem:tameness_heuristic}
While string and special biserial algebras are unconditionally tame by \cite{WaldWaschbusch1985, ButlerRingel1987}, our dictionary provides a geometric reason. Continuous 2-parameter families of finite-dimensional indecomposables at a fixed dimension vector would require continuous 2-parameter local systems on a single homotopy class. However, the surface model restricts finite-dimensional indecomposables to 1-dimensional immersed arcs \(\gamma_w \simeq [0,1]\) (whose trivial cohomology \(H^1([0,1]; \K^\times) = \{1\}\) forces rigidity) or primitive closed curves \(\gamma_v \simeq S^1\) (whose rank-\(1\) local systems are parameterized by \(H^1(S^1; \K^\times) \cong \K^\times\)). This 1-dimensional cohomological bound geometrically realizes the definition of tameness, leaving a purely topological proof of Drozd's tame-wild dichotomy for wild algebras as future work.
\end{remark}

\vspace{-0.5em}
\section{Superdecomposable Pure-Injectives and the Perfect Core}
\label{sec:superdecomposables}

In this section, we study the pure-injective modules governed by the perfect core
\[
\mathcal{P} \;=\; \bigcap_\alpha D^\alpha\bigl(\Zg(\Lambda)\bigr).
\]
We first construct the pointed evaluation map \(\Psi\) and establish the homeomorphism between the geometric perfect core \(\mathcal{P}_{\mathrm{geo}}\) and the algebraic perfect core \(\mathcal{P}\). We then use a prefix-free dyadic block construction in \(\mathcal{P}_{\mathrm{geo}}\) and a self-contained monomial-basis proof of free distributive generation inside the localized pp-lattice \(\mathcal{L}_i(\Lambda)\) to verify the hypothesis of Puninski's cardinality-free two-chain criterion \cite[Theorem~2.3]{Puninski2008}, resolving the open problem on the existence of superdecomposable pure-injective modules over uncountable base fields for all non-domestic string and special biserial algebras.

\subsection{Superdecomposable Modules and the Perfect Core}
\label{subsec:superdecomp_core}

\begin{definition}[Superdecomposable Modules, \cite{Ziegler1984, Krause1998, Prest2009}]
\label{def:superdecomposable}
A non-zero pure-injective module \(M \in \Mod{\Lambda}\) is \emph{superdecomposable} if it contains no non-zero indecomposable direct summands.
\end{definition}

Every point of the Ziegler spectrum \(\Zg(\Lambda)\) is by definition an indecomposable pure-injective module, so no single point of \(\Zg(\Lambda)\) is itself superdecomposable. Rather, a superdecomposable pure-injective module is a large pure-injective module whose Ziegler support is a non-empty closed subset of the perfect core \(\mathcal{P}(\Zg(\Lambda)) = \bigcap_\alpha D^\alpha(\Zg(\Lambda))\) devoid of isolated points. By Ziegler \cite[Theorem~7.8]{Ziegler1984}, if a ring \(R\) admits a superdecomposable pure-injective module, then the lattice of pp-formulas \(\mathcal{L}(R)\) has infinite width. The converse holds unconditionally when \(R\) is countable \cite[Theorem~7.8]{Ziegler1984}, and holds over arbitrary fields whenever \(\mathcal{L}(R)\) satisfies Puninski's free distributive two-chain criterion \cite[Theorem~2.3]{Puninski2008}. For string algebras, infinite width is equivalent to \(\KG(\Lambda) = \infty\) by Prest--Puninski \cite{PrestPuninski2002} and Schr\"oer \cite{Schroer2000b}.

\subsection{The Pointed Evaluation Map \(\Psi\) and Binary Corridor Trees}
\label{subsec:evaluation_map_psi}

For each vertex \(i \in Q_0\), let \(\mathcal{P}_{\mathrm{geo}}(i) \coloneqq \bigcap_\alpha \Dgeo^\alpha(\overline{\mathcal{F}_i})\) denote the geometric perfect core of the fan anchored at \(d_i\).

\begin{construction}[The Pointed Evaluation Map \(\Psi\)]
\label{const:evaluation_map}
For any pointed complete aperiodic geodesic \((\ell, p_0) \in \mathcal{P}_{\mathrm{geo}}(i)\) (in particular, any pointed leaf of a minimal arational geodesic lamination \(\mathcal{L} \subset \mathcal{S}^\circ \setminus \forb\)), the sequence of transversal crossings of \(\ell\) with the cellular dissection \(\Delta\) centered at \(p_0 \in d_i\) records a bi-infinite (or one-sided infinite) admissible string \(\overleftrightarrow{w}_{(\ell, p_0)}\) avoiding \(\rho \cup \rho^{-1}\). We define the \emph{pointed evaluation map}
\[
\Psi_i \colon \mathcal{P}_{\mathrm{geo}}(i) \;\longrightarrow\; \mathcal{P}\bigl(\Zg(\Lambda)\bigr) \cap \overline{H_\Lambda(S(i))}, \qquad \Psi_i(\ell, p_0) \;\coloneqq\; \PE\Bigl(M\bigl(\overleftrightarrow{w}_{(\ell, p_0)}\bigr)\Bigr),
\]
as the pure-injective envelope of the infinite-dimensional string module \(M(\overleftrightarrow{w}_{(\ell, p_0)})\) anchored at the canonical basis vector \(v_0 \in M(\overleftrightarrow{w}_{(\ell, p_0)})e_i\) corresponding to \(p_0 \in d_i\). Because all dense leaves of a minimal arational lamination \(\mathcal{L}\) share the exact same finite subwords (and hence generate the same basic corridor open sets), the isomorphism class of \(\Psi_i(\ell, p_0)\) in \(\Zg(\Lambda)\) depends only on \(\mathcal{L}\) and not on the specific leaf or basepoint, descending to the inverse \(\Psi = \Phi^{-1} \colon \mathcal{P}_{\mathrm{geo}}(\mathcal{C}_{\mathrm{corr}}) \xrightarrow{\sim} \mathcal{P}(\Zg(\Lambda))\) of the global corridor homeomorphism.
\end{construction}

\subsection{Homeomorphism of the Perfect Cores and Geometric Obstruction}
\label{subsec:perfect_core_homeomorphism}

\begin{theorem}[Homeomorphism of the Perfect Cores]
\label{thm:perfect_core_homeomorphism}
Let \(\Lambda\) be a non-domestic string algebra modeled on \((\mathcal{S}, \mathcal{M}, \forb)\). The global evaluation map restricts to a canonical homeomorphism
\[
\Psi = \Phi^{-1} \colon \mathcal{P}_{\mathrm{geo}}(\mathcal{C}_{\mathrm{corr}}) \;\xrightarrow{\;\sim\;}\; \mathcal{P}\bigl(\Zg(\Lambda)\bigr)
\]
between the geometric perfect core of admissible geodesics/laminations in \(\mathcal{C}_{\mathrm{corr}}\) and the algebraic perfect core of the Ziegler spectrum \(\Zg(\Lambda)\).
\end{theorem}

\begin{proof}
By Theorem~\ref{thm:classification}(1), for every ordinal \(\alpha \ge 2\), the global corridor map \(\Phi \colon C \xrightarrow{\sim} \mathcal{C}_{\mathrm{corr}}\) restricts to a topological homeomorphism \(D^\alpha(\Zg(\Lambda)) = D^\alpha(C) \xrightarrow{\sim} \Dgeo^\alpha(\mathcal{C}_{\mathrm{corr}})\). Taking the intersection over all ordinals \(\alpha\) yields \(\mathcal{P}(\Zg(\Lambda)) = \bigcap_\alpha D^\alpha(C) \xrightarrow{\sim} \bigcap_\alpha \Dgeo^\alpha(\mathcal{C}_{\mathrm{corr}}) = \mathcal{P}_{\mathrm{geo}}(\mathcal{C}_{\mathrm{corr}})\), with inverse \(\Psi = \Phi^{-1}\).
\end{proof}

\begin{proposition}[Geometric Obstruction to Finite Corridor Isolation]
\label{prop:geometric_obstruction}
Let \(\Lambda\) be a non-domestic string algebra, and let \(N = \Psi(\mathcal{L}) \in \mathcal{P}(\Zg(\Lambda))\) be the indecomposable pure-injective module associated with a minimal arational lamination \(\mathcal{L} \subset \mathcal{S}^\circ \setminus \forb\) (or a Cantor lamination subshift \(W\)). Then no basic open set \((\phi/\psi)\) defined by finite pp-formulas can isolate \(N\) inside \(\mathcal{P}(\Zg(\Lambda))\).
\end{proposition}

\begin{proof}
By Lemma~\ref{lem:corridors_form_basis}, any basic Ziegler open neighborhood of \(N\) in \(C\) contains a pointed string test corridor \(\mathcal{C}_{u, \Psi}\) of some finite template radius \(R < \infty\). Because \(\mathcal{L}\) is a minimal arational lamination (or a Cantor lamination subshift \(W\)), its intersection with a transversal arc is a perfect Cantor set devoid of isolated points, so infinitely many distinct leaves and continuations share the radius-\(R\) prefix \(u\) while avoiding \(\Psi\) and diverging beyond radius \(R\). Thus \((\phi/\psi) \cap \mathcal{P}(\Zg(\Lambda))\) contains uncountably many points of \(\mathcal{P}(\Zg(\Lambda))\) and cannot isolate \(N\).
\end{proof}

\subsection{Resolution of Puninski's Open Question over Uncountable Fields}
\label{subsec:uncountable_superdecomp}

Puninski \cite{Puninski2004} proved that the pp-lattice of any non-domestic string algebra has infinite width. However, extracting a superdecomposable pure-injective module from infinite width via Ziegler's theorem \cite[Theorem~7.8(2)]{Ziegler1984} requires the base field \(\K\) to be countable, leaving the case of uncountable fields \(\K\) explicitly open \cite[Abstract and \S 4]{Puninski2008} (see also Prest's MathSciNet review MR2434248 of \cite{Puninski2008}). We now resolve this problem over fields \(\K\) of arbitrary cardinality.

\begin{lemma}[Density of the Transverse Leaf Order]
\label{lem:leaf_density}
Let \(\mathcal{L} \subset \mathcal{S}^\circ \setminus \forb\) be a minimal arational geodesic lamination, and let \(\tau \subset \mathcal{S}^\circ\) be a complete transversal arc intersecting \(\mathcal{L}\). Then:
\begin{enumerate}
    \item \(\mathcal{L} \cap \tau\) is a non-empty compact perfect subset of \(\tau\).
    \item The linear order of \(\mathcal{L} \cap \tau\) along \(\tau\) (after removing the two extreme boundary endpoints of \(\mathcal{L} \cap \tau\), or on any interior open subinterval) is a dense linear order.
\end{enumerate}
\end{lemma}

\begin{proof}
(1) By minimality of \(\mathcal{L}\), every leaf is dense in \(\mathcal{L}\); if \(\mathcal{L} \cap \tau\) contained an isolated point \(p\), the leaf through \(p\) would be isolated in \(\mathcal{L}\) near \(p\), contradicting density. (2) For any open subarc \(\tau_0 \subset \tau\) intersecting \(\mathcal{L}\), \(\mathcal{L} \cap \tau_0\) is infinite since \(\mathcal{L} \cap \tau\) is perfect; picking three points \(p_1 < p_2 < p_3\) in \(\mathcal{L} \cap \tau_0\) shows that the middle point \(p_2\) lies strictly between two leaves of \(\mathcal{L}\) in \(\tau_0\), establishing density. (Geometrically, by the Casson--Bleiler structure theorem \cite[\S I.4.2]{CassonBleiler1988}, see also \cite[\S 4.1]{MjRafi2018}, on a finite-type surface the complementary regions of a minimal arational lamination are finitely many ideal polygons, so only the finitely many boundary leaves of those ideal polygons have immediate successors in \(\mathcal{L} \cap \tau\), while the interior points form a dense order isomorphic to \((\mathbb{Q}, <)\).)
\end{proof}

\noindent\textbf{Convention on the Prest--Puninski String Order \(\le_{\mathrm{str}}\).}
Following Prest and Puninski \cite[\S 2]{PrestPuninski2002}, on the set of one-sided strings starting at vertex \(i\) in a fixed walk direction, \(u <_{\mathrm{str}} v\) if and only if either \(u = v \alpha w'\) for a direct arrow \(\alpha\) (direct extensions decrease) or \(v = u \beta^{-1} w'\) for an inverse arrow \(\beta^{-1}\) (inverse extensions increase), or at the first syllable where \(u\) and \(v\) diverge, \(u\) takes a direct syllable (or terminates) while \(v\) takes an inverse syllable (or terminates). Under this convention, the mapping \(u \mapsto \phi_u(x)\) is strictly order-preserving on one-sided string formulas:
\[
u \le_{\mathrm{str}} v \quad \Longleftrightarrow \quad \phi_u(M) \subseteq \phi_v(M) \;\text{ for all } M \in \Mod{\Lambda},
\]
and the maximal inverse extension \(u^\uparrow\) is the \(\le_{\mathrm{str}}\)-maximum element having prefix \(u\).

\begin{lemma}[Monomial-Basis Proof of Free Distributive Generation]
\label{lem:free_generation_monomial}
Let \(\Lambda\) be a non-domestic string algebra over an algebraically closed field \(\K\) of arbitrary cardinality. The localized pp-lattice \(\mathcal{L}_i(\Lambda)\) contains a sublattice that is \emph{freely generated} as a distributive lattice \(\mathcal{C}^- \otimes \mathcal{C}^+\) by two chains \(\mathcal{C}^-\) and \(\mathcal{C}^+\) of pointed string pp-formulas, each order-isomorphic to \((\mathbb{Q}, <)\); restricting to closed bounded subchains \(\widehat{\mathcal{C}}^\pm \subset \mathcal{C}^\pm\) (each containing a densely ordered open interior isomorphic to \((\mathbb{Q}, <)\)) yields a free distributive product \(\widehat{\mathcal{C}}^- \otimes \widehat{\mathcal{C}}^+\) inside an explicit closed interval \([a_0 \wedge b_0,\, a_1 + b_1]\) of \(\mathcal{L}_i(\Lambda)\).
\end{lemma}

\begin{proof}
Let \(v_1, v_2\) be two communicating primitive bands at vertex \(i\) with bridge strings \(u_{12}, u_{21}\), and set \(w_1 = v_1 u_{12} v_2 u_{21}\) and \(w_2 = v_1^2 u_{12} v_2 u_{21}\). Then \(U_1 \coloneqq w_1 w_2\) and \(U_2 \coloneqq w_2 w_1\) are distinct admissible closed blocks of equal syllable length \(L = |w_1| + |w_2|\) starting with the exact same initial syllable in \(Q_{1,\mathrm{walk}}^+(i)\) and ending with the exact same terminal syllable in \(Q_{1,\mathrm{walk}}^-(i)\) across \(d_i\), and diverging internally before completion. Relabeling if necessary, assume \(U_1 <_{\mathrm{str}} U_2\). Define three blocks of equal length \(2L\):
\[
B_0 \;\coloneqq\; U_1 U_1, \qquad B_1 \;\coloneqq\; U_1 U_2, \qquad B_2 \;\coloneqq\; U_2 U_2,
\]
which satisfy \(B_0 <_{\mathrm{str}} B_1 <_{\mathrm{str}} B_2\). Define the prefix-free dyadic block family on the forward arm \(Q_{1,\mathrm{walk}}^+(i)\):
\[
\mathcal{D}^+ \;\coloneqq\; \Bigl\{\, B_{d_1} B_{d_2} \dots B_{d_m} B_1 \;\Big|\; m \ge 0,\; d_1, \dots, d_m \in \{0, 2\} \,\Bigr\} \;\subset\; \{U_1, U_2\}^*,
\]
and symmetrically on the backward arm \(Q_{1,\mathrm{walk}}^-(i)\), repeat the identical construction with the communicating bands \(v_1^{-1}, v_2^{-1}\) (bridged by \(u_{21}^{-1}, u_{12}^{-1}\)) to define a prefix-free dyadic family \(\mathcal{D}^- \subset \operatorname{Str}^-(i)\) which is order-isomorphic to \((\mathbb{Q}, <)\). Because \(B_1\) appears exclusively as the terminal block of every word in \(\mathcal{D}^+\), \emph{no word in \(\mathcal{D}^+\) is a proper prefix of any other word in \(\mathcal{D}^+\).} Consequently, any two distinct words \(W, W' \in \mathcal{D}^+\) diverge internally at their first differing block \(B_0 <_{\mathrm{str}} B_1 <_{\mathrm{str}} B_2\), where the comparison is determined by the fixed internal syllable divergence between \(U_1\) and \(U_2\). Therefore, the map
\[
B_{d_1} B_{d_2} \dots B_{d_m} B_1 \;\longmapsto\; \sum_{k=1}^m \frac{d_k}{2} 2^{-k} + 2^{-(m+1)}
\]
is a strict order-isomorphism from \((\mathcal{D}^+, <_{\mathrm{str}})\) onto the dyadic rationals in \((0,1)\), which is a countable dense linear order without endpoints, hence order-isomorphic to \((\mathbb{Q}, <)\) by Cantor's theorem. Thus \(\mathcal{C}^+ \coloneqq \{ \phi_{W^{\uparrow}} \mid W \in \mathcal{D}^+ \}\) and \(\mathcal{C}^- \coloneqq \{ \phi_{V^{\uparrow}} \mid V \in \mathcal{D}^- \}\) are two chains in \(\mathcal{L}_i(\Lambda)\), each order-isomorphic to \((\mathbb{Q}, <)\).

To prove that the canonical lattice homomorphism \(\iota \colon \mathcal{C}^- \otimes \mathcal{C}^+ \to \mathcal{L}_i(\Lambda)\) is injective, we verify Birkhoff's criterion for free products of two chains in a distributive lattice \cite[\S III.7]{Birkhoff1967}, \cite[\S IV.1]{Gratzer2003}: it suffices to show that for all \(a' < a\) in \(\mathcal{C}^-\) and \(b' < b\) in \(\mathcal{C}^+\),
\[
(a' \wedge b) \vee (a \wedge b') \;<\; a \wedge b \quad\text{in }\mathcal{L}_i(\Lambda).
\]
Fix such \(a' < a\) and \(b' < b\), writing \(a = \phi_{V_r^\uparrow}\), \(a' = \phi_{V_{r'}^\uparrow}\) with \(V_{r'} <_{\mathrm{str}} V_r\) in \(\mathcal{D}^-\) and \(b = \phi_{W_s^\uparrow}\), \(b' = \phi_{W_{s'}^\uparrow}\) with \(W_{s'} <_{\mathrm{str}} W_s\) in \(\mathcal{D}^+\). Let \(M_{r,s}\) be the finite-dimensional string module whose anchor basis vector \(v_0\) has backward arm \(V_r^\uparrow\) (read from \(v_0\) in the \(Q_{1,\mathrm{walk}}^-(i)\)-direction) and forward arm \(W_s^\uparrow\) (read from \(v_0\) in the \(Q_{1,\mathrm{walk}}^+(i)\)-direction). Equivalently, \(M_{r,s}\) is the string module obtained by identifying the terminal basis vector of \(M(V_r^\uparrow)\) with the initial basis vector of \(M(W_s^\uparrow)\), so that \(v_0\) is the common vertex \(i\). Because \(\mathcal{D}^\pm\) are prefix-free, \(V_{r'}^\uparrow \neq V_r^\uparrow\) and \(W_{s'}^\uparrow \neq W_s^\uparrow\); hence, by Prest--Puninski \cite[Lemma~2.3]{PrestPuninski2002} and our covariant convention,
\[
v_0 \in (a \wedge b)(M_{r,s}), \quad\text{while}\quad v_0 \notin a'(M_{r,s}) \text{ and } v_0 \notin b'(M_{r,s}).
\]
Since \(v_0 \in b(M_{r,s})\) but \(v_0 \notin a'(M_{r,s})\), we have \(v_0 \notin (a' \wedge b)(M_{r,s})\); similarly \(v_0 \notin (a \wedge b')(M_{r,s})\). Because pp-definable subspaces of \(M_{r,s}e_i\) are monomial \cite{CrawleyBoevey1989}, \(v_0 \notin [(a' \wedge b) + (a \wedge b')](M_{r,s})\). Hence \((a' \wedge b) \vee (a \wedge b') < a \wedge b\) in \(\mathcal{L}_i(\Lambda)\), verifying Birkhoff's criterion and the injectivity of \(\iota\).

The image of \(\iota\) is thus a freely generated distributive sublattice of \(\mathcal{L}_i(\Lambda)\) isomorphic to \(\mathcal{C}^- \otimes \mathcal{C}^+\). To bound the chains inside an explicit closed interval of \(\mathcal{L}_i(\Lambda)\), fix the two elements \(W_- \coloneqq B_0 B_1 <_{\mathrm{str}} W_+ \coloneqq B_2 B_1\) in \(\mathcal{D}^+\) (corresponding to the dyadic rationals \(1/4 < 3/4\)) and the corresponding elements \(V_- <_{\mathrm{str}} V_+\) in \(\mathcal{D}^-\). The closed subchains \(\widehat{\mathcal{C}}^+ \coloneqq \{\phi_{W^\uparrow} \mid W \in \mathcal{D}^+,\; W_- \le_{\mathrm{str}} W \le_{\mathrm{str}} W_+\} \subset \mathcal{C}^+\) and \(\widehat{\mathcal{C}}^- \coloneqq \{\phi_{V^\uparrow} \mid V \in \mathcal{D}^-,\; V_- \le_{\mathrm{str}} V \le_{\mathrm{str}} V_+\} \subset \mathcal{C}^-\) are bounded chains with minimum elements \(a_0 = \phi_{V_-^\uparrow},\, b_0 = \phi_{W_-^\uparrow}\) and maximum elements \(a_1 = \phi_{V_+^\uparrow},\, b_1 = \phi_{W_+^\uparrow}\), whose open interiors are order-isomorphic to \((\mathbb{Q}, <)\). By the Birkhoff criterion proved above, \(\widehat{\mathcal{C}}^- \otimes \widehat{\mathcal{C}}^+\) embeds into the interval \([a_0 \wedge b_0,\, a_1 + b_1] \subset \mathcal{L}_i(\Lambda)\) as a freely generated distributive sublattice.
\end{proof}

\begin{theorem}[Theorem D: Superdecomposable Pure-Injectives over Arbitrary Base Fields]
\label{thm:superdecomp_any_field}
Let \(\Lambda\) be a non-domestic string algebra over an algebraically closed field \(\K\) of \emph{arbitrary cardinality}. Then \(\Lambda\) admits a superdecomposable pure-injective module.
\end{theorem}

\begin{proof}
By Puninski's cardinality-free two-chain criterion \cite[Theorem~2.3]{Puninski2008}, if the pp-lattice \(\mathcal{L}(R)\) of a ring \(R\) contains an interval in which two bounded chains, each containing a densely ordered subchain, freely generate a distributive sublattice, then \(R\) possesses a superdecomposable pure-injective module (the construction uses only first-order compactness on the freely generated distributive lattice of chains, with no countability hypothesis on \(R\) or \(\K\)). By Lemma~\ref{lem:free_generation_monomial}, the interval \([a_0 \wedge b_0,\, a_1 + b_1] \subset \mathcal{L}_i(\Lambda)\) contains the free distributive product \(\widehat{\mathcal{C}}^- \otimes \widehat{\mathcal{C}}^+\) of two bounded chains containing subchains isomorphic to \((\mathbb{Q}, <)\). Consequently, \(\Lambda\) admits a superdecomposable pure-injective module.
\end{proof}

\begin{remark}[Resolution of Puninski's Open Problem]
\label{rem:puninski_resolved}
Theorem~\ref{thm:superdecomp_any_field} resolves in the affirmative the open question posed by Puninski \cite[Abstract and \S 4]{Puninski2008} (see also Prest's MathSciNet review MR2434248 of \cite{Puninski2008}) of whether every non-domestic string algebra over an uncountable field admits a superdecomposable pure-injective module.
\end{remark}

\subsection{Extension to Special Biserial and Brauer Graph Algebras}
\label{subsec:special_biserial_extension}

Recall \cite{WaldWaschbusch1985, ButlerRingel1987} that in any special biserial algebra \(\Lambda = \K\mathcal{Q}/I\), the non-monomial commutativity relations \(p_r - c_r q_r \in I\) (\(1 \le r \le m\)) define biserial projective-injective modules \(P_r\) with simple socles \(\operatorname{soc}(P_r) = \K \bar{p}_r\). Setting \(I_0 \coloneqq \bigoplus_{r=1}^m \operatorname{soc}(P_r)\), the quotient \(\overline{\Lambda} \coloneqq \Lambda / I_0\) is a finite-dimensional monomial string algebra.

\begin{corollary}[Extension to Special Biserial and Brauer Graph Algebras]
\label{cor:special_biserial_extension}
Let \(\Lambda\) be a special biserial algebra (in particular, any Brauer graph algebra) over an algebraically closed field \(\K\) of arbitrary cardinality, with canonical string quotient \(\overline{\Lambda} = \Lambda / I_0\). Then:
\begin{enumerate}
    \item \textbf{Clopen Decomposition:} \(\Zg(\Lambda) = \Zg(\overline{\Lambda}) \sqcup \{P_1, \dots, P_m\}\), where each biserial projective-injective module \(P_r = e_{s(p_r)}\Lambda\) is an isolated clopen point cut out by \((x = x e_{s(p_r)} / x \bar{p}_r = 0)\). Indeed, the idempotent generator \(x = e_{s(p_r)} \in P_r e_{s(p_r)}\) witnesses non-triviality of this open set on \(P_r\) (since \(e_{s(p_r)} \cdot \bar{p}_r = \bar{p}_r \neq 0\) spans \(\operatorname{soc}(P_r)\)), while any \(N \in \Zg(\Lambda)\) belonging to this open set satisfies \(N I_0 \neq 0\) and hence, by the argument below, \(N \cong P_r\). Consequently, \(D^\alpha(\Zg(\Lambda)) = D^\alpha(\Zg(\overline{\Lambda}))\) for all ordinals \(\alpha \ge 1\), and \(\KG(\Lambda) = \KG(\overline{\Lambda}) = \stgeo(\mathcal{S}_{\overline{\Lambda}}, \mathcal{M}_{\overline{\Lambda}}, \forb_{\overline{\Lambda}})\).
    \item \textbf{Superdecomposables and Shape Entropy:} \(\Lambda\) is non-domestic if and only if \(\overline{\Lambda}\) is non-domestic, in which case \(\Lambda\) admits a superdecomposable pure-injective module over \(\K\) of arbitrary cardinality, and Theorems~\ref{thm:ergodic_classification} and~\ref{thm:entropy_equivalence} hold on \(\mathcal{P}(\Zg(\Lambda)) = \mathcal{P}(\Zg(\overline{\Lambda}))\).
\end{enumerate}
\end{corollary}

\begin{proof}
(1) Because \(\Lambda \twoheadrightarrow \overline{\Lambda} = \Lambda / I_0\) is a surjective ring homomorphism with finitely generated two-sided ideal kernel \(I_0\), restriction of scalars embeds \(\Mod{\overline{\Lambda}}\) as the definable subcategory \(\{ M \in \Mod{\Lambda} \mid M I_0 = 0 \} \subseteq \Mod{\Lambda}\), inducing a closed topological embedding \(\Zg(\overline{\Lambda}) \hookrightarrow \Zg(\Lambda)\) \cite[Theorem~5.1.11 and Corollary~5.1.12]{Prest2009}. Conversely, if \(N \in \Zg(\Lambda)\) satisfies \(N I_0 \neq 0\), then \(m_0 \bar{p}_r \neq 0\) for some \(r \in \{1, \dots, m\}\) and \(m_0 \in N e_{s(p_r)}\). The \(\Lambda\)-homomorphism \(f \colon P_r = e_{s(p_r)}\Lambda \to N\) given by \(f(e_{s(p_r)}\lambda) = m_0 \lambda\) is non-zero on the simple essential socle \(\operatorname{soc}(P_r) = \K \bar{p}_r\), hence is a monomorphism. Since \(P_r\) is injective in \(\Mod{\Lambda}\), \(f\) splits, forcing \(N \cong P_r\) because \(N\) is indecomposable. The basic open set \((x = x e_{s(p_r)} / x \bar{p}_r = 0)\) equals \(\{P_r\}\), so each \(P_r\) is isolated and clopen in \(\Zg(\Lambda)\), giving \(\Zg(\Lambda) = \Zg(\overline{\Lambda}) \sqcup \{P_1, \dots, P_m\}\) and \(D^\alpha(\Zg(\Lambda)) = D^\alpha(\Zg(\overline{\Lambda}))\) for all \(\alpha \ge 1\).

(2) Immediate from (1) together with Theorems~\ref{thm:classification}, \ref{thm:superdecomp_any_field}, \ref{thm:ergodic_classification}, and~\ref{thm:entropy_equivalence} applied to \(\overline{\Lambda}\). Indeed, because \(I_0\) is finitely generated, every finitely presented \(\overline{\Lambda}\)-module is finitely presented as a \(\Lambda\)-module, so the full exact definable embedding \(\Mod{\overline{\Lambda}} \hookrightarrow \Mod{\Lambda}\) preserves pure-injectivity \cite[Theorem~5.1.11]{Prest2009} as well as direct-sum decompositions, mapping any superdecomposable pure-injective \(\overline{\Lambda}\)-module to a superdecomposable pure-injective \(\Lambda\)-module.
\end{proof}

\vspace{-0.5em}
\section{The Stable Rank of the Radical Filtration}
\label{sec:stable_rank}

In this section, we relate the geometric pruning derivatives \(\Dgeo^k(\mathcal{C}_{\mathrm{corr}})\) established in Section~\ref{sec:pruning_kg} to the transfinite powers of the Jacobson radical of the module category \(\md{\Lambda}\), connecting our geometric surface model directly to the foundational results of Schr\"oer \cite{Schroer2000a, Schroer2000b}, Laking--Prest--Puninski \cite{LPP2018}, and Gupta, Kuber, and Sardar \cite{GuptaKuberSardar2022}.

\subsection{The Algebraic Radical Filtration}
\label{subsec:radical_filtration}

\begin{definition}[Radical Filtration and Stable Rank]
\label{def:radical_filtration}
Let \(\Lambda\) be an artin algebra. The \emph{Jacobson radical} \(\rad_\Lambda\) of \(\md{\Lambda}\) is the two-sided ideal of morphisms generated by all non-invertible morphisms between indecomposables. The transfinite powers are defined by \(\rad_\Lambda^1 \coloneqq \rad_\Lambda\), \(\rad_\Lambda^{\alpha+1}(M, N) \coloneqq \sum_{X \in \md{\Lambda}} \rad_\Lambda^\alpha(X, N) \circ \rad_\Lambda(M, X)\) (composition read as applying the right factor first; \(X\) ranges over indecomposables), and \(\rad_\Lambda^\lambda \coloneqq \bigcap_{\alpha < \lambda} \rad_\Lambda^\alpha\) for limit ordinals \(\lambda\). When \((\rad_\Lambda^\omega)^k\) stabilizes (in particular, whenever \(\rad_\Lambda^\omega\) is nilpotent), one has \(\rad_\Lambda^{\omega \cdot k} = (\rad_\Lambda^\omega)^k\) for all \(k \ge 1\) \cite[\S 2]{LPP2018}. The \emph{stable rank} \(\operatorname{st}(\Lambda)\) is the smallest ordinal \(\nu\) such that \(\rad_\Lambda^\nu = \rad_\Lambda^{\nu+1}\), and the \emph{stable radical} is \(\rad_\Lambda^{\operatorname{st}(\Lambda)} = \rad_\Lambda^\infty\).
\end{definition}

By Schr\"oer \cite{Schroer2000a, Schroer2000b}, if \(\Lambda\) is a domestic string algebra, then \(\rad_\Lambda^\omega\) is nilpotent of finite index, so \(\operatorname{st}(\Lambda) < \omega^2\) and \(\rad_\Lambda^{\operatorname{st}(\Lambda)} = 0\). If \(\Lambda\) is non-domestic, then \(\rad_\Lambda^{\operatorname{st}(\Lambda)} \neq 0\) \cite{Schroer2000a, GuptaKuberSardar2022}.

\subsection{Radical Powers and Elementary Pivots}
\label{subsec:radical_pivots}

\begin{lemma}[Morphisms in \(\rad_\Lambda^k\), Pivot Chains, and \(\rad_\Lambda^\omega\)]
\label{lem:radical_pivots}
Let \(\Lambda\) be a string algebra modeled on \((\mathcal{S}, \mathcal{M}, \forb)\), and let \(M(w), M(w')\) be finite-dimensional string modules. Every morphism \(f \in \rad_\Lambda^k(M(w), M(w'))\) decomposes as a finite \(\K\)-linear combination of compositions of at least \(k\) irreducible morphisms between string modules (corresponding geometrically to sequences of at least \(k\) elementary boundary pivots), plus a morphism in \(\rad_\Lambda^\omega(M(w), M(w'))\). Moreover, \(f \in \rad_\Lambda^\omega\) if and only if \(f\) factors through arbitrarily long finite pivot sequences (or through a band tube), which geometrically means that \(f\) factors through a spiraling geodesic layer in \(\Dgeo^1(\mathcal{C}_{\mathrm{corr}})\).
\end{lemma}

\begin{proof}
By Proposition~\ref{prop:irreducible_pivot}, irreducible morphisms between string modules correspond bijectively to elementary boundary pivots. By Schr\"oer \cite{Schroer1999, Schroer2000a}, any morphism between a string module and a band module lies in \(\rad_\Lambda^\omega\), and a graph map between string modules lies in \(\rad_\Lambda^k \setminus \rad_\Lambda^\omega\) if and only if it factors as a chain of at least \(k\) irreducible morphisms (i.e., \(k\) elementary pivots). Taking the intersection \(\rad_\Lambda^\omega = \bigcap_{k < \omega} \rad_\Lambda^k\), a morphism lies in \(\rad_\Lambda^\omega\) if and only if it factors through arbitrarily long pivot sequences or through band modules, corresponding topologically to accumulation onto \(\Dgeo^1(\mathcal{C}_{\mathrm{corr}})\) (Theorem~\ref{thm:topological_convergence}).
\end{proof}

\subsection{Geometric Stabilization and the Radical Filtration}
\label{subsec:pruning_stabilization}

\begin{theorem}[Topological Characterization of Radical Stabilization]
\label{thm:stable_rank}
Let \(\Lambda\) be a string algebra modeled on \((\mathcal{S}, \mathcal{M}, \forb)\) with geometric pruning depth \(\stgeo = \stgeo(\mathcal{S}, \mathcal{M}, \forb) = \KG(\Lambda)\).
\begin{enumerate}
    \item \(\Lambda\) is representation-finite if and only if \(\stgeo = 0\), in which case \(\operatorname{st}(\Lambda) = N < \omega\), where \(N\) is the nilpotency index of \(\rad_\Lambda\) (the maximal irredundant pivot path length in the finite geometric fan \(\mathcal{F}\)), and \(\rad_\Lambda^N = 0\).
    \item \(\Lambda\) is \(n\)-domestic for some \(n \ge 1\) if and only if \(2 \le \stgeo = n + 1 < \omega\). In this case, each geometric pruning layer \(\Dgeo^k(\mathcal{C}_{\mathrm{corr}})\) (\(1 \le k \le \stgeo\)) corresponds to morphisms factoring through the \(\rad_\Lambda^{\omega \cdot k}\)-layer, and the infinite radical \(\rad_\Lambda^\omega\) is nilpotent of index \(k_0\) with \(\stgeo + 1 \le k_0 \le 2\stgeo - 1\). Consequently, \(\rad_\Lambda^{\operatorname{st}(\Lambda)} = 0\), and the stable rank satisfies
    \[
    \operatorname{st}(\Lambda) \;\in\; \bigl(\omega \cdot \stgeo, \, \omega \cdot (2\stgeo - 1)\bigr] \;\subsetneq\; [\omega \cdot \stgeo + 1, \, \omega^2).
    \]
    Equivalently, \(\operatorname{st}(\Lambda) = \omega \cdot \stgeo + c\) for some non-zero ordinal \(1 \le c \le \omega \cdot (\stgeo - 1)\).
    \item \(\Lambda\) is non-domestic if and only if \(\stgeo = \infty\). In this case, the non-empty perfect core \(\mathcal{P}_{\mathrm{geo}} \neq \emptyset\) obstructs the vanishing of the transfinite radical, so \(\rad_\Lambda^{\operatorname{st}(\Lambda)} \neq 0\) (and the powers \((\rad_\Lambda^\omega)^k\) are non-zero for all \(k \ge 1\)).
\end{enumerate}
\end{theorem}

\begin{proof}
(1) By Definition~\ref{def:cb_derivative} and Theorem~\ref{thm:classification}, \(\stgeo = \KG(\Lambda) = 0\) if and only if \(D^1(\Zg(\Lambda)) = \emptyset\), meaning \(\Lambda\) is representation-finite. Then \(\overline{\mathcal{F}} = \mathcal{F}\) is a finite poset, and \(\rad_\Lambda\) is nilpotent of index \(N = \operatorname{st}(\Lambda) < \omega\). For example, if \(\Lambda = \K A_N\) is the path algebra of a linearly oriented \(A_N\) quiver on \(N\) vertices (representation-finite), then the longest chain of irreducible morphisms with non-zero composition in \(\md{\Lambda}\) has length \(N-1\), so \(\stgeo = 0\) and \(\operatorname{st}(\Lambda) = N\).

(2) By Proposition~\ref{prop:n_domestic}, \(\Lambda\) is \(n\)-domestic if and only if \(2 \le \stgeo = n + 1 < \omega\) and the communication graph of closed geodesics has maximal chain length \(n-1\). By Lemma~\ref{lem:radical_pivots}, each transition from a preprojective spiraling cone into a closed geodesic tube or across a heteroclinic bridge contributes a factor in \(\rad_\Lambda^\omega\). Because the directed chain of closed geodesics has length at most \(n\), each closed geodesic can be entered and exited at most once in a finite-to-finite morphism (contributing at most two \(\rad_\Lambda^\omega\)-steps per band), plus one final step exiting into the preinjective cone, we get \((\rad_\Lambda^\omega)^{2n+1} = 0\). Thus \(\rad_\Lambda^\omega\) is nilpotent of index \(k_0 \le 2n + 1 = 2\stgeo - 1\). Since the top pruning layer \(\Dgeo^{\stgeo} \setminus \Dgeo^{\stgeo + 1} \neq \emptyset\) contributes nonzero morphisms in \((\rad_\Lambda^\omega)^{\stgeo} \setminus (\rad_\Lambda^\omega)^{\stgeo + 1}\), we have \(k_0 \ge \stgeo + 1\). Hence \(\operatorname{st}(\Lambda) \in \bigl(\omega \cdot \stgeo, \omega \cdot (2\stgeo - 1)\bigr]\), and \(\rad_\Lambda^{\operatorname{st}(\Lambda)} = 0\) \cite{Schroer2000a, LPP2018}.

(3) By Theorem~\ref{thm:exponential_growth}, \(\Lambda\) is non-domestic if and only if \(\stgeo = \infty\) and \(\mathcal{P}_{\mathrm{geo}} \neq \emptyset\). Because finite arcs can wind through communicating closed geodesics in arbitrarily long alternating chains, there exist non-zero morphisms between finite string modules factoring through \((\rad_\Lambda^\omega)^\alpha\) for arbitrarily high ordinals \(\alpha\), so \(\rad_\Lambda^{\operatorname{st}(\Lambda)} \neq 0\) \cite{Schroer2000a, GuptaKuberSardar2022}.
\end{proof}

\subsection{Example: Radical Stabilization of the Kronecker Algebra}
\label{subsec:stable_rank_examples}

\begin{example}[Radical Filtration of the Kronecker Algebra]
\label{ex:kronecker_radical}
Consider the \(1\)-domestic Kronecker algebra \(\Lambda = \K(1 \rightrightarrows 2)\) modeled on the two-cusped annulus, for which \(\stgeo = \KG(\Lambda) = 2\) (Example~\ref{ex:pruning_kronecker}):
\begin{enumerate}
    \item Morphisms between preprojective modules \(P_k \to P_{k+m}\) (or between preinjective modules) factor through \(m\) elementary pivots and lie in \(\rad_\Lambda^m \setminus \rad_\Lambda^{m+1}\).
    \item Morphisms from preprojective to regular band modules (or regular to preinjective) factor through the infinite spiraling rays \(\mathcal{L}_{\mathrm{spiral}} \in \Dgeo^1(\overline{\mathcal{F}})\) and lie in \(\rad_\Lambda^\omega \setminus \rad_\Lambda^{\omega \cdot 2}\).
    \item Morphisms from preprojective \(P_k\) to preinjective \(I_l\) factor through the central closed geodesic family \(\mathcal{L}_{\mathrm{band}} \in \Dgeo^2(\overline{\mathcal{F}})\) via two infinite spiraling transitions, lying in \((\rad_\Lambda^\omega)^2 = \rad_\Lambda^{\omega \cdot 2}\), so \((\rad_\Lambda^\omega)^2 \neq 0\). Because \(\Dgeo^3(\overline{\mathcal{F}}) = \emptyset\), one has \((\rad_\Lambda^\omega)^3 = \rad_\Lambda^{\omega \cdot 3} = 0\). Consequently, \(\operatorname{st}(\Lambda) = \omega \cdot 2 + c\) for some \(c \in [1, \omega]\), matching the sharp nilpotency index \((\rad_\Lambda^\omega)^3 = 0\) computed in \cite[\S 4]{Schroer2000a}.
\end{enumerate}
\end{example}

\vspace{-0.6em}
\section{Surface Symmetries, Mapping Class Groups, and the Ziegler Spectrum}
\label{sec:mapping_class_group}

In this section, we study how topological symmetries of the surface model \((\mathcal{S}, \mathcal{M}, \forb)\) act on the corridor space \(\mathcal{C}_{\mathrm{corr}}\) and the Ziegler spectrum \(\Zg(\Lambda)\). We show that general bordered mapping classes transport the corridor topology while preserving all transfinite CB-invariants, and that mapping classes preserving the oriented dissection \(\vec{\Delta}\) induce exact model-theoretic autoequivalences of \(\Zg(\Lambda)\).

\subsection{Mapping Class Group Action on the Corridor Space}
\label{subsec:mcg_action}

\begin{definition}[Bordered Mapping Class Group]
\label{def:mcg}
The \emph{bordered mapping class group} \(\operatorname{Mod}(\mathcal{S}, \mathcal{M}, \forb)\) is the group of isotopy classes of orientation-preserving homeomorphisms \(h \colon \mathcal{S} \xrightarrow{\sim} \mathcal{S}\) that preserve the marked segments \(\mathcal{M}\) and forbidden segments \(\forb\) setwise.
\end{definition}

Any mapping class \(h\) carries the cellular dissection \(\Delta\) to a transported dissection \(h(\Delta)\), defining a string algebra \(\Lambda_{h(\Delta)}\) modeled on the same surface.

\begin{proposition}[Topological Transport of the Corridor Space]
\label{prop:topological_transport}
Let \(h \in \operatorname{Mod}(\mathcal{S}, \mathcal{M}, \forb)\). The action of \(h\) on admissible arcs induces a canonical homeomorphism of global corridor spaces
\[
h_* \colon \bigl(\mathcal{C}_{\mathrm{corr}}^\Delta, \tau_{\mathrm{corr}}^\Delta\bigr) \;\xrightarrow{\;\sim\;}\; \bigl(\mathcal{C}_{\mathrm{corr}}^{h(\Delta)}, \tau_{\mathrm{corr}}^{h(\Delta)}\bigr)
\]
that commutes with the geometric pruning derivatives \(h_* \circ \Dgeo^\alpha = \Dgeo^\alpha \circ h_*\) for all ordinals \(\alpha\). Consequently, \(h_*\) preserves the pruning depth \(\stgeo\) and maps \(\mathcal{P}_{\mathrm{geo}}^\Delta\) homeomorphically onto \(\mathcal{P}_{\mathrm{geo}}^{h(\Delta)}\), yielding \(\KG(\Lambda_\Delta) = \KG(\Lambda_{h(\Delta)})\).
\end{proposition}

\begin{proof}
Because \(h\) preserves \(\mathcal{M}\) and \(\forb\), it maps admissible arcs to admissible arcs and test corridors to test corridors. Choosing geodesic representatives with respect to the fixed metric, the induced map \(h_*\) is a bijection on the basic open sets of \(\tau_{\mathrm{corr}}^\Delta\) and \(\tau_{\mathrm{corr}}^{h(\Delta)}\), hence a global homeomorphism. Since the Cantor-Bendixson pruning derivative \(\Dgeo\) is a purely topological invariant, \(h_*\) commutes with \(\Dgeo^\alpha\). By Theorem~\ref{thm:classification}, the pruning depths match.
\end{proof}

\subsection{Oriented-Dissection Symmetries and Autoequivalences of \(\Zg(\Lambda)\)}
\label{subsec:autoequivalences}

\begin{theorem}[Dissection Symmetries Induce Model-Theoretic Autoequivalences]
\label{thm:mcg_autoequivalence}
Let \(\operatorname{Mod}(\mathcal{S}, \mathcal{M}, \forb, \vec{\Delta})\) be the subgroup of \(\operatorname{Mod}(\mathcal{S}, \mathcal{M}, \forb)\) that preserves the oriented cellular dissection \(\vec{\Delta}\) setwise. Every \(h \in \operatorname{Mod}(\mathcal{S}, \mathcal{M}, \forb, \vec{\Delta})\) induces a quiver automorphism preserving \(\rho\), lifting to an exact autoequivalence of \(\Mod{\Lambda}\) and a self-homeomorphism
\[
\Phi_h \;\coloneqq\; \Psi \circ h_* \circ \Phi \colon \Zg(\Lambda) \;\xrightarrow{\;\sim\;}\; \Zg(\Lambda)
\]
(extended across band tubes via \(\gamma_v \mapsto \gamma_{h(v)}\)). The map \(\Phi_h\) commutes with all CB-derivatives \(D^\alpha\), preserves \(\KG(\Lambda)\) and \(\operatorname{st}(\Lambda)\), and induces an automorphism of the lattice of definable subcategories of \(\Mod{\Lambda}\).
\end{theorem}

\begin{proof}
By Definition~\ref{def:generalized_surface}, the quiver \(\mathcal{Q}\) and relations \(\rho\) are uniquely determined by the oriented ribbon gluing. Thus, \(h\) permutes the dissection arcs while preserving their incidences, inducing a quiver automorphism \(\sigma_h \colon \mathcal{Q} \xrightarrow{\sim} \mathcal{Q}\) preserving \(\rho\). Restriction of scalars along \(\sigma_h^{-1}\) gives an exact autoequivalence of \(\Mod{\Lambda}\). It restricts to \(\Phi_h\) on \(C\), maps band tubes homeomorphically, and since it preserves basic open sets, it commutes with \(D^\alpha\) and preserves definable subcategories \cite{Ziegler1984, Prest2009}.
\end{proof}

\begin{remark}[General Mapping Classes and Derived Equivalences]
\label{rem:general_mapping_classes}
For general \(h \notin \operatorname{Mod}(\mathcal{S}, \mathcal{M}, \forb, \vec{\Delta})\), \(\Lambda_{h(\Delta)}\) need not be isomorphic to \(\Lambda_\Delta\). In the gentle case \cite{OPS2018}, \(h\) induces a derived equivalence \(D^b(\Lambda_\Delta) \cong D^b(\Lambda_{h(\Delta)})\), and Proposition~\ref{prop:topological_transport} gives a geometric explanation for the invariance of \(\KG\) under derived equivalence.
\end{remark}

\subsection{Dehn Twists and Twisting of pp-Definable Subgroups}
\label{subsec:dehn_twists}

\begin{definition}[Dehn Twist]
\label{def:dehn_twist}
Let \(c \subset \mathcal{S}^\circ \setminus \forb\) be an essential simple closed curve corresponding to a primitive band \(v \in \operatorname{Ba}(\Lambda)\). The \emph{Dehn twist} \(\tau_c \in \operatorname{Mod}(\mathcal{S}, \mathcal{M}, \forb)\) is obtained by cutting \(\mathcal{S}\) along \(c\), rotating one side by \(2\pi\), and gluing back.
\end{definition}

\begin{figure}[htbp]
\centering
\begin{tikzpicture}[scale=0.85]
    \begin{scope}[shift={(-3.6,0)}]
        \draw[thick] (0,0) circle (2.2cm);
        \draw[thick] (0,0) circle (0.75cm);
        \draw[blue, very thick] (0,0) circle (1.35cm);
        
        \node[blue, above] at (0,1.65) {\small \(c\)};
        
        \draw[red!85!black, thick, domain=0:10.5, samples=500, variable=\x]
            plot ({ (1.35 + 0.85*exp(-0.32*\x)) * cos(90 - \x*(180/pi)) },
                  { (1.35 + 0.85*exp(-0.32*\x)) * sin(90 - \x*(180/pi)) });
                  
        \node[red!85!black, font=\small] at (0, -2.6) {$\mathcal{L}_{\mathrm{spiral}}$};
    \end{scope}

    \draw[-{Stealth[length=2.8mm]}, very thick] (-0.8,0) -- (0.8,0) node[midway, above] {\(\tau_c\)};

    \begin{scope}[shift={(3.6,0)}]
        \draw[thick] (0,0) circle (2.2cm);
        \draw[thick] (0,0) circle (0.75cm);
        \draw[blue, very thick] (0,0) circle (1.35cm);
        
        \node[blue, above] at (0,1.65) {\small \(c\)};
        
        \draw[red!85!black, thick, domain=0:16.6, samples=600, variable=\x]
            plot ({ (1.35 + 0.85*exp(-0.20*\x)) * cos(90 - \x*(180/pi)) },
                  { (1.35 + 0.85*exp(-0.20*\x)) * sin(90 - \x*(180/pi)) });
                  
        \node[red!85!black, font=\small] at (0, -2.6) {$\tau_c(\mathcal{L}_{\mathrm{spiral}})$};
    \end{scope}
\end{tikzpicture}
\caption[The action of a Dehn twist on a spiraling geodesic lamination]{The action of a Dehn twist \(\tau_c\) on a spiraling geodesic lamination. The core closed geodesic \(c\) is fixed, while the spiraling leaf accumulates with a shifted winding filtration.}
\label{fig:dehn_twist}
\end{figure}
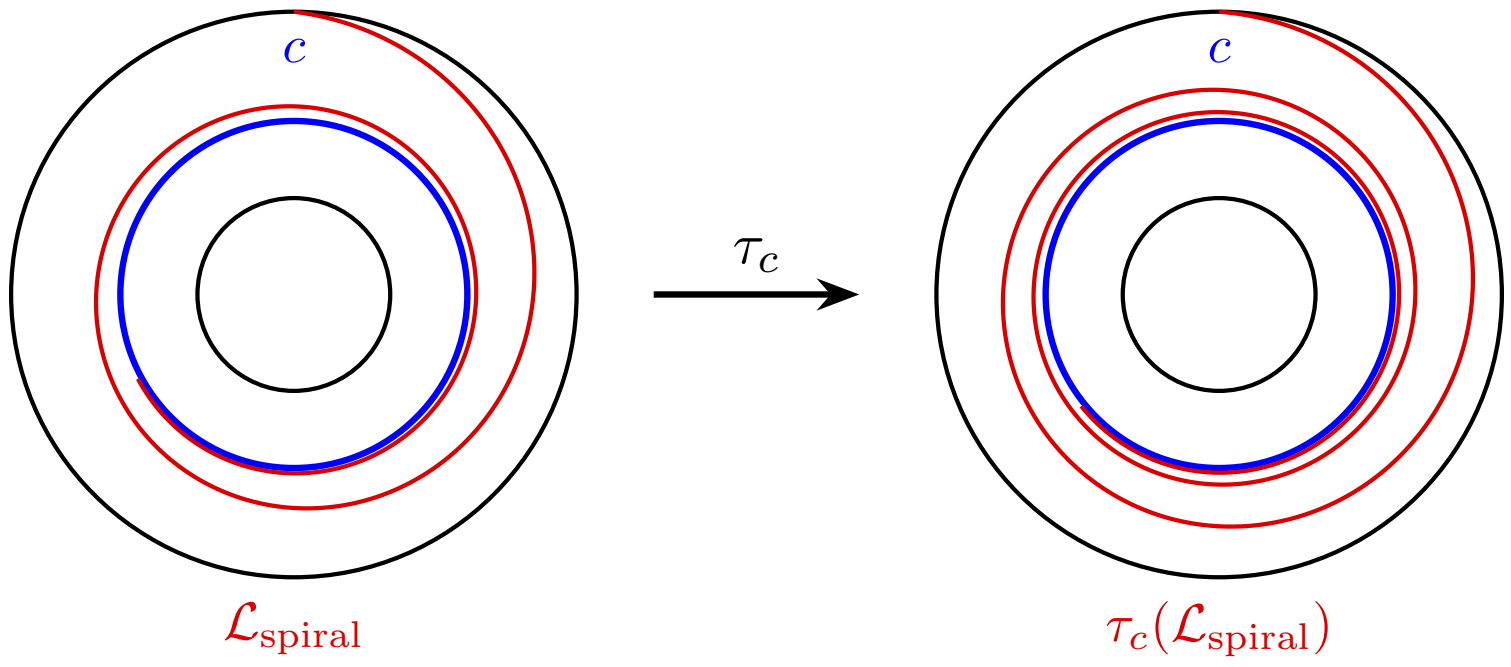

\begin{example}[Dehn Twists on the Kronecker Algebra]
\label{ex:dehn_twist_kronecker}
Consider the Kronecker algebra \(\Lambda = \K(1 \rightrightarrows 2)\) on the two-cusped annulus (Example~\ref{ex:kronecker_torus}). The core curve \(c = \gamma_v\) is the unique primitive band \(v = \beta^{-1}\alpha\). (Here \(\tau_c\) acts on the corridor space \(\mathcal{C}_{\mathrm{corr}}^\Delta\) via the canonical identification \(\mathcal{C}_{\mathrm{corr}}^{h(\Delta)} \cong \mathcal{C}_{\mathrm{corr}}^\Delta\) arising from the symmetry of the Kronecker algebra.) The Dehn twist \(\tau_c\) shears the annulus, adding one full wind around \(c\) to every intersecting geodesic (Figure~\ref{fig:dehn_twist}). On the finite arcs of CB-rank \(0\), \(\tau_c\) shifts the winding number (\(\gamma_n \mapsto \gamma_{n+1}\)), realizing (with appropriate orientation) the Auslander-Reiten translation \(\tau^{\pm 1}\) on the string hammock (Lemma~\ref{lem:topological_ar}). In the limits \(\Dgeo^1\) and \(\Dgeo^2\), \(\tau_c\) fixes the core geodesic \(c = \gamma_v\) (the band generic module \(G_v\)) and shifts the finite-prefix pp-definable filtration inside the spiraling periodic ray and coray modules. On non-domestic surfaces with intersecting curves \(c_1, c_2\), the twist \(\tau_{c_1}\) permutes the minimal arational laminations inside \(\mathcal{P}_{\mathrm{geo}}\).
\end{example}

\vspace{-0.9em}

\vspace{-0.4em}
\section{Measured Laminations, Ergodic Decomposition, and Asymptotic Profiles}
\label{sec:measured}

In this section, we enrich the topological dictionary with quantitative ergodic invariants. By equipping minimal arational laminations with Thurston's transverse invariant measures \cite{Thurston1988, CassonBleiler1988}, we show that asymptotic dimension vectors and pp-subquotient growth frequencies along finite-prefix exhaustions of arational pure-injective modules are governed by the Birkhoff Ergodic Theorem across the Choquet simplex \(\mathcal{M}_{\mathrm{erg}}(\mathcal{L})\).

\subsection{The Thurston Intersection and Corridor Measure Maps}
\label{subsec:thurston_maps}

Let \(\mathcal{ML}(\mathcal{S}^\circ, \forb)\) denote Thurston's space of measured geodesic laminations supported on \(\mathcal{S}^\circ \setminus \forb\). A measured lamination \((\mathcal{L}, \mu)\) consists of a geodesic lamination \(\mathcal{L} \subset \mathcal{S}^\circ \setminus \forb\) equipped with a non-zero, transverse-invariant Borel measure \(\mu\) on transverse arcs with \(\operatorname{supp}(\mu) = \mathcal{L}\). The space \(\mathcal{ML}(\mathcal{S}^\circ, \forb)\) is topologized by weak-\(*\) convergence on transverse arcs, and its projectivization modulo positive scaling \(\R_{>0}\) is denoted \(\mathcal{PML}(\mathcal{S}^\circ, \forb)\).

Normalized by total dissection weight \(|\mu|_\Delta \coloneqq \sum_{j \in Q_0} \mu(d_j) > 0\), each projective class \([\mu] \in \mathcal{PML}(\mathcal{S}^\circ, \forb)\) has a unique representative satisfying \(|\mu|_\Delta = 1\). We associate two quantitative profiles to \((\mathcal{L}, \mu)\):
\begin{enumerate}
    \item \textbf{The Normalized Intersection Vector:} The map \(\Phi_\Delta \colon \mathcal{ML}(\mathcal{S}^\circ, \forb) \to \R_{\ge 0}^{Q_0}\) defined by \(\Phi_\Delta(\mu)_i \coloneqq \mu(d_i)\) for \(i \in Q_0\), with projectivized representative \(\underline{d}(\mu) \coloneqq \Phi_\Delta(\mu) / |\mu|_\Delta \in [0,1]^{Q_0}\).
    \item \textbf{The Corridor Measure Profile:} The map \(\widehat{\Phi}_\Delta \colon \mathcal{ML}(\mathcal{S}^\circ, \forb) \to \prod_{i \in Q_0} \R_{\ge 0}^{\mathcal{L}_i(\Lambda)}\) assigning to each pointed string pp-pair \((\phi/\psi) \in \mathcal{L}_i(\Lambda)\) the transverse measure \(\widehat{\Phi}_\Delta(\mu)_{\phi/\psi} \coloneqq \mu(\mathcal{C}_{\phi/\psi} \cap d_i)\) of the leaves traversing the geometric test corridor \(\mathcal{C}_{\phi/\psi}\) across \(d_i\), with normalized frequency \(\rho_{\phi/\psi}(\mu) \coloneqq \mu(\mathcal{C}_{\phi/\psi} \cap d_i) / |\mu|_\Delta \in [0,1]\).
\end{enumerate}

\begin{lemma}[Injectivity of the Corridor Measure Map]
\label{lem:corridor_measure_injectivity}
The projectivized corridor measure map
\[
\mathbb{P}\widehat{\Phi}_\Delta \colon \mathcal{PML}(\mathcal{S}^\circ, \forb) \;\hookrightarrow\; \mathbb{P}\Bigl(\prod_{i \in Q_0} \R_{\ge 0}^{\mathcal{L}_i(\Lambda)}\Bigr)
\]
is a topological embedding.
\end{lemma}

\begin{proof}
Suppose \((\mathcal{L}, \mu)\) and \((\mathcal{L}', \mu')\) are normalized so that \(|\mu|_\Delta = |\mu'|_\Delta = 1\) and satisfy \(\widehat{\Phi}_\Delta(\mu) = \widehat{\Phi}_\Delta(\mu')\). Because \(\operatorname{supp}(\mu) = \mathcal{L}\), a finite two-sided string template \(u\) of radius \(n\) is traversed by a leaf of \(\mathcal{L}\) if and only if the corresponding cylinder corridor \(\mathcal{C}_u \cap d_i\) is a non-empty open subset of \(\mathcal{L} \cap d_i\), which holds if and only if \(\mu(\mathcal{C}_u \cap d_i) > 0\). Thus \(\widehat{\Phi}_\Delta(\mu) = \widehat{\Phi}_\Delta(\mu')\) implies that \(\mathcal{L}\) and \(\mathcal{L}'\) traverse the exact same admissible words of every radius \(n \ge 1\). By Theorem~\ref{thm:theorem_B_sec5}, a complete geodesic is uniquely determined by its infinite dissection crossing sequence, so \(\mathcal{L} = \mathcal{L}'\) (for non-minimal \(\mathcal{L}\), the same argument applies separately to each minimal component and each spiraling leaf, using the fact that any two distinct leaves diverge at some finite radius).

On the common support \(\mathcal{L}\), Lemma~\ref{lem:corridors_form_basis} shows that the intersections \(\mathcal{C}_{\phi/\psi} \cap (\mathcal{L} \cap d_i)\) of pointed string test corridors form a clopen basis for the Cantor cross-section \(\mathcal{L} \cap d_i\). The finite disjoint unions of cylinder corridors form an algebra of subsets of \(\mathcal{L}\cap d_i\) that generates the Borel \(\sigma\)-algebra. Since \(\mu\) and \(\mu'\) agree on each cylinder corridor, they agree on this algebra, and hence \(\mu=\mu'\) by the Carath\'eodory extension theorem. Thus \(\mu\) and \(\mu'\) agree on every transversal \(d_i\) (\(i \in Q_0\)). Because the complement of \(\Delta = \bigsqcup_{i \in Q_0} d_i\) in \(\mathcal{S}^\circ \setminus \forb\) consists of ideal polygons bounded by \(\mathcal{M}\) and \(\forb\) that contain no complete geodesics in their interior, \(\bigsqcup_{i \in Q_0} d_i\) is a complete transversal for every lamination \(\mathcal{L} \subset \mathcal{S}^\circ \setminus \forb\), so \(\mu = \mu'\) globally by transverse invariance. Continuity of \(\mathbb{P}\widehat{\Phi}_\Delta\) follows from the clopenness of cylinder corridors on \(\mathcal{L} \cap d_i\).
\end{proof}

\begin{remark}[When is each of \(\mathbb{P}\Phi_\Delta\) and \(\mathbb{P}\widehat{\Phi}_\Delta\) an embedding?]
\label{rem:when_which_embedding}
The corridor measure map \(\mathbb{P}\widehat{\Phi}_\Delta\) is a topological embedding unconditionally (Lemma~\ref{lem:corridor_measure_injectivity}). It combines vertex weights with corner arc weights (via test corridors of every template) and therefore eliminates all turn ambiguities, at the cost of embedding into a larger projective space.

The vertex intersection map \(\mathbb{P}\Phi_\Delta \colon \mathcal{PML}(\mathcal{S}^\circ, \forb) \to \mathbb{P}(\R_{\ge 0}^{Q_0})\) alone is not in general injective. Two measured laminations can pass through a polygon \(P\) of the dissection \(\Delta\) with \(\ge 4\) sides along different corner pairings while having identical transverse weights on every dissection arc \(d_j\); relations of length \(k \ge 3\) in \(\Lambda\) are precisely the data that produce such polygons. However, \(\mathbb{P}\Phi_\Delta\) is a topological embedding when \(\Delta\) restricts to an ideal triangulation of the lamination-supporting subsurface (in particular, when \(\Lambda\) is gentle), in which case the edge weights \(i(\mu, d_j)\) uniquely determine the corner arc weights via \(a_{jk} = \frac{1}{2}(i(\mu, d_j) + i(\mu, d_k) - i(\mu, d_l))\) \cite[\S 4]{Thurston1988}.
\end{remark}

\subsection{The Finite-Level Identity}
\label{subsec:finite_level_identity}

\begin{proposition}[The Finite-Level Identity]
\label{prop:finite_level_identity}
Let \(w \in \operatorname{Str}(\Lambda)\) be a finite admissible string with geodesic arc \(\gamma_w\), and let \((\phi/\psi)\) be a pointed string pp-pair at vertex \(i\) of maximal template radius \(L = L(\phi/\psi)\). Then
\[
\dim_\K \bigl(M(w)e_i\bigr) \;=\; I(\gamma_w, d_i),
\]
\[
\dim_\K \left(\frac{\phi(M(w))}{\psi(M(w))}\right) \;=\; I\bigl(\gamma_w, \mathcal{C}_{\phi/\psi}\bigr) \;+\; O\bigl(L(\phi/\psi)\bigr),
\]
where \(I(\gamma_w, d_i)\) and \(I(\gamma_w, \mathcal{C}_{\phi/\psi})\) denote the geometric intersection numbers with \(d_i\) and the interior crossings through \(\mathcal{C}_{\phi/\psi}\), respectively. The error term satisfies \(\bigl|\dim_\K(\phi(M(w))/\psi(M(w))) - I(\gamma_w, \mathcal{C}_{\phi/\psi})\bigr| \le 2L(\phi/\psi)\) and accounts strictly for basis vectors whose anchor position lies within \(L(\phi/\psi)\) syllables of either endpoint of \(w\).
\end{proposition}

\begin{proof}
By Theorem~\ref{thm:bls_realization} and Definition~\ref{def:oriented_intersections}, \(\gamma_w\) is in minimal position with respect to \(\Delta\), so the canonical basis of \(M(w)e_i\) is in exact bijection with \(\gamma_w \cap d_i\), giving \(\dim_\K(M(w)e_i) = I(\gamma_w, d_i)\). By Crawley-Boevey \cite{CrawleyBoevey1989}, \(\dim_\K(\phi(M(w))/\psi(M(w)))\) equals the number of canonical basis vectors \(v_r \in M(w)e_i\) belonging to \(\phi(M(w)) \setminus \psi(M(w))\). Whether \(v_r \in \phi(M(w)) \setminus \psi(M(w))\) depends only on the radius-\(L\) neighborhood of \(v_r\) inside \(w\). For every anchor \(v_r\) at syllable distance strictly greater than \(L\) from both endpoints, its radius-\(L\) neighborhood is complete, so by Lemma~\ref{lem:geometric_pp_evaluation}, \(v_r \in \phi(M(w)) \setminus \psi(M(w))\) if and only if the corresponding crossing traverses the corridor \(\mathcal{C}_{\phi/\psi}\). Because at most \(L\) crossings can occur within \(L\) syllables of each endpoint, the two counts differ by at most \(2L = O(L(\phi/\psi))\).
\end{proof}

\subsection{Asymptotic Dimension Vectors in the Uniquely Ergodic Case}
\label{subsec:uniquely_ergodic}

A minimal arational lamination \(\mathcal{L}\) is \emph{uniquely ergodic} if it supports a unique transverse invariant measure \(\mu\) up to scaling.

\begin{theorem}[Asymptotic Profiles in the Uniquely Ergodic Case]
\label{thm:unique_ergodic_profile}
Let \(N = \Psi(\mathcal{L}) \in \mathcal{P}(\Zg(\Lambda))\) be associated with a uniquely ergodic minimal arational lamination \((\mathcal{L}, [\mu])\) (normalized so that \(|\mu|_\Delta = 1\)). For \emph{every} leaf \(\ell \subset \mathcal{L}\) and \emph{every} sequence of finite subword prefixes \(w_n\) of \(\vec{w}_\ell\) with \(|w_n| \to \infty\),
\begin{align*}
\lim_{n \to \infty} \frac{\underline{\dim}\, M(w_n)}{\dim_\K M(w_n)} 
  &= \underline{d}(\mu) = \Phi_\Delta(\mu), \\
\lim_{n \to \infty} \frac{\dim_\K\bigl(\phi(M(w_n))/\psi(M(w_n))\bigr)}{\dim_\K M(w_n)} 
  &= \rho_{\phi/\psi}(\mu) = \widehat{\Phi}_\Delta(\mu)_{\phi/\psi}.
\end{align*}
\end{theorem}

\begin{proof}
Let \(X = \mathcal{L} \cap \bigsqcup_{j \in Q_0} d_j\) be the compact Cantor cross-section, and let \(T \colon X \to X\) be the first-return map. Unique ergodicity of \((\mathcal{L}, \mu)\) is equivalent to unique ergodicity of \((X, T, \mu)\). By Oxtoby's theorem for uniquely ergodic systems, for every clopen cylinder set \(A \subset X\), the empirical visitation averages \(\frac{1}{N}\sum_{k=0}^{N-1} \mathbf{1}_A(T^k x)\) converge uniformly for all \(x \in X\) to \(\mu(A)/\mu(X) = \mu(A)\). Applying Proposition~\ref{prop:finite_level_identity} with \(A = \mathcal{L} \cap d_i\) and \(A = \mathcal{C}_{\phi/\psi} \cap (\mathcal{L} \cap d_i)\), and noting that \(\dim_\K M(w_n) = |w_n| + 1\) while \(O(L(\phi/\psi)) / \dim_\K M(w_n) \to 0\), yields the claimed limits for every leaf. (Here \(|w_n|\) equals the number of first-return visits to \(X\) along the leaf \(\ell\), up to the bounded additive error \(O(L(\phi/\psi))\) arising from the boundary effect of Proposition~\ref{prop:finite_level_identity}.)
\end{proof}

\begin{corollary}[Completeness on uniquely ergodic points]
\label{cor:unique_ergodic_complete}
Let \([\mathcal{L},\mu],[\mathcal{L}',\mu']\in\mathcal{PML}(\mathcal{S}^\circ,\forb)\) be uniquely ergodic minimal arational measured laminations (normalized so that \(|\mu|_\Delta=|\mu'|_\Delta=1\)). If \(\mathbb{P}\widehat{\Phi}_\Delta(\mu)=\mathbb{P}\widehat{\Phi}_\Delta(\mu')\), then \(\mathcal{L}=\mathcal{L}'\) and \([\mu]=[\mu']\).
\end{corollary}

\begin{proof}
By Lemma~\ref{lem:corridor_measure_injectivity}, \(\mathbb{P}\widehat{\Phi}_\Delta\) is injective, so this forces \([\mathcal{L},\mu]=[\mathcal{L}',\mu']\) in \(\mathcal{PML}\). The support of \([\mathcal{L},\mu]\) is \(\mathcal{L}\), so \(\mathcal{L}=\mathcal{L}'\); and since each \(\mathcal{L}\) is uniquely ergodic with normalized transverse measure, \([\mu]=[\mu']\).
\end{proof}

\subsection{Ergodic Classification for Non-Uniquely Ergodic Laminations}
\label{subsec:ergodic_classification}

In general, the space of normalized transverse invariant measures on a minimal arational lamination \(\mathcal{L}\) is a finite-dimensional Choquet simplex whose extreme points \(\mathcal{M}_{\mathrm{erg}}(\mathcal{L})\) are the ergodic transverse probability measures.

\begin{theorem}[Ergodic Decomposition of Generic Exhaustions]
\label{thm:ergodic_classification}
Let \(N = \Psi(\mathcal{L}) \in \mathcal{P}(\Zg(\Lambda))\) be associated with a minimal arational lamination \(\mathcal{L} \subset \mathcal{S}^\circ \setminus \forb\). For each \([\mu] \in \mathcal{M}_{\mathrm{erg}}(\mathcal{L})\) (normalized by \(|\mu|_\Delta = 1\)), there exists a \(\mu\)-full-measure, topologically dense union of generic leaves \(\mathcal{L}_\mu \subseteq \mathcal{L}\) such that for every leaf \(\ell \in \mathcal{L}_\mu\) and every sequence of finite prefix exhaustions \(w_n\) of \(\vec{w}_\ell\) with \(|w_n| \to \infty\):
\begin{align*}
\lim_{n \to \infty} \frac{\underline{\dim}\, M(w_n)}{\dim_\K M(w_n)} &= \underline{d}(\mu) = \Phi_\Delta(\mu), \\
\lim_{n \to \infty} \frac{\dim_\K \bigl(\phi(M(w_n))/\psi(M(w_n))\bigr)}{\dim_\K M(w_n)} &= \rho_{\phi/\psi}(\mu) = \widehat{\Phi}_\Delta(\mu)_{\phi/\psi}
\end{align*}
for all pointed string pp-pairs \((\phi/\psi) \in \mathcal{L}_i(\Lambda)\). Moreover, the map \([\mu] \mapsto \rho_\bullet(\mu) = \mathbb{P}\widehat{\Phi}_\Delta(\mu)\) is injective on \(\mathcal{M}_{\mathrm{erg}}(\mathcal{L})\), so distinct ergodic classes partition \(\mathcal{L}\) into disjoint generic leaf sets \(\mathcal{L}_\mu \cap \mathcal{L}_{\mu'} = \emptyset\) with strictly distinct asymptotic pp-growth profiles.
\end{theorem}

\begin{proof}
Let \((X, T)\) be the first-return dynamical system on \(X = \mathcal{L} \cap \bigsqcup_{j \in Q_0} d_j\). Because there are only countably many pointed string pp-pairs, applying the Birkhoff Pointwise Ergodic Theorem to \(\mathbf{1}_{\mathcal{C}_{\phi/\psi} \cap d_i}\) for each pair and taking the countable intersection of the resulting full-\(\mu\)-measure generic sets produces a \(T\)-invariant generic set \(X_\mu \subseteq X\) with \(\mu(X_\mu) = 1\). Let \(\mathcal{L}_\mu \subseteq \mathcal{L}\) be the union of leaves passing through \(X_\mu\). Since \(\operatorname{supp}(\mu) = \mathcal{L}\), \(\mathcal{L}_\mu\) is dense. For any leaf \(\ell \in \mathcal{L}_\mu\), Proposition~\ref{prop:finite_level_identity} implies that the empirical ratios differ from the Birkhoff time averages by \(O(L(\phi/\psi))/|w_n| \to 0\), yielding the limits. Finally, by Lemma~\ref{lem:corridor_measure_injectivity}, \(\mathbb{P}\widehat{\Phi}_\Delta\) is injective, so \([\mu] \neq [\mu']\) implies \(\rho_\bullet(\mu) \neq \rho_\bullet(\mu')\) and hence \(\mathcal{L}_\mu \cap \mathcal{L}_{\mu'} = \emptyset\).
\end{proof}

\vspace{-0.5em}
\section{Algebraic pp-Entropy, Shape-Injectivity, and Symbolic Lamination Dynamics}
\label{sec:entropy}

Because every arational pure-injective module \(N \in \mathcal{P}(\Zg(\Lambda))\) satisfies \(\CB(N) = \infty\), the Cantor-Bendixson rank cannot distinguish points inside the perfect core. In this section, we introduce a purely algebraic invariant of an indecomposable pure-injective module \(N \in \mathcal{P}_i \coloneqq \mathcal{P}(\Zg(\Lambda)) \cap \overline{H_\Lambda(S(i))}\)---the \emph{localized pp-shape entropy} \(h^{\mathrm{shape}}(N)\)---measuring the exponential growth rate of distinct pp-definable subgroups of \(N e_i\) cut out by radius-\(n\) templates. By proving the Basis Identity Lemma and the Shape-Injectivity Theorem, we show that \(h^{\mathrm{shape}}(N)\) coincides with the topological coding entropy \(h_{\mathrm{cod}}(\mathcal{L})\) of the corresponding minimal arational lamination \(\mathcal{L}\).

\subsection{Explicit Interval Formulas and Localized pp-Shape Entropy}
\label{subsec:explicit_interval_formulas}

Fix a vertex \(i \in Q_0\) and an oriented dissection arc \(\vec{d}_i\). By Definition~\ref{def:generalized_surface}(1), the incident syllables at vertex \(i\) partition into two disjoint sets \(Q_{1,\mathrm{walk}}^-(i) \sqcup Q_{1,\mathrm{walk}}^+(i)\) corresponding to the two opposite walk directions across the ribbon strip \(R_i\). For each integer \(n \ge 1\), let \(\operatorname{Str}_n^\pm(i)\) denote the finite set of one-sided admissible strings of length \(n\) starting at vertex \(i\) with initial syllable in \(Q_{1,\mathrm{walk}}^\pm(i)\) respectively, and let
\[
\operatorname{Str}_{2n}(i) \;\coloneqq\; \Bigl\{\, u = (u^-, u^+) \;\Big|\; u^- \in \operatorname{Str}_n^-(i),\; u^+ \in \operatorname{Str}_n^+(i),\; (u^-)^{-1}u^+ \in \operatorname{Str}(\Lambda) \,\Bigr\}
\]
denote the set of two-sided admissible templates of radius \(n\) centered at vertex \(i\).

We also consider the set \(\widetilde{\operatorname{Str}}_{2n}(i)\) of oriented words \(w\) of total length \(2n\) whose midpoint vertex is \(i\). Because every word \(w \in \widetilde{\operatorname{Str}}_{2n}(i)\) crosses \(d_i\) at its midpoint in one of the two orientations, either \(w = (u^-)^{-1} u^+\) (crossing from \(-\) to \(+\)) or \(w^{-1} = (u^-)^{-1} u^+\) (crossing from \(+\) to \(-\)). Thus the word-inversion involution \(w \mapsto \bar{w} \coloneqq w^{-1}\) on \(\widetilde{\operatorname{Str}}_{2n}(i)\) preserves the underlying two-sided arm pair \(\{u^-, u^+\}\), with fixed points given by the \emph{involutive palindromes} \(P_n \subseteq \widetilde{\operatorname{Str}}_{2n}(i)\) satisfying \(w = \bar{w}\).

Following the covariant Prest--Puninski string order \(\le_{\mathrm{str}}\) \cite[\S 2]{PrestPuninski2002} (Section~\ref{subsec:uncountable_superdecomp}), on \(\operatorname{Str}_n^\pm(i)\) we have \(v_1 \le_{\mathrm{str}} v_2 \iff \phi_{v_1}(M) \subseteq \phi_{v_2}(M)\) for all \(M \in \Mod{\Lambda}\). For each \(v \in \operatorname{Str}_n^\pm(i)\), the cone of all one-sided strings having prefix \(v\) has a unique \(\le_{\mathrm{str}}\)-maximum element \(v^\uparrow \coloneqq v \beta_1^{-1}\dots\beta_s^{-1}\) (the maximal inverse extension of \(v\)), and an immediate lower-complement bound in the finite linear order \((\operatorname{Str}_n^\pm(i), \le_{\mathrm{str}})\):
\[
v^\sharp \;\coloneqq\; \Bigl(\max_{\le_{\mathrm{str}}}\bigl\{\, v' \in \operatorname{Str}_n^\pm(i) \;\big|\; v' <_{\mathrm{str}} v \,\bigr\}\Bigr)^\uparrow
\]
(setting \(\phi_{v^\sharp}(x) \equiv 0\) if \(v\) is the \(\le_{\mathrm{str}}\)-minimal element).

For each two-sided template \(u = (u^-, u^+) \in \operatorname{Str}_{2n}(i)\) (or oriented word \(w \in \widetilde{\operatorname{Str}}_{2n}(i)\) with arms \((u^-, u^+)\)), we define the explicitly bounded interval pp-pair \(\phi_{u^\uparrow}(x) / \psi_{u^\downarrow}(x)\) at sort \(i\) by:
\begin{align}
\phi_{u^{\uparrow}}(x) &\;\coloneqq\; \phi_{(u^-)^{\uparrow}}(x) \;\wedge\; \phi_{(u^+)^{\uparrow}}(x), \label{eq:phi_u_up} \\
\psi_{u^{\downarrow}}(x) &\;\coloneqq\; \bigl(\phi_{(u^-)^{\sharp}}(x) \wedge \phi_{(u^+)^{\uparrow}}(x)\bigr) \;+\; \bigl(\phi_{(u^-)^{\uparrow}}(x) \wedge \phi_{(u^+)^{\sharp}}(x)\bigr). \label{eq:psi_u_down}
\end{align}

\begin{definition}[Localized pp-Shape Entropy]
\label{def:shape_entropy}
Let \(N \in \mathcal{P}_i\). For each \(n \ge 1\), let
\[
\mathcal{T}_n(N) \;\coloneqq\; \Bigl\{\, w \in \widetilde{\operatorname{Str}}_{2n}(i) \;\Big|\; \frac{\phi_{w^{\uparrow}}(N)}{\psi_{w^{\downarrow}}(N)} \neq 0 \,\Bigr\}
\]
denote the set of active radius-\(n\) oriented templates on \(N\) (which is finite since \(\widetilde{\operatorname{Str}}_{2n}(i)\) is finite), with cardinality \(N_n(N) \coloneqq |\mathcal{T}_n(N)|\) and combinatorial word-growth rate \(h(N) \coloneqq \lim_{n \to \infty} \frac{1}{2n}\log N_n(N)\).

Two active templates \(w, w' \in \mathcal{T}_n(N)\) are \emph{\(N\)-shape-equivalent}, written \(w \sim_N w'\), if they cut out the exact same pair of pp-definable \(\K\)-subspaces inside \(N e_i\):
\[
w \sim_N w' \quad \Longleftrightarrow \quad \phi_{w^{\uparrow}}(N) = \phi_{w'^{\uparrow}}(N) \;\;\text{ and }\;\; \psi_{w^{\downarrow}}(N) = \psi_{w'^{\downarrow}}(N).
\]
Letting \(S_n(N) \coloneqq |\mathcal{T}_n(N) / \!\sim_N|\) be the number of distinct pp-definable subgroup pairs cut out on \(N e_i\) at radius \(n\) (and abbreviating \(P_n(N)\coloneqq P_n\cap\mathcal{T}_n(N)\) for the active palindromic templates), the \emph{localized pp-shape entropy} of \(N\) at sort \(i\) is defined as
\[
h^{\mathrm{shape}}(N) \;\coloneqq\; \lim_{n \to \infty} \frac{1}{2n} \log S_n(N).
\]
\end{definition}

\subsection{The Basis Identity Lemma and the Shape-Injectivity Theorem}
\label{subsec:shape_injectivity}

\begin{lemma}[Basis Identity Lemma]
\label{lem:basis_identity}
Let \(N = \Psi(\mathcal{L}) = \PE(M(\overleftrightarrow{w})) \in \mathcal{P}_i\) be an arational pure-injective module. For each canonical basis vector \(e_k \in M(\overleftrightarrow{w})e_i\), let \(w_n^-(k) \in \operatorname{Str}_n^-(i)\) and \(w_n^+(k) \in \operatorname{Str}_n^+(i)\) denote the length-\(n\) prefixes of \(\overleftrightarrow{w}\) exiting \(e_k\) along \(Q_{1,\mathrm{walk}}^-(i)\) and \(Q_{1,\mathrm{walk}}^+(i)\) respectively. Then for every template \(u = (u^-, u^+) \in \operatorname{Str}_{2n}(i)\):
\begin{enumerate}
    \item The subspaces \(\phi_{u^{\uparrow}}(M(\overleftrightarrow{w}))\) and \(\psi_{u^{\downarrow}}(M(\overleftrightarrow{w}))\) are spanned by subsets of the canonical basis \(\{e_k\}\) of \(M(\overleftrightarrow{w})e_i\).
    \item A canonical basis vector \(e_k \in M(\overleftrightarrow{w})e_i\) belongs to the difference \(\phi_{u^{\uparrow}}(M(\overleftrightarrow{w})) \setminus \psi_{u^{\downarrow}}(M(\overleftrightarrow{w}))\) if and only if
    \[
    \bigl(w_n^-(k),\, w_n^+(k)\bigr) \;=\; \bigl(u^-,\, u^+\bigr).
    \]
    In particular, when \(w \in \mathcal{T}_n(N)\) (i.e.\ \(\phi_{w^\uparrow}(N)/\psi_{w^\downarrow}(N) \neq 0\)), at least one such occurrence exists.
\end{enumerate}
\end{lemma}

\begin{proof}
(1) By Crawley-Boevey \cite{CrawleyBoevey1989} and Prest--Puninski \cite{PrestPuninski2002}, every pointed string formula \(\phi_v(M(\overleftrightarrow{w}))\) is spanned by a subset of the canonical basis of \(M(\overleftrightarrow{w})e_i\). Because intersections and sums of monomial subspaces with respect to a fixed basis remain monomial, \(\phi_{u^\uparrow}(M(\overleftrightarrow{w}))\) and \(\psi_{u^\downarrow}(M(\overleftrightarrow{w}))\) are spanned by subsets of the canonical basis \(\{e_k\}\).

(2) Because \(Q_{1,\mathrm{walk}}^-(i) \cap Q_{1,\mathrm{walk}}^+(i) = \emptyset\), any graph map \(M((u^\pm)^\uparrow) \to M(\overleftrightarrow{w})\) mapping the anchor to \(e_k\) must send the initial syllable of \(u^\pm\) into the unique arm of \(\overleftrightarrow{w}\) at \(e_k\) whose initial syllable lies in \(Q_{1,\mathrm{walk}}^\pm(i)\). Thus the evaluation of \(\phi_{(u^-)^\uparrow}\) and \(\phi_{(u^-)^\sharp}\) at \(e_k\) depends exclusively on the backward ray, while the evaluation of \(\phi_{(u^+)^\uparrow}\) and \(\phi_{(u^+)^\sharp}\) depends exclusively on the forward ray.

By Prest and Puninski \cite[\S 2 and Lemma~2.3]{PrestPuninski2002}, a one-sided ray \(w^\pm(k)\) at \(e_k\) with length-\(n\) prefix \(w_n^\pm(k)\) satisfies \(e_k \in \phi_{(u^\pm)^\uparrow}(M(\overleftrightarrow{w}))\) if and only if \(w^\pm(k) \le_{\mathrm{str}} (u^\pm)^\uparrow\). Because \((u^\pm)^\uparrow\) is the maximal inverse extension, this holds if and only if \(w_n^\pm(k) \le_{\mathrm{str}} u^\pm\). Similarly, \(e_k \in \phi_{(u^\pm)^\sharp}(M(\overleftrightarrow{w}))\) if and only if \(w_n^\pm(k) <_{\mathrm{str}} u^\pm\). Consequently:
\begin{align*}
e_k \in \phi_{u^\uparrow}\bigl(M(\overleftrightarrow{w})\bigr) &\quad\Longleftrightarrow\quad w_n^-(k) \le_{\mathrm{str}} u^- \;\;\text{ and }\;\; w_n^+(k) \le_{\mathrm{str}} u^+, \\[4pt]
e_k \in \psi_{u^\downarrow}\bigl(M(\overleftrightarrow{w})\bigr) &\quad\Longleftrightarrow\quad 
\begin{aligned}[t]
&\bigl(w_n^-(k) <_{\mathrm{str}} u^- \text{ and } w_n^+(k) \le_{\mathrm{str}} u^+\bigr) \\
&\text{ or } \bigl(w_n^-(k) \le_{\mathrm{str}} u^- \text{ and } w_n^+(k) <_{\mathrm{str}} u^+\bigr).
\end{aligned}
\end{align*}
Subtracting the two conditions shows that \(e_k \in \phi_{u^\uparrow}(M(\overleftrightarrow{w})) \setminus \psi_{u^\downarrow}(M(\overleftrightarrow{w}))\) if and only if \(w_n^-(k) = u^-\) and \(w_n^+(k) = u^+\). (The backward and forward implications both hold: the prefix \(u^\pm\) forces \(w^\pm(k) \le_{\mathrm{str}} (u^\pm)^\uparrow\) by maximality of \((u^\pm)^\uparrow\), giving membership in \(\phi_{u^\uparrow}\); and since \((u^\pm)^\sharp\) has prefix strictly below \(u^\pm\) in \(\le_{\mathrm{str}}\), the prefix \(u^\pm\) rules out membership in \(\phi_{(u^\pm)^\sharp}\), giving non-membership in \(\psi_{u^\downarrow}\).)
\end{proof}

\begin{theorem}[Shape-Injectivity Theorem]
\label{thm:shape_injectivity}
Let \(N = \Psi(\mathcal{L}) \in \mathcal{P}_i\) be an arational pure-injective module. For every active oriented template \(w \in \mathcal{T}_n(N)\), its \(N\)-shape equivalence class is given precisely by the inversion fiber \([w]_N = \{w, \bar{w}\}\) (where \(\bar{w} = w^{-1}\)). Letting \(P_n(N) = \{w \in \mathcal{T}_n(N) \mid w = \bar{w}\}\) denote the palindromic fixed-point subset, we have the exact identity
\[
S_n(N) \;=\; \frac{N_n(N) + |P_n(N)|}{2}, \qquad \text{and consequently} \qquad h^{\mathrm{shape}}(N) \;=\; h(N).
\]
\end{theorem}

\begin{proof}
Let \(w \in \mathcal{T}_n(N)\) have two-sided arm pair \(u = (u^-, u^+) \in \operatorname{Str}_{2n}(i)\). Reversing the reading orientation of \(w\) gives \(\bar{w} = w^{-1} \in \mathcal{T}_n(N)\), which has the exact same two-sided arm pair. Since the formulas \(\phi_{u^\uparrow}\) and \(\psi_{u^\downarrow}\) in \eqref{eq:phi_u_up}--\eqref{eq:psi_u_down} are symmetric in the two conjunction/addition factors, \(\phi_{w^\uparrow} = \phi_{\bar{w}^\uparrow}\) and \(\psi_{w^\downarrow} = \psi_{\bar{w}^\downarrow}\) identically, so \(\bar{w} \in [w]_N\).

Conversely, suppose \(w, w' \in \mathcal{T}_n(N)\) satisfy \(w \sim_N w'\), with respective arm pairs \(u = (u^-, u^+)\) and \(u' = (u'^-, u'^+)\). By definition of \(w \sim_N w'\), we have \(\phi_{u^\uparrow}(N) = \phi_{u'^\uparrow}(N)\) and \(\psi_{u^\downarrow}(N) = \psi_{u'^\downarrow}(N)\). Because \(M(\overleftrightarrow{w}) \hookrightarrow \PE(M(\overleftrightarrow{w})) = N\) is a pure embedding, for every pp-formula \(\chi \in \mathcal{L}_i(\Lambda)\) we have \(\chi(M(\overleftrightarrow{w})) = \chi(N) \cap M(\overleftrightarrow{w})\). Intersecting with \(M(\overleftrightarrow{w})\) therefore yields
\[
\phi_{u^\uparrow}\bigl(M(\overleftrightarrow{w})\bigr) \setminus \psi_{u^\downarrow}\bigl(M(\overleftrightarrow{w})\bigr) \;=\; \phi_{u'^\uparrow}\bigl(M(\overleftrightarrow{w})\bigr) \setminus \psi_{u'^\downarrow}\bigl(M(\overleftrightarrow{w})\bigr).
\]
Because \(w \in \mathcal{T}_n(N)\), we have \(\phi_{u^\uparrow}(N)/\psi_{u^\downarrow}(N) \neq 0\), so by purity \(\phi_{u^\uparrow}(M(\overleftrightarrow{w}))/\psi_{u^\downarrow}(M(\overleftrightarrow{w})) \neq 0\). By Lemma~\ref{lem:basis_identity}(1), both subspaces are spanned by subsets of the canonical basis of \(M(\overleftrightarrow{w})e_i\), so there exists at least one canonical basis vector \(e_k \in M(\overleftrightarrow{w})e_i\) lying in \(\phi_{u^\uparrow}(M(\overleftrightarrow{w})) \setminus \psi_{u^\downarrow}(M(\overleftrightarrow{w}))\). By the equality above, \(e_k\) also lies in \(\phi_{u'^\uparrow}(M(\overleftrightarrow{w})) \setminus \psi_{u'^\downarrow}(M(\overleftrightarrow{w}))\). Applying Lemma~\ref{lem:basis_identity}(2) to both memberships at the same basis vector \(e_k\) forces
\[
(u^-, u^+) \;=\; \bigl(w_n^-(k), w_n^+(k)\bigr) \;=\; (u'^-, u'^+).
\]
Thus \(w\) and \(w'\) have the exact same arm pair, which means \(w' \in \{w, \bar{w}\}\). Hence \([w]_N = \{w, \bar{w}\}\). The involution \(w \mapsto \bar{w}\) partitions \(\mathcal{T}_n(N)\) into \(|P_n(N)|\) singleton orbits and \((N_n(N) - |P_n(N)|)/2\) two-element orbits, so \(S_n(N) = (N_n(N) + |P_n(N)|)/2\). Taking \(\frac{1}{2n}\log(\cdot)\) yields \(h^{\mathrm{shape}}(N) = h(N)\).
\end{proof}

\subsection{Equivalence with Topological Lamination Entropy and the Variational Principle}
\label{subsec:variational_principle}

Let \(\Sigma_m(\mathcal{L})\) denote the set of all admissible words of syllable length \(m \ge 1\) traversed by leaves of the minimal recurrent lamination \(\mathcal{L}\), let \(\Sigma_{2n}(\mathcal{L}; i) \subseteq \widetilde{\operatorname{Str}}_{2n}(i)\) denote the subset of words of length \(2n\) centered at \(d_i\), and let
\[
h_{\mathrm{cod}}(\mathcal{L}) \;\coloneqq\; \lim_{m \to \infty} \frac{1}{m} \log \bigl|\Sigma_m(\mathcal{L})\bigr|
\]
be the topological coding entropy of the first-return subshift \((X, T)\) on the global Cantor cross-section \(X = \mathcal{L} \cap \bigsqcup_{j \in Q_0} d_j\).

\begin{theorem}[Theorem E(1): Equivalence of Shape Entropy and Lamination Entropy]
\label{thm:entropy_equivalence}
For any arational pure-injective module \(N = \Psi(\mathcal{L}) \in \mathcal{P}_i\), the limit defining \(h^{\mathrm{shape}}(N)\) exists, is independent of the supporting vertex \(i \in Q_0\), and equals both the topological coding entropy \(h_{\mathrm{cod}}(\mathcal{L})\) and the supremum of Kolmogorov--Sinai measure-theoretic entropies across the Choquet simplex \(\mathcal{M}_{\mathrm{erg}}(\mathcal{L})\):
\[
h^{\mathrm{shape}}(N) \;=\; h_{\mathrm{cod}}(\mathcal{L}) \;=\; \sup_{[\mu] \in \mathcal{M}_{\mathrm{erg}}(\mathcal{L})} h_{\mathbb{P}_\mu}(T).
\]
\end{theorem}

\begin{proof}
By Lemma~\ref{lem:basis_identity}(2) and minimality of \(\mathcal{L}\) (every leaf \(\overleftrightarrow{w}\) is dense in \(\mathcal{L}\) and therefore realizes all finite subwords in \(\Sigma_{2n}(\mathcal{L}; i)\) at its vertex-\(i\) anchors), an oriented template \(w \in \widetilde{\operatorname{Str}}_{2n}(i)\) belongs to \(\mathcal{T}_n(N)\) if and only if \(w \in \Sigma_{2n}(\mathcal{L}; i)\). Thus \(N_n(N) = |\Sigma_{2n}(\mathcal{L}; i)|\) for all \(n \ge 1\). Because the full language \(\Sigma(\mathcal{L})\) is a shift-invariant closed subshift, the sequence \(a_m = \log|\Sigma_m(\mathcal{L})|\) is subadditive (\(a_{m+k} \le a_m + a_k\)), so Fekete's subadditivity lemma guarantees that the limit \(h_{\mathrm{cod}}(\mathcal{L}) = \lim_{m \to \infty} \frac{1}{m}\log|\Sigma_m(\mathcal{L})|\) exists. Moreover, by minimality of \(\mathcal{L}\) and compactness of \(X\), the clopen slice \(X_i = \mathcal{L} \cap d_i \subseteq X\) is syndetic at the syllable level: there exists \(K \ge 1\) such that every leaf segment of length \(K\) syllables crosses \(d_i\) at least once, so every \(u \in \Sigma_{2n}(\mathcal{L})\) extends to some \(w \in \Sigma_{2(n+K)}(\mathcal{L}; i)\) at an offset \(r \in \{-K, \dots, K\}\) from the midpoint of \(w\). Hence
\[
\bigl|\Sigma_{2n}(\mathcal{L}; i)\bigr| \;\le\; \bigl|\Sigma_{2n}(\mathcal{L})\bigr| \;\le\; (2K + 1)\,\bigl|\Sigma_{2(n+K)}(\mathcal{L}; i)\bigr|.
\]
Taking \(\frac{1}{2n}\log(\cdot)\) and letting \(n \to \infty\) gives \(h(N) = \lim_{n \to \infty} \frac{1}{2n}\log N_n(N) = h_{\mathrm{cod}}(\mathcal{L})\), and Theorem~\ref{thm:shape_injectivity} yields \(h^{\mathrm{shape}}(N) = h(N) = h_{\mathrm{cod}}(\mathcal{L})\).

Finally, because the symbolic coding map from the global Cantor cross-section \((X, T)\) of \(\mathcal{L}\) to the subshift \((\Sigma(\mathcal{L}), \sigma)\) is a continuous factor map that is injective outside the boundary-leaf set of the complementary ideal polygons (which is at most countable and therefore carries zero measure for every non-atomic transverse measure \(\mu\)), it is a measure-theoretic isomorphism for every \([\mu] \in \mathcal{M}_{\mathrm{erg}}(\mathcal{L})\). Applying the classical Goodwyn--Dinaburg--Goodman Variational Principle \cite[Theorem~8.6]{Walters1982} to \((X, T)\) yields \(h_{\mathrm{cod}}(\mathcal{L}) = \sup_{[\mu] \in \mathcal{M}_{\mathrm{erg}}(\mathcal{L})} h_{\mathbb{P}_\mu}(T)\).
\end{proof}

\begin{remark}[Coherence of Algebraic Shape Entropy and Ergodic Profiles]
\label{rem:sec11_sec12_coherence}
Theorems~\ref{thm:ergodic_classification} and \ref{thm:entropy_equivalence} illuminate the dual dynamical structure on the perfect core. The localized pp-shape entropy \(h^{\mathrm{shape}}(N)\) is a purely topological invariant of \(N\): it counts how many distinct pp-definable subgroup pairs are non-trivially cut out inside \(N e_i\), independent of how frequently a leaf visits each corridor. Conversely, the asymptotic dimension vectors and pp-growth profiles of Section~\ref{sec:measured} weight each corridor by its asymptotic visitation frequency along finite-prefix exhaustions, thereby separating the generic leaf sets of distinct ergodic transverse measures. The Variational Principle bridges the two invariants directly.
\end{remark}

\subsection{The Localized pp-Shape Pressure Function, Microcanonical \texorpdfstring{\(K_0\)}{K0}-Entropy, and Iso-Entropic Separation}
\label{subsec:pp_pressure}

While the localized pp-shape entropy \(h^{\mathrm{shape}}(N) \in [0, \infty)\) is a single scalar---and vanishes identically on simple (non-self-intersecting) minimal arational surface laminations \(\mathcal{L} \subset \mathcal{S}^\circ \setminus \forb\), whose first-return maps are interval exchange transformations of linear complexity \(N_n(N) = O(n)\)---the Shape-Injectivity Theorem (Theorem~\ref{thm:shape_injectivity}) allows us to stratify \(\mathcal{S}_n(N) \coloneqq \mathcal{T}_n(N)/\!\sim_N\) by the Grothendieck group \(K_0(\md{\Lambda}) \cong \Z^{Q_0}\) and promote \(h^{\mathrm{shape}}(N)\) to a convex pressure function on \(\Hom_{\Z}(K_0(\md{\Lambda}), \R) \cong \R^{Q_0}\). Throughout this subsection, \(\mathcal{L} \in \mathcal{P}_{\mathrm{geo}}\) denotes any minimal recurrent geodesic system in the geometric perfect core (either a simple minimal arational lamination or a minimal Cantor geodesic subshift, consistent with Section~\ref{sec:superdecomposables} and Section~\ref{sec:measured}). We write \(\Sigma^{Q_0} \coloneqq \{d \in \R_{\ge 0}^{Q_0} : \sum_{j \in Q_0} d_j = 1\}\) for the standard probability simplex in \(\R^{Q_0}\), avoiding conflict with the cellular dissection \(\Delta\).

\begin{lemma}[Well-Definedness of Witness Modules and Grothendieck Classes on \(\mathcal{S}_n(N)\)]
\label{lem:witness_module_shape}
Let \(N = \Psi(\mathcal{L}) \in \mathcal{P}_i\) be an arational pure-injective module. For each \(n \ge 1\) and each \(N\)-shape equivalence class \(\mathfrak{s} = [w]_N \in \mathcal{S}_n(N)\), the finite-dimensional string module \(M(\mathfrak{s}) \coloneqq M(w)\) is uniquely determined up to \(\Lambda\)-module isomorphism by the subgroup pair \(\mathfrak{s} = \bigl(\phi_{w^\uparrow}(N), \psi_{w^\downarrow}(N)\bigr)\) alone. In particular, every \(N\)-shape class \(\mathfrak{s} \in \mathcal{S}_n(N)\) carries a canonical Grothendieck class
\[
\bigl[M(\mathfrak{s})\bigr] \;=\; \underline{\dim}\, M(\mathfrak{s}) \;\in\; K_0\bigl(\md{\Lambda}\bigr) \;\cong\; \Z^{Q_0}
\]
satisfying \(\sum_{j \in Q_0} \dim_\K\bigl(M(\mathfrak{s})e_j\bigr) = 2n + 1\).
\end{lemma}

\begin{proof}
By Theorem~\ref{thm:shape_injectivity}, for any two active templates \(w, w' \in \mathcal{T}_n(N)\), we have \(w \sim_N w'\) if and only if \(w' \in \{w, \bar{w}\}\), where \(\bar{w} = w^{-1}\). By Definition~\ref{def:string_band_modules}, reversing an admissible string yields an isomorphic string module \(M(w) \cong M(w^{-1})\). Therefore, the isomorphism class of \(M(\mathfrak{s}) \coloneqq M(w)\) in \(\md{\Lambda}\) and its dimension vector \(\underline{\dim}\, M(\mathfrak{s}) \in \Z^{Q_0}\) depend only on the \(N\)-shape equivalence class \(\mathfrak{s} \in \mathcal{S}_n(N)\). Because \(w\) has syllable length \(2n\), \(M(\mathfrak{s})\) has total \(\K\)-dimension \(2n + 1\).
\end{proof}

Lemma~\ref{lem:witness_module_shape} allows us to define three purely algebraic invariants of \(N \in \mathcal{P}_i\) from \(\mathcal{L}_i(\Lambda)\) and \(K_0(\md{\Lambda})\) without reference to the surface model:

\begin{definition}[Asymptotic Grothendieck Support, Microcanonical \(K_0\)-Shape Entropy, and Pressure Function]
\label{def:pp_pressure}
Let \(N \in \mathcal{P}_i\), and for each \(n \ge 1\) let \(\mathcal{D}_n^{\mathrm{alg}}(N) \coloneqq \bigl\{\frac{[M(\mathfrak{s})]}{2n+1} \;\big|\; \mathfrak{s} \in \mathcal{S}_n(N)\bigr\} \subset \Sigma^{Q_0} \subset K_0(\md{\Lambda})_\R\) be the set of normalized Grothendieck classes of the radius-\(n\) witness modules of \(N\).
\begin{enumerate}
    \item \textbf{Asymptotic Grothendieck Support and Intrinsic \(K_0\)-Rank:} The \emph{asymptotic Grothendieck support} of \(N\) is the Kuratowski upper limit of \(\mathcal{D}_n^{\mathrm{alg}}(N)\) in \(K_0(\md{\Lambda})_\R\)\,\allowbreak\(\cong \R^{Q_0}\):
    \[
    \operatorname{Supp}_{K_0}^\infty(N) \;\coloneqq\; \bigcap_{m=1}^\infty \overline{\bigcup_{n=m}^\infty \mathcal{D}_n^{\mathrm{alg}}(N)},
    \qquad
    \operatorname{rk}^\infty_{K_0}(N) \;\coloneqq\; \dim_\R \operatorname{span}_\R\bigl(\operatorname{Supp}_{K_0}^\infty(N)\bigr).
    \]
    \item \textbf{Microcanonical \(K_0\)-Shape Entropy:} For any direction \(d \in \Sigma^{Q_0}\) and tolerance \(\varepsilon > 0\), let
    \[
    S_n(N; d, \varepsilon) \;\coloneqq\; \Bigl\vert{}\Bigl\{\, \mathfrak{s} \in \mathcal{S}_n(N) \;\Big\vert{}\; \Bigl\Vert{}\frac{[M(\mathfrak{s})]}{2n+1} - d\Bigr\Vert{}_\infty \le \varepsilon \,\Bigr\}\Bigr\vert{},
    \]
    and define the \emph{microcanonical \(K_0\)-shape entropy} of \(N\) at \(d\) by
    \[
    \mathfrak{h}_{K_0}(N; d) \;\coloneqq\; \lim_{\varepsilon \downarrow 0} \limsup_{n \to \infty} \frac{1}{2n} \log S_n(N; d, \varepsilon) \;\in\; \{-\infty\} \cup [0, \infty).
    \]
    \item \textbf{Localized pp-Shape Partition and Pressure Functions:} For any additive Grothendieck weight \(\theta = (\theta_j)_{j \in Q_0} \in \Hom_{\Z}(K_0(\md{\Lambda}), \R) \cong \R^{Q_0}\), paired with \([M] \in K_0(\md{\Lambda})\) via \(\langle \theta, [M] \rangle \coloneqq \sum_{j \in Q_0} \theta_j \dim_\K(Me_j)\), define the partition sums
    \[
    Z_n^{\mathrm{shape}}(N; \theta) \;\coloneqq\; \sum_{\mathfrak{s} \in \mathcal{S}_n(N)} \exp\Bigl(\bigl\langle \theta,\, [M(\mathfrak{s})] \bigr\rangle\Bigr),
    \]
    \(Z_n(N; \theta) \coloneqq \sum_{w \in \mathcal{T}_n(N)} \exp\bigl(\langle \theta, [M(w)] \rangle\bigr)\), and \(Z_n^{\mathrm{pal}}(N; \theta) \coloneqq \sum_{w \in P_n(N)} \exp\bigl(\langle \theta, [M(w)] \rangle\bigr)\). The \emph{localized pp-shape pressure function} of \(N\) is
    \[
    P^{\mathrm{shape}}_N(\theta) \;\coloneqq\; \lim_{n \to \infty} \frac{1}{2n} \log Z_n^{\mathrm{shape}}(N; \theta).
    \]
\end{enumerate}
\end{definition}

Note that at \(\theta = 0\), \(Z_n^{\mathrm{shape}}(N; 0) = S_n(N)\) and \(P^{\mathrm{shape}}_N(0) = h^{\mathrm{shape}}(N)\).

\begin{proposition}[Weighted Shape-Injectivity Identity, Existence, and Vertex-Invariance]
\label{prop:weighted_shape_injectivity}
Let \(N = \Psi(\mathcal{L}) \in \mathcal{P}_i\) be an arational pure-injective module. For every \(\theta \in \R^{Q_0}\):
\begin{enumerate}
    \item \textbf{Weighted Shape-Injectivity Identity:} For every \(n \ge 1\),
    \[
    Z_n^{\mathrm{shape}}(N; \theta) \;=\; \frac{Z_n(N; \theta) + Z_n^{\mathrm{pal}}(N; \theta)}{2},
    \]
    and in particular \(\frac{1}{2} Z_n(N; \theta) \le Z_n^{\mathrm{shape}}(N; \theta) \le Z_n(N; \theta)\).
    \item \textbf{Existence and Vertex-Invariance:} The limit \(P^{\mathrm{shape}}_N(\theta) = \lim_{n \to \infty} \frac{1}{2n} \log Z_n^{\mathrm{shape}}(N; \theta)\) exists and is finite for all \(\theta \in \R^{Q_0}\), is a convex \(1\)-Lipschitz function with respect to \(\Vert{}\theta\Vert{}_\infty\), and is independent of the choice of vertex \(i \in Q_0\) in the support of \(N\).
\end{enumerate}
\end{proposition}

\begin{proof}
(1) By Theorem~\ref{thm:shape_injectivity}, the involution \(w \mapsto \bar{w} = w^{-1}\) partitions \(\mathcal{T}_n(N)\) into \(|P_n(N)|\) singleton classes \(\{w\}\) (\(w \in P_n(N)\)) and \((N_n(N) - |P_n(N)|)/2\) two-element classes \(\{w, \bar{w}\}\) (\(w \notin P_n(N)\)). By Lemma~\ref{lem:witness_module_shape}, \(\langle \theta, [M(w)] \rangle = \langle \theta, [M(\bar{w})] \rangle\) is constant on each equivalence class \(\mathfrak{s} = [w]_N\). Summing \(\exp(\langle \theta, [M(\mathfrak{s})] \rangle)\) over \(\mathfrak{s} \in \mathcal{S}_n(N)\) counts each \(w \in P_n(N)\) with weight \(1\) and each pair \(\{w, \bar{w}\}\) (\(w \notin P_n(N)\)) once instead of twice, giving \(2 Z_n^{\mathrm{shape}}(N; \theta) = Z_n(N; \theta) + Z_n^{\mathrm{pal}}(N; \theta)\). Since \(0 \le Z_n^{\mathrm{pal}}(N; \theta) \le Z_n(N; \theta)\), the inequality \(\frac{1}{2} Z_n(N; \theta) \le Z_n^{\mathrm{shape}}(N; \theta) \le Z_n(N; \theta)\) follows immediately.

(2) Let \(\Sigma_m(\mathcal{L})\) denote the set of all admissible words \(u\) of syllable length \(m \ge 1\) traversed by leaves of \(\mathcal{L}\) (without restricting the midpoint vertex), and define the full subshift partition sum
\[
Z_m^{\mathrm{full}}(\mathcal{L}; \theta) \;\coloneqq\; \sum_{u \in \Sigma_m(\mathcal{L})} \exp\Bigl(\bigl\langle \theta,\, [M(u)] \bigr\rangle\Bigr).
\]
If \(u \in \Sigma_{m+k}(\mathcal{L})\) is the concatenation \(u = u_1 u_2\) of \(u_1 \in \Sigma_m(\mathcal{L})\) and \(u_2 \in \Sigma_k(\mathcal{L})\) sharing one junction vertex \(j_0 \in Q_0\), then \(\underline{\dim}\, M(u_1 u_2) = \underline{\dim}\, M(u_1) + \underline{\dim}\, M(u_2) - e_{j_0}\), so
\[
\bigl\langle \theta,\, [M(u_1 u_2)] \bigr\rangle \;\le\; \bigl\langle \theta,\, [M(u_1)] \bigr\rangle \;+\; \bigl\langle \theta,\, [M(u_2)] \bigr\rangle \;+\; \|\theta\|_\infty.
\]
Consequently, \(Z_{m+k}^{\mathrm{full}}(\mathcal{L}; \theta) \le e^{\|\theta\|_\infty} Z_m^{\mathrm{full}}(\mathcal{L}; \theta) Z_k^{\mathrm{full}}(\mathcal{L}; \theta)\), so the sequence \(b_m \coloneqq \log Z_m^{\mathrm{full}}(\mathcal{L}; \theta) + \|\theta\|_\infty\) is subadditive (\(b_{m+k} \le b_m + b_k\)). By Fekete's subadditivity lemma, the limit \(P(\mathcal{L}; \theta) \coloneqq \lim_{m \to \infty} \frac{1}{m}\log Z_m^{\mathrm{full}}(\mathcal{L}; \theta)\) exists and is finite (bounded by \(h_{\mathrm{cod}}(\mathcal{L}) \pm \|\theta\|_\infty\)).

To compare \(Z_n(N; \theta)\) (where words \(w \in \mathcal{T}_n(N) = \Sigma_{2n}(\mathcal{L}; i)\) are anchored at vertex \(i\) at their midpoint) with \(Z_{2n}^{\mathrm{full}}(\mathcal{L}; \theta)\), observe first that \(\Sigma_{2n}(\mathcal{L}; i) \subseteq \Sigma_{2n}(\mathcal{L})\), so \(Z_n(N; \theta) \le Z_{2n}^{\mathrm{full}}(\mathcal{L}; \theta)\). Conversely, by the syllable-level syndeticity of \(X_i = \mathcal{L} \cap d_i \subseteq X\) established in the proof of Theorem~\ref{thm:entropy_equivalence}, there exists an integer \(K \ge 1\) such that every word \(u \in \Sigma_{2n}(\mathcal{L})\) extends on both ends along a leaf of \(\mathcal{L}\) to a word \(w \in \Sigma_{2(n+K)}(\mathcal{L}; i)\) centered at vertex \(i\) such that \(u\) appears as the length-\(2n\) subword of \(w\) at some syllable offset \(r \in \{-K, \dots, K\}\) from the midpoint of \(w\). Because for each \(w \in \Sigma_{2(n+K)}(\mathcal{L}; i)\) and each offset \(r \in \{-K, \dots, K\}\), the length-\(2n\) subword of \(w\) at offset \(r\) is unique, at most \(2K + 1\) distinct words \(u \in \Sigma_{2n}(\mathcal{L})\) can extend to the same \(w\) at an allowed offset. Combined with \(\bigl|\langle \theta, [M(w)] \rangle - \langle \theta, [M(u)] \rangle\bigr| \le 2K\|\theta\|_\infty\), we obtain
\[
Z_{2n}^{\mathrm{full}}(\mathcal{L}; \theta) \;\le\; (2K+1)\, e^{2K\|\theta\|_\infty} \, Z_{n+K}(N; \theta).
\]
Taking \(\frac{1}{2n}\log(\cdot)\), letting \(n \to \infty\), and using part (1) gives
\[
P^{\mathrm{shape}}_N(\theta) \;=\; \lim_{n \to \infty} \frac{1}{2n} \log Z_n(N; \theta) \;=\; \lim_{n \to \infty} \frac{1}{2n} \log Z_{2n}^{\mathrm{full}}(\mathcal{L}; \theta) \;=\; P(\mathcal{L}; \theta),
\]
which exists and depends only on \(\mathcal{L}\) (and not on the supporting vertex \(i \in Q_0\)). Convexity of \(P^{\mathrm{shape}}_N(\theta)\) follows from H\"older's inequality applied to \(Z_n^{\mathrm{shape}}(N; \theta)\), and \(1\)-Lipschitz continuity follows from \(\bigl|\langle \theta - \theta', [M(\mathfrak{s})] \rangle\bigr| \le (2n+1)\|\theta - \theta'\|_\infty\).
\end{proof}

\begin{remark}[Purely Algebraic Definition of \(\operatorname{Supp}_{K_0}^\infty(N)\), \(\mathfrak{h}_{K_0}(N; d)\), and \(P^{\mathrm{shape}}_N(\theta)\)]
\label{rem:algebraic_independence_pressure}
All three invariants in Definition~\ref{def:pp_pressure} are defined purely in terms of the module category \(\Mod{\Lambda}\), the localized pp-lattice \(\mathcal{L}_i(\Lambda)\), and the Grothendieck group \(K_0(\md{\Lambda})\), without reference to the surface \((\mathcal{S}, \mathcal{M}, \forb)\):
\begin{enumerate}
    \item The active radius-\(n\) interval pairs \(\mathcal{T}_n(N)\) and their \(N\)-equivalence classes \(\mathcal{S}_n(N) = \mathcal{T}_n(N)/\!\sim_N\) are determined solely by the lattice of pp-definable subspaces of \(N e_i\).
    \item By construction (Section~\ref{subsec:geometric_pp}), for each pointed string formula \(\phi_u(x)\) with \(u \in \operatorname{Str}_{2n}(i)\), the pointed string module \((M(u), v_0)\) satisfies \(v_0 \in \phi_u(M(u))\) and every \(m \in \phi_u(M)\) (\(M \in \Mod{\Lambda}\)) is of the form \(f(v_0)\) for some homomorphism \(f \colon M(u) \to M\). This is the defining universal property of a free realization of \(\phi_u(x)\) at the distinguished generator \(v_0\) \cite[\S 1.2]{Prest2009} within the class of pointed string modules. By Lemma~\ref{lem:witness_module_shape}, the unpointed isomorphism class \(M(\mathfrak{s}) \cong M(u)\) and Grothendieck class \([M(\mathfrak{s})] \in K_0(\md{\Lambda})\) depend only on the \(N\)-shape equivalence class \(\mathfrak{s} \in \mathcal{S}_n(N)\).
\end{enumerate}
\end{remark}

\begin{theorem}[Algebraic Ruelle--Walters Variational Principle]
\label{thm:algebraic_variational_principle}
Let \(N = \Psi(\mathcal{L}) \in \mathcal{P}_i\) be an arational pure-injective module associated with a minimal recurrent lamination \(\mathcal{L} \in \mathcal{P}_{\mathrm{geo}}\). Let \(X = \mathcal{L} \cap \bigsqcup_{j \in Q_0} d_j\) be the global dissection Cantor cross-section with first-return map \(T \colon X \to X\), and for each \([\mu] \in \mathcal{M}_{\mathrm{erg}}(\mathcal{L})\) (normalized by \(|\mu|_\Delta = \sum_{j \in Q_0} \mu(d_j) = 1\)), let \(\mathbb{P}_\mu\) denote the induced \(T\)-invariant ergodic Borel probability measure on \(X\) with asymptotic dimension vector \(\underline{d}(\mu) = \Phi_\Delta(\mu) \in \Sigma^{Q_0}\) (Theorem~\ref{thm:ergodic_classification}). Then for every weight \(\theta \in \Hom_{\Z}(K_0(\md{\Lambda}), \R) \cong \R^{Q_0}\):
\begin{equation}
\label{eq:algebraic_variational_principle}
P^{\mathrm{shape}}_N(\theta) \;=\; \sup_{[\mu] \in \mathcal{M}_{\mathrm{erg}}(\mathcal{L})} \Bigl(\, h_{\mathbb{P}_\mu}(T) \;+\; \bigl\langle \theta,\, \underline{d}(\mu) \bigr\rangle \,\Bigr).
\end{equation}
\end{theorem}

\begin{proof}
Let \(\Sigma(\mathcal{L}) \subset (Q_1 \cup Q_1^{-1})^{\Z}\) be the two-sided symbolic subshift of bi-infinite admissible words coding the leaves of \(\mathcal{L}\), equipped with the left shift \(\sigma \colon \Sigma(\mathcal{L}) \to \Sigma(\mathcal{L})\). Define the locally constant cylinder potential \(f_\theta \colon \Sigma(\mathcal{L}) \to \R\) by \(f_\theta(\overleftrightarrow{w}) \coloneqq \theta_{s(l_1)}\), where \(l_1\) is the \(0\)-th syllable of \(\overleftrightarrow{w}\). For any finite word \(u = l_m \dots l_1 \in \Sigma_m(\mathcal{L})\) of syllable length \(m \ge 1\) and any \(\overleftrightarrow{w} \in [u] \subset \Sigma(\mathcal{L})\), the \(m\)-step Birkhoff sum \((S_m f_\theta)(\overleftrightarrow{w}) = \sum_{k=0}^{m-1} f_\theta(\sigma^k \overleftrightarrow{w})\) counts the initial vertices \(s(l_1), \dots, s(l_m)\) of the \(m\) syllables of \(u\). By Proposition~\ref{prop:finite_level_identity}, the canonical basis of \(M(u)\) consists of these \(m\) vertices together with the terminal vertex \(t(l_m)\), so
\[
\Bigl|\, \bigl(S_m f_\theta\bigr)(\overleftrightarrow{w}) \;-\; \bigl\langle \theta,\, [M(u)] \bigr\rangle \,\Bigr| \;=\; \bigl|\theta_{t(l_m)}\bigr| \;\le\; \|\theta\|_\infty
\]
uniformly for all \(m \ge 1\) and \(u \in \Sigma_m(\mathcal{L})\). Consequently, the topological pressure of \(f_\theta\) with respect to \(\sigma\) satisfies
\[
\begin{aligned}
P_{\mathrm{top}}(\sigma, f_\theta) 
&= \lim_{m \to \infty} \frac{1}{m} \log \sum_{u \in \Sigma_m(\mathcal{L})} \sup_{\overleftrightarrow{w} \in [u]} \exp\bigl((S_m f_\theta)(\overleftrightarrow{w})\bigr) \\
&= \lim_{m \to \infty} \frac{1}{m} \log Z_m^{\mathrm{full}}(\mathcal{L}; \theta) \;=\; P^{\mathrm{shape}}_N(\theta)
\end{aligned}
\]
by Proposition~\ref{prop:weighted_shape_injectivity}(2). By the Ruelle--Walters Variational Principle \cite[Theorem~9.10]{Walters1982},
\[
P_{\mathrm{top}}(\sigma, f_\theta) \;=\; \sup_{\nu \in \mathcal{M}_{\mathrm{erg}}(\Sigma(\mathcal{L}), \sigma)} \left(\, h_\nu(\sigma) \;+\; \int_{\Sigma(\mathcal{L})} f_\theta \, d\nu \,\right).
\]
As shown in the proof of Theorem~\ref{thm:entropy_equivalence}, the symbolic coding map \((X, T) \to (\Sigma(\mathcal{L}), \sigma)\) is a measure-theoretic isomorphism for every non-atomic \(T\)-invariant Borel probability measure \(\mathbb{P}_\mu\) (\([\mu] \in \mathcal{M}_{\mathrm{erg}}(\mathcal{L})\)), intertwining \(T\) with \(\sigma\) and mapping the partition \(X = \bigsqcup_{j \in Q_0} (\mathcal{L} \cap d_j)\) to the \(0\)-th vertex partition of \(\Sigma(\mathcal{L})\). Under this identification, \(h_\nu(\sigma) = h_{\mathbb{P}_\mu}(T)\) and \(\int_{\Sigma(\mathcal{L})} f_\theta \, d\nu = \sum_{j \in Q_0} \theta_j \, \mu(d_j) = \langle \theta, \underline{d}(\mu) \rangle\), yielding \eqref{eq:algebraic_variational_principle}.
\end{proof}

\begin{theorem}[Legendre--Fenchel Duality, Asymptotic Grothendieck Support, and Microcanonical \(K_0\)-Entropy]
\label{thm:legendre_duality}
Let \(N = \Psi(\mathcal{L}) \in \mathcal{P}_i\) be an arational pure-injective module, and let
\[
\mathcal{D}(N) \;\coloneqq\; \operatorname{conv}\Bigl\{\, \underline{d}(\mu) \;\Big|\; [\mu] \in \mathcal{M}_{\mathrm{erg}}(\mathcal{L}) \,\Bigr\} \;\subset\; \Sigma^{Q_0} \subset \R^{Q_0}
\]
denote the compact convex simplex of asymptotic dimension vectors of \(N\) (Section~\ref{sec:measured}). Then:
\begin{enumerate}
    \item \textbf{Identity of the Asymptotic Grothendieck Support, Recession Subdifferential, and Ergodic Simplex:} For every \(\theta \in \R^{Q_0}\), the recession limit \\ \(\sigma_{\mathcal{D}(N)}(\theta) \coloneqq \lim_{\beta \to +\infty} \beta^{-1} P^{\mathrm{shape}}_N(\beta\theta)\) equals the support function \(\max_{d \in \mathcal{D}(N)} \langle \theta, d \rangle\), and its subdifferential \(\partial\sigma_{\mathcal{D}(N)}(0)\) equals the closed convex hull of the intrinsic algebraic Grothendieck support \(\operatorname{Supp}_{K_0}^\infty(N)\) (Definition~\ref{def:pp_pressure}(1)):
    \[
    \begin{gathered}
    \operatorname{conv}\bigl(\operatorname{Supp}_{K_0}^\infty(N)\bigr)
      \;=\; \partial\sigma_{\mathcal{D}(N)}(0) \;=\; \mathcal{D}(N), \\[2pt]
    \dim_{\mathrm{aff}}\bigl(\mathcal{D}(N)\bigr)
      \;=\; \operatorname{rk}^\infty_{K_0}(N) \;-\; 1.
    \end{gathered}
    \]
    Whenever the cellular dissection \(\Delta\) restricts to an ideal triangulation on the lamination-supporting subsurface of \(\mathcal{S}^\circ \setminus \forb\) (in particular, for every gentle algebra; see Remark~\ref{rem:when_which_embedding}), the topological dimension of the Choquet simplex \(\mathcal{M}_1(\mathcal{L}) = \overline{\operatorname{conv}}\,\mathcal{M}_{\mathrm{erg}}(\mathcal{L})\) is given purely algebraically by \(\dim_{\mathrm{top}}(\mathcal{M}_1(\mathcal{L})) = \operatorname{rk}^\infty_{K_0}(N) - 1\); under this same triangulation condition, \(\mathcal{L}\) is uniquely ergodic if and only if \(\operatorname{rk}^\infty_{K_0}(N) = 1\).
    \item \textbf{Subdifferential at \(\theta = 0\) and Maximal-Entropy Measures:} The subdifferential \(\partial P^{\mathrm{shape}}_N(0) \coloneqq \{ d \in \R^{Q_0} \mid P^{\mathrm{shape}}_N(\theta) \ge h^{\mathrm{shape}}(N) + \langle \theta, d \rangle \;\forall \theta \in \R^{Q_0} \}\) equals the convex hull of the asymptotic dimension vectors of the measures of maximal entropy on \(\mathcal{L}\):
    \[
    \partial P^{\mathrm{shape}}_N(0) \;=\; \operatorname{conv}\Bigl\{\, \underline{d}(\mu) \;\Big|\; [\mu] \in \mathcal{M}_{\mathrm{erg}}(\mathcal{L}),\; h_{\mathbb{P}_\mu}(T) = h^{\mathrm{shape}}(N) \,\Bigr\}.
    \]
    In particular, if \(\mathcal{L}\) is uniquely ergodic with unique normalized measure \([\mu]\), then \(P^{\mathrm{shape}}_N\) is differentiable at \(\theta = 0\) and \(\nabla P^{\mathrm{shape}}_N(0) = \underline{d}(\mu) = \Phi_\Delta(\mu)\).
    \item \textbf{Legendre--Fenchel Transform and Microcanonical \(K_0\)-Shape Entropy:} For every \(d \in \Sigma^{Q_0}\), the microcanonical \(K_0\)-shape entropy \(\mathfrak{h}_{K_0}(N; d)\) (Definition~\ref{def:pp_pressure}(2)) satisfies the sandwich identity
    \[
    \sup_{\substack{[\mu] \in \mathcal{M}_{\mathrm{erg}}(\mathcal{L}) \\ \underline{d}(\mu) = d}} h_{\mathbb{P}_\mu}(T)
    \;\le\;
    \mathfrak{h}_{K_0}(N; d)
    \;\le\;
    \sup_{\substack{\mu \in \mathcal{M}_1(\mathcal{L}) \\ \underline{d}(\mu) = d}} h_{\mathbb{P}_\mu}(T)
    \;=\;
    \inf_{\theta \in \R^{Q_0}} \Bigl(\, P^{\mathrm{shape}}_N(\theta) \;-\; \bigl\langle \theta,\, d \bigr\rangle \,\Bigr)
    \]
    (with \(\mathfrak{h}_{K_0}(N; d) = -\infty\) for all \(d \notin \mathcal{D}(N)\)). In particular, at every extreme point \(d \in \operatorname{ext}(\mathcal{D}(N))\) (and globally on \(\mathcal{D}(N)\) whenever \(\mathcal{L}\) is uniquely ergodic), the Legendre--Fenchel transform of \(P^{\mathrm{shape}}_N\) computes \(\mathfrak{h}_{K_0}(N; d)\) identically:
    \[
    \mathfrak{h}_{K_0}(N; d) \;=\; \inf_{\theta \in \R^{Q_0}} \Bigl(\, P^{\mathrm{shape}}_N(\theta) \;-\; \bigl\langle \theta,\, d \bigr\rangle \,\Bigr).
    \]
\end{enumerate}
\end{theorem}

\begin{proof}
(1) Replacing \(\theta\) by \(\beta\theta\) (\(\beta > 0\)) in Theorem~\ref{thm:algebraic_variational_principle} and using \(0 \le h_{\mathbb{P}_\mu}(T) \le h^{\mathrm{shape}}(N) < \infty\) gives
\[
\begin{aligned}
\beta \sup_{[\mu] \in \mathcal{M}_{\mathrm{erg}}(\mathcal{L})} \bigl\langle \theta,\, \underline{d}(\mu) \bigr\rangle 
  &\;\le\; P^{\mathrm{shape}}_N(\beta\theta) \\
  &\;\le\; h^{\mathrm{shape}}(N) \;+\; \beta \sup_{[\mu] \in \mathcal{M}_{\mathrm{erg}}(\mathcal{L})} \bigl\langle \theta,\, \underline{d}(\mu) \bigr\rangle.
\end{aligned}
\]
Dividing by \(\beta\) and letting \(\beta \to +\infty\) yields \(\sigma_{\mathcal{D}(N)}(\theta) = \max_{d \in \mathcal{D}(N)} \langle \theta, d \rangle\), so \(\partial\sigma_{\mathcal{D}(N)}(0) = \mathcal{D}(N)\) by Hahn--Banach separation. To prove \(\operatorname{conv}(\operatorname{Supp}_{K_0}^\infty(N)) = \mathcal{D}(N)\), let \([\mu] \in \mathcal{M}_{\mathrm{erg}}(\mathcal{L})\). By Theorem~\ref{thm:ergodic_classification}, for any generic leaf \(\ell \in \mathcal{L}_\mu\), its centered radius-\(n\) prefixes \(w_n \in \mathcal{T}_n(N)\) satisfy \(\lim_{n \to \infty} \frac{[M(\mathfrak{s}_{w_n})]}{2n+1} = \underline{d}(\mu)\), so \(\underline{d}(\mu) \in \operatorname{Supp}_{K_0}^\infty(N)\) and hence \(\mathcal{D}(N) \subseteq \operatorname{conv}(\operatorname{Supp}_{K_0}^\infty(N))\). Conversely, if \(d \in \operatorname{Supp}_{K_0}^\infty(N)\) is the limit of \(\frac{[M(\mathfrak{s}_k)]}{2n_k+1}\) for \(\mathfrak{s}_k = [w_{n_k}]_N \in \mathcal{S}_{n_k}(N)\) (\(n_k \to \infty\)), then by the Finite-Level Identity (Proposition~\ref{prop:finite_level_identity}) applied to the empirical measures \(\frac{1}{2n_k}\sum_{r=0}^{2n_k-1}\delta_{\sigma^r(\overleftrightarrow{w}_k)}\) and the weak-\(*\) compactness of \(\mathcal{M}_1(\Sigma(\mathcal{L})) \cong \mathcal{M}_1(\mathcal{L})\), any weak-\(*\) cluster point \(\nu_\infty \in \mathcal{M}_1(\mathcal{L})\) is shift-invariant with \(\underline{d}(\nu_\infty) = d\); since \(\mathcal{M}_1(\mathcal{L}) = \overline{\operatorname{conv}}\,\mathcal{M}_{\mathrm{erg}}(\mathcal{L})\) by Choquet's theorem, \(d \in \mathcal{D}(N)\). Thus \(\operatorname{Supp}_{K_0}^\infty(N) \subseteq \mathcal{D}(N)\) and \(\operatorname{conv}(\operatorname{Supp}_{K_0}^\infty(N)) = \mathcal{D}(N)\). Because taking the convex hull preserves the real linear span, \(\operatorname{span}_\R(\operatorname{Supp}_{K_0}^\infty(N)) = \operatorname{span}_\R(\mathcal{D}(N))\); since \(\mathcal{D}(N)\) lies in the affine hyperplane \(\sum_{j \in Q_0} d_j = 1\) not passing through \(0\), \(\dim_{\mathrm{aff}}(\mathcal{D}(N)) = \dim_\R \operatorname{span}_\R(\mathcal{D}(N)) - 1 = \operatorname{rk}^\infty_{K_0}(N) - 1\). When \(\Delta\) restricts to an ideal triangulation on \(\operatorname{supp}(\mathcal{L})\), Remark~\ref{rem:when_which_embedding} ensures that \(\Phi_\Delta \colon \mathcal{M}_1(\mathcal{L}) \to \mathcal{D}(N)\) is an affine homeomorphism, giving \(\dim_{\mathrm{top}}(\mathcal{M}_1(\mathcal{L})) = \operatorname{rk}^\infty_{K_0}(N) - 1\).

(2) Because the symbolic subshift \((\Sigma(\mathcal{L}), \sigma)\) is expansive, the Kolmogorov--Sinai entropy map \(\mu \mapsto h_{\mathbb{P}_\mu}(T)\) is affine and upper semi-continuous on \(\mathcal{M}_1(\mathcal{L})\) \cite[Theorem~8.2]{Walters1982}. By Theorem~\ref{thm:algebraic_variational_principle}, \(P^{\mathrm{shape}}_N(\theta) = \max_{\mu \in \mathcal{M}_1(\mathcal{L})} \bigl(h_{\mathbb{P}_\mu}(T) + \langle \theta, \underline{d}(\mu) \rangle\bigr)\), so \(d \in \partial P^{\mathrm{shape}}_N(0)\) if and only if \(d = \underline{d}(\mu)\) for some \(\mu \in \mathcal{M}_1(\mathcal{L})\) with \(h_{\mathbb{P}_\mu}(T) = h^{\mathrm{shape}}(N)\) \cite[Theorem~9.12]{Walters1982}. When \(\mathcal{L}\) is uniquely ergodic, \(\partial P^{\mathrm{shape}}_N(0) = \{\underline{d}(\mu)\}\) is a singleton, so \(\nabla P^{\mathrm{shape}}_N(0) = \underline{d}(\mu)\) \cite[Theorem~25.1]{Rockafellar1970}.

(3) The image \(\bigl\{\bigl(\underline{d}(\mu), h_{\mathbb{P}_\mu}(T)\bigr)\bigr\}_{\mu \in \mathcal{M}_1(\mathcal{L})} \subset \R^{Q_0} \times \R\) of the compact convex set \(\mathcal{M}_1(\mathcal{L})\) under the continuous affine map \(\mu \mapsto (\underline{d}(\mu), h_{\mathbb{P}_\mu}(T))\) (with upper semi-continuous second coordinate) has a compact convex hypograph, so its upper boundary \(H(d) \coloneqq \sup\{h_{\mathbb{P}_\mu}(T) \mid \mu \in \mathcal{M}_1(\mathcal{L}),\; \underline{d}(\mu) = d\}\) (with \(H(d) = -\infty\) for \(d \notin \mathcal{D}(N)\)) is concave and upper semi-continuous on \(\R^{Q_0}\). By Theorem~\ref{thm:algebraic_variational_principle} and Fenchel--Moreau biconjugation \cite[Theorem~12.2]{Rockafellar1970}, \(H(d) = \inf_{\theta \in \R^{Q_0}} \bigl(P^{\mathrm{shape}}_N(\theta) - \langle \theta, d \rangle\bigr)\). For any \(\theta \in \R^{Q_0}\) and \(\mathfrak{s} \in \mathcal{S}_n(N)\) with \(\bigl\|\frac{[M(\mathfrak{s})]}{2n+1} - d\bigr\|_\infty \le \varepsilon\), we have \(\langle \theta, [M(\mathfrak{s})] \rangle \ge (2n+1)(\langle \theta, d \rangle - \|\theta\|_1 \varepsilon)\), so
\[
Z_n^{\mathrm{shape}}(N; \theta) \;\ge\; S_n(N; d, \varepsilon) \exp\bigl((2n+1)(\langle \theta, d \rangle - \|\theta\|_1 \varepsilon)\bigr).
\]
Taking \(\frac{1}{2n}\log(\cdot)\), letting \(n \to \infty\) and \(\varepsilon \downarrow 0\), and optimizing over \(\theta \in \R^{Q_0}\) gives \(\mathfrak{h}_{K_0}(N; d) \le \inf_\theta (P^{\mathrm{shape}}_N(\theta) - \langle \theta, d \rangle) = H(d)\). Conversely, for any ergodic measure \([\mu] \in \mathcal{M}_{\mathrm{erg}}(\mathcal{L})\) with \(\underline{d}(\mu) = d\) and any \(\varepsilon > 0\), Birkhoff's Ergodic Theorem and the Shannon--McMillan--Breiman Theorem \cite[Theorem~9.3]{Walters1982} (combined with the syndeticity bound of Proposition~\ref{prop:weighted_shape_injectivity}) imply that the set of radius-\(n\) classes \(\mathfrak{s} \in \mathcal{S}_n(N)\) satisfying \(\bigl\|\frac{[M(\mathfrak{s})]}{2n+1} - d\bigr\|_\infty \le \varepsilon\) carries \(\mathbb{P}_\mu\)-measure \(\to 1\) and therefore has cardinality \(S_n(N; d, \varepsilon) \ge \exp(2n(h_{\mathbb{P}_\mu}(T) - o(1)))\), proving \(\mathfrak{h}_{K_0}(N; d) \ge h_{\mathbb{P}_\mu}(T)\). At any extreme point \(d \in \operatorname{ext}(\mathcal{D}(N))\), the fiber \(\underline{d}^{-1}(d) = \{\mu \in \mathcal{M}_1(\mathcal{L}) \mid \underline{d}(\mu) = d\}\) is a compact convex extremal face of \(\mathcal{M}_1(\mathcal{L})\); by Bauer's Maximum Principle, the upper semi-continuous affine function \(\mu \mapsto h_{\mathbb{P}_\mu}(T)\) attains its maximum on \(\underline{d}^{-1}(d)\) at an extreme point of \(\underline{d}^{-1}(d)\), which is an extreme point of \(\mathcal{M}_1(\mathcal{L})\) and hence ergodic, so the upper and lower bounds coincide.
\end{proof}

Finally, we prove that \(P^{\mathrm{shape}}_N(\theta)\) strictly separates uncountable families of arational pure-injective modules sharing the same shape entropy \(h^{\mathrm{shape}}\). Throughout Theorem~\ref{thm:pressure_separation}, \(v_1\) is based at vertex \(i_1 = i = s(v_1) = t(v_1)\), the bridge \(u_{12}\) runs from \(i_1\) to \(i_2 \coloneqq s(v_2) = t(v_2) \in Q_0\), and \(u_{21}\) runs from \(i_2\) back to \(i_1\); the normalized cyclic dimension vector of \(v_r\) (\(r = 1, 2\)) with respect to its base vertex \(i_r\) is defined by \(d^{\mathrm{cyc}}(v_r)_j \coloneqq |v_r|^{-1}(\dim_\K M(v_r)e_j - \delta_{j, i_r}) \in \Sigma^{Q_0}\).

\begin{theorem}[Strict Separation of Iso-Entropic Pure-Injective Modules by the Pressure Function]
\label{thm:pressure_separation}
Let \(\Lambda\) be a non-domestic string algebra possessing two primitive bands \(v_1, v_2\) communicating in both directions via bridge strings \(u_{12}, u_{21}\) as above, such that \(d^{\mathrm{cyc}}(v_1) \neq d^{\mathrm{cyc}}(v_2)\) in \(\R^{Q_0}\). Define the bridge-completed equal-length composite blocks
\[
A \;\coloneqq\; v_1^{2|v_2|} \, u_{12} \, v_2^{|v_1|} \, u_{21}, \qquad B \;\coloneqq\; v_1^{|v_2|} \, u_{12} \, v_2^{2|v_1|} \, u_{21}
\]
of common syllable length \(L \coloneqq 3|v_1||v_2| + |u_{12}| + |u_{21}|\). Then every concatenation in \(\{A, B\}^{\Z}\) is an admissible bi-infinite string, \(d^{\mathrm{cyc}}(A) - d^{\mathrm{cyc}}(B) = \frac{|v_1||v_2|}{L}(d^{\mathrm{cyc}}(v_1) - d^{\mathrm{cyc}}(v_2)) \neq 0\), and:
\begin{enumerate}
    \item \textbf{Uncountable Zero-Entropy Sturmian Separation:} For each irrational slope \(\alpha \in (0, 1) \setminus \Q\), let \(\Omega_\alpha \subset \{0, 1\}^{\Z}\) be the minimal Sturmian subshift of slope \(\alpha\), and let \(N_\alpha \coloneqq \PE(M(\Theta(\overleftrightarrow{\omega}_\alpha))) \in \mathcal{P}(\Zg(\Lambda))\) be the pure-injective module obtained via the block substitution \(\Theta \colon 0 \mapsto A,\, 1 \mapsto B\). Then every \(N_\alpha\) has identical shape entropy
    \[
    h^{\mathrm{shape}}(N_\alpha) \;=\; 0 \qquad \text{for all } \alpha \in (0, 1) \setminus \Q,
    \]
    whereas for every \(\theta \in \R^{Q_0}\),
    \[
    P^{\mathrm{shape}}_{N_\alpha}(\theta) \;=\; \Bigl\langle \theta,\, (1-\alpha)\,d^{\mathrm{cyc}}(A) + \alpha\,d^{\mathrm{cyc}}(B) \Bigr\rangle,
    \qquad
    \nabla P^{\mathrm{shape}}_{N_\alpha}(0) \;=\; (1-\alpha)\,d^{\mathrm{cyc}}(A) + \alpha\,d^{\mathrm{cyc}}(B).
    \]
    In particular, \(\alpha \mapsto P^{\mathrm{shape}}_{N_\alpha}\) is strictly injective on \((0, 1) \setminus \Q\), and \(\{N_\alpha\}_{\alpha \in (0,1)\setminus\Q}\) is an uncountable family of pairwise non-isomorphic pure-injective modules in \(\mathcal{P}(\Zg(\Lambda))\) indistinguishable by \(h^{\mathrm{shape}}\) that are completely separated by \(P^{\mathrm{shape}}_N\).
    \item \textbf{Positive-Entropy Separation:} Let \(Y \subset \{0, 1\}^{\Z}\) be any minimal subshift of positive topological entropy \(h_{\mathrm{top}}(Y) > 0\), and let \(\Omega \subset \{0, 1\}^{\Z}\) be the minimal subshift generated by applying the recognizable constant-length-\(4\) substitution \(\psi \colon 0 \mapsto 0010,\, 1 \mapsto 0001\) to \(Y\) (or more generally, any minimal subshift of positive entropy \(h_0 = h_{\mathrm{top}}(\Omega) > 0\) satisfying \(\sup_{\nu \in \mathcal{M}_1(\Omega)} \nu([1]) < 1/2\)), with bit-flipped minimal subshift \(\overline{\Omega}\) (\(0 \leftrightarrow 1\)). Then \(N_\Omega, N_{\overline{\Omega}} \in \mathcal{P}(\Zg(\Lambda))\) have identical positive shape entropy
    \[
    h^{\mathrm{shape}}(N_\Omega) \;=\; \frac{h_0}{L} \;=\; h^{\mathrm{shape}}(N_{\overline{\Omega}}) \;>\; 0,
    \]
    but disjoint ergodic dimension simplices \(\mathcal{D}(N_\Omega) \cap \mathcal{D}(N_{\overline{\Omega}}) = \emptyset\) and strictly distinct pressure functions \(P^{\mathrm{shape}}_{N_\Omega}(\theta) \neq P^{\mathrm{shape}}_{N_{\overline{\Omega}}}(\theta)\).
\end{enumerate}
\end{theorem}

\begin{proof}
Because \(v_1, v_2\) communicate via \(u_{12}, u_{21}\), the cyclic word \(v_1^p u_{12} v_2^q u_{21}\) is admissible for all \(p, q \ge 1\). Since both \(A\) and \(B\) start with the initial syllable of \(v_1\) and end with the terminal syllable of \(u_{21}\), and the junction \(u_{21} v_1\) is admissible, every bi-infinite concatenation in \(\{A, B\}^{\Z}\) is an admissible string avoiding \(\rho \cup \rho^{-1}\). Moreover, \(|A| = 2|v_2||v_1| + |u_{12}| + |v_1||v_2| + |u_{21}| = 3|v_1||v_2| + |u_{12}| + |u_{21}| = |B| = L\). Along one period of \(A\) (resp.\ \(B\)), the bridge words \(u_{12}\) and \(u_{21}\) occur once each, so their contributions cancel identically in \(d^{\mathrm{cyc}}(A) - d^{\mathrm{cyc}}(B)\), leaving
\[
\begin{aligned}
d^{\mathrm{cyc}}(A) - d^{\mathrm{cyc}}(B) 
&= \frac{(2|v_2| - |v_2|)|v_1|\,d^{\mathrm{cyc}}(v_1) + (|v_1| - 2|v_1|)|v_2|\,d^{\mathrm{cyc}}(v_2)}{L} \\
&= \frac{|v_1||v_2|}{L}\bigl(d^{\mathrm{cyc}}(v_1) - d^{\mathrm{cyc}}(v_2)\bigr) \;\neq\; 0
\end{aligned}
\]
at every vertex \(j \in Q_0\). Because \(v_1, v_2\) are distinct primitive bands, \(A\) and \(B\) are distinct words of equal length \(L\), so \(\Theta \colon \{0, 1\}^{\Z} \to \{A, B\}^{\Z}\) is a constant-length-\(L\) injective block code. For any minimal subshift \(\Omega \subset \{0, 1\}^{\Z}\), the syllable-level shift closure \(\Sigma_\Omega \coloneqq \bigcup_{r=0}^{L-1} \sigma^r(\Theta(\Omega))\) is a minimal aperiodic subshift in \(\mathcal{P}_{\mathrm{geo}}\) with syllable coding entropy \(h_{\mathrm{cod}}(\Sigma_\Omega) = \frac{1}{L} h_{\mathrm{top}}(\Omega)\), and each ergodic measure \(\nu \in \mathcal{M}_{\mathrm{erg}}(\Omega)\) with \(1\)-frequency \(\alpha = \nu([1]) \in [0, 1]\) induces an ergodic measure \(\mu_\nu\) on \(\Sigma_\Omega\) whose normalized asymptotic dimension vector at every vertex \(j \in Q_0\) is
\[
\underline{d}(\mu_\nu) \;=\; (1-\alpha)\,d^{\mathrm{cyc}}(A) \;+\; \alpha\,d^{\mathrm{cyc}}(B).
\]
(1) When \(\Omega = \Omega_\alpha\) is the Sturmian minimal subshift of irrational slope \(\alpha \in (0, 1) \setminus \Q\), \(\Omega_\alpha\) is uniquely ergodic with \(\nu([1]) = \alpha\) and \(h_{\mathrm{top}}(\Omega_\alpha) = 0\). By Theorems~\ref{thm:entropy_equivalence}, \ref{thm:algebraic_variational_principle}, and \ref{thm:legendre_duality}(2), \(h^{\mathrm{shape}}(N_\alpha) = 0\) and \(P^{\mathrm{shape}}_{N_\alpha}(\theta) = \langle \theta, (1-\alpha)d^{\mathrm{cyc}}(A) + \alpha d^{\mathrm{cyc}}(B) \rangle\). Since \(d^{\mathrm{cyc}}(A) \neq d^{\mathrm{cyc}}(B)\), the affine map \(\alpha \mapsto (1-\alpha)d^{\mathrm{cyc}}(A) + \alpha d^{\mathrm{cyc}}(B)\) is strictly injective on \((0, 1) \setminus \Q\), so \(P^{\mathrm{shape}}_{N_\alpha} \neq P^{\mathrm{shape}}_{N_{\alpha'}}\) (and hence \(N_\alpha \not\cong N_{\alpha'}\)) for all \(\alpha \neq \alpha'\).

(2) Under the recognizable constant-length-\(4\) substitution \(\psi \colon 0 \mapsto 0010,\, 1 \mapsto 0001\), every \(4\)-block contains three \(0\)s and one \(1\), so for any minimal subshift \(Y \subset \{0,1\}^{\Z}\) with \(h_{\mathrm{top}}(Y) > 0\), its shift-orbit closure \(\Omega\) is minimal with \(h_0 = h_{\mathrm{top}}(\Omega) = \frac{1}{4}h_{\mathrm{top}}(Y) > 0\) and satisfies \(\nu([1]) = \frac{1}{4} < \frac{1}{2}\) for every \(\nu \in \mathcal{M}_1(\Omega)\). More generally, for any minimal \(\Omega\) with \(\sup_{\nu \in \mathcal{M}_1(\Omega)} \nu([1]) = \alpha_+ < 1/2\), the bit-flip involution \(0 \leftrightarrow 1\) is a topological conjugacy from \(\Omega\) to \(\overline{\Omega}\), so \(h_{\mathrm{top}}(\overline{\Omega}) = h_{\mathrm{top}}(\Omega) = h_0\) and \(h^{\mathrm{shape}}(N_\Omega) = h_0/L = h^{\mathrm{shape}}(N_{\overline{\Omega}}) > 0\), while every \(\bar{\nu} \in \mathcal{M}_1(\overline{\Omega})\) satisfies \(\bar{\nu}([1]) = 1 - \nu([1]) \ge 1 - \alpha_+ > 1/2\). Since \(d^{\mathrm{cyc}}(A) \neq d^{\mathrm{cyc}}(B)\) exactly, \(\mathcal{D}(N_\Omega)\) and \(\mathcal{D}(N_{\overline{\Omega}})\) lie on strictly opposite sides of the midpoint \(\frac{1}{2}(d^{\mathrm{cyc}}(A) + d^{\mathrm{cyc}}(B))\) along the line segment \([d^{\mathrm{cyc}}(A), d^{\mathrm{cyc}}(B)]\), so \(\mathcal{D}(N_\Omega) \cap \mathcal{D}(N_{\overline{\Omega}}) = \emptyset\). By Theorem~\ref{thm:legendre_duality}(1), \(P^{\mathrm{shape}}_{N_\Omega} \neq P^{\mathrm{shape}}_{N_{\overline{\Omega}}}\).
\end{proof}

\subsection{The Algebraic pp-Bifurcation Index \texorpdfstring{\(\kappa_{\mathrm{pp}}(N)\)}{kappa\_pp(N)}, Gauss--Bonnet Simplicity Obstruction, and Iso-Pressure Separation}
\label{subsec:pp_bifurcation_index}

While the pressure function \(P^{\mathrm{shape}}_N(\theta)\) reconstructs the ergodic dimension simplex \(\mathcal{D}(N)\) and separates iso-entropic modules with distinct dimension simplices (Theorem~\ref{thm:pressure_separation}), it is defined via an exponential limit \(\lim_{n \to \infty} \frac{1}{2n}\log(\cdot)\). Consequently, for every zero-entropy arational pure-injective module \(N\) (including all simple minimal arational surface laminations \(\mathcal{L} \subset \mathcal{S}^\circ \setminus \forb\) and all Sturmian modules \(N_\alpha\)), \(P^{\mathrm{shape}}_N(\theta) = \max_{\mu \in \mathcal{M}_1(\mathcal{L})} \langle \theta, \underline{d}(\mu) \rangle\) depends solely on the first-order syllable frequencies \(\underline{d}(\mu)\) and is blind to linear and polynomial complexity of the localized pp-lattice \(\mathcal{L}_i(\Lambda)\).

We now introduce a finer algebraic invariant of \(\mathcal{L}_i(\Lambda)\)---the \emph{pp-bifurcation index} \(\kappa_{\mathrm{pp}}(N)\)---that counts one-step interval splittings in \(\mathcal{L}_i(\Lambda)\), provides a purely algebraic Gauss--Bonnet obstruction to whether the surface curve \(\Phi(N)\) is simple (non-self-intersecting), and strictly separates infinite hierarchies of arational pure-injective modules whose pressure functions coincide identically on all of \(\R^{Q_0}\).

Fix a vertex \(i \in Q_0\), a walk direction \(\epsilon \in \{-, +\}\), and \(n \ge 1\). Let
\[
\mathcal{U}_n^\epsilon(N; i) \;\coloneqq\; \Bigl\{\, v \in \operatorname{Str}_n^\epsilon(i) \;\Big|\; \frac{\phi_{v^\uparrow}(N)}{\phi_{v^\sharp}(N)} \neq 0 \,\Bigr\}
\]
denote the set of active one-sided length-\(n\) templates of \(N\) at \((i, \epsilon)\) (denoted \(\mathcal{U}_n^\epsilon\) to distinguish them from the two-sided templates \(\mathcal{T}_n(N)\) of Definition~\ref{def:shape_entropy}). For any \(v \in \operatorname{Str}_n^\epsilon(i)\) terminating at vertex \(j = t(v)\), the last syllable of \(v\) enters the ribbon strip \(R_j\) on one half-side of \(d_j\), so any non-backtracking extension syllable must exit \(R_j\) on the opposite half-side of \(d_j\). By the string algebra axioms (Definition~\ref{def:string_algebra}(1)--(2)), that half-side contains at most one direct arrow \(\alpha_v \in Q_1\) (\(s(\alpha_v) = j\), \(v\alpha_v \in \operatorname{Str}_{n+1}^\epsilon(i)\)) and at most one inverse arrow \(\beta_v^{-1} \in Q_1^{-1}\) (\(t(\beta_v) = j\), \(v\beta_v^{-1} \in \operatorname{Str}_{n+1}^\epsilon(i)\)). Under the covariant Prest--Puninski order \(\le_{\mathrm{str}}\) on \(\operatorname{Str}_{n+1}^\epsilon(i)\), when both extensions exist they are adjacent and satisfy
\[
(v\beta_v^{-1})^\uparrow \;=\; v^\uparrow, \qquad (v\alpha_v)^\uparrow \;=\; (v\beta_v^{-1})^\sharp, \qquad (v\alpha_v)^\sharp \;=\; v^\sharp.
\]
Hence the one-sided pp-interval \([\phi_{v^\sharp}(N), \phi_{v^\uparrow}(N)]\) in \(\mathcal{L}_i(\Lambda)\) is filtered at step \(n+1\) by the single canonical intermediate pp-definable subspace \(\phi_{(v\alpha_v)^\uparrow}(N)\):
\begin{equation}
\label{eq:one_step_pp_filtration}
\phi_{v^\sharp}(N) \;=\; \phi_{(v\alpha_v)^\sharp}(N) \;\subseteq\; \phi_{(v\alpha_v)^\uparrow}(N) \;=\; \phi_{(v\beta_v^{-1})^\sharp}(N) \;\subseteq\; \phi_{(v\beta_v^{-1})^\uparrow}(N) \;=\; \phi_{v^\uparrow}(N).
\end{equation}

\begin{definition}[One-Step pp-Forks and the Algebraic pp-Bifurcation Index]
\label{def:pp_bifurcation_index}
Let \(N \in \mathcal{P}(\Zg(\Lambda))\) be an arational pure-injective module.
\begin{enumerate}
    \item For each \(i \in Q_0\), \(\epsilon \in \{-, +\}\), and \(n \ge 1\), the set of \emph{depth-\(n\) pp-forks} of \(N\) at \((i, \epsilon)\) is the set of active one-sided templates \(v \in \mathcal{U}_n^\epsilon(N; i)\) for which both inclusions in \eqref{eq:one_step_pp_filtration} are strictly proper on \(N e_i\):
    \[
    \operatorname{Fork}_n^\epsilon(N; i) \;\coloneqq\; \Bigl\{\, v \in \mathcal{U}_n^\epsilon(N; i) \;\Big|\; \phi_{v^\sharp}(N) \;\subsetneq\; \phi_{(v\alpha_v)^\uparrow}(N) \;\subsetneq\; \phi_{v^\uparrow}(N) \,\Bigr\}.
    \]
    \item The \emph{algebraic pp-bifurcation index} of \(N\) is the asymptotic total forking rate across all sorts \(i \in Q_0\):
    \[
    \kappa_{\mathrm{pp}}(N) \;\coloneqq\; \limsup_{n \to \infty} \frac{1}{n} \sum_{i \in Q_0} \sum_{k=1}^{n-1} \Bigl(\bigl|\operatorname{Fork}_k^+(N; i)\bigr| \;+\; \bigl|\operatorname{Fork}_k^-(N; i)\bigr|\Bigr) \;\in\; [0, \infty].
    \]
\end{enumerate}
\end{definition}

\begin{theorem}[Exact Bifurcation Identity, Gauss--Bonnet Spike Theorem, and Iso-Pressure Separation]
\label{thm:bifurcation_index}
Let \(\Lambda\) be a non-domestic string algebra modeled on \((\mathcal{S}, \mathcal{M}, \forb)\), and let \(N = \Psi(\mathcal{L}) \in \mathcal{P}(\Zg(\Lambda))\) be an arational pure-injective module. Then:
\begin{enumerate}
    \item \textbf{Exact Algebraic Bifurcation Identity:} For every vertex \(i \in Q_0\), walk direction \(\epsilon \in \{-, +\}\), and depth \(n \ge 1\),
    \[
    \bigl|\mathcal{U}_{n+1}^\epsilon(N; i)\bigr| \;-\; \bigl|\mathcal{U}_n^\epsilon(N; i)\bigr| \;=\; \bigl|\operatorname{Fork}_n^\epsilon(N; i)\bigr|,
    \qquad
    \bigl|\mathcal{U}_n^\epsilon(N; i)\bigr| \;=\; \bigl|\mathcal{U}_1^\epsilon(N; i)\bigr| \;+\; \sum_{k=1}^{n-1} \bigl|\operatorname{Fork}_k^\epsilon(N; i)\bigr|.
    \]
    \item \textbf{Gauss--Bonnet Spike Theorem for Simple Geodesic Laminations:} Suppose \(\mathcal{L} \subset \mathcal{S}^\circ \setminus \forb\) is a \emph{simple} (non-self-intersecting) minimal arational geodesic lamination. Then for every \(n \ge 1\) and \(\epsilon \in \{-, +\}\), the depth-\(n\) pp-forks \(\bigsqcup_{i \in Q_0} \operatorname{Fork}_n^\epsilon(N; i)\) inject into the set of ideal boundary spikes of the complementary ideal polygons of \(\mathcal{S}^\circ \setminus \mathcal{L}\). Consequently, for all \(n \ge 1\):
    \[
    1 \;\le\; \sum_{i \in Q_0} \bigl|\operatorname{Fork}_n^\epsilon(N; i)\bigr| \;\le\; \frac{3}{\pi}\operatorname{Area}(\mathcal{S}^\circ) \;=\; 3\bigl(|\mathcal{M}| - 2\chi(\mathcal{S})\bigr),\footnote{Each complementary region $R_r$ of $\mathcal{S}^\circ \setminus \mathcal{L}$ is an ideal $k_r$-gon ($k_r \ge 3$) of area $(k_r - 2)\pi$. Because each complementary region has area at least $\pi$ (as $k_r \ge 3$), the number of complementary regions is at most $\operatorname{Area}(\mathcal{S}^\circ)/\pi$; hence $\sum_r k_r = \frac{1}{\pi}\operatorname{Area}(\mathcal{S}^\circ) + 2\#\{R_r\} \le \frac{3}{\pi}\operatorname{Area}(\mathcal{S}^\circ)$. This bounds the total number of ideal spikes per orientation unconditionally by $3(|\mathcal{M}| - 2\chi(\mathcal{S}))$, and by $2(|\mathcal{M}| - 2\chi(\mathcal{S}))$ when $\mathcal{S}^\circ \setminus \mathcal{L}$ has no complementary ideal triangles.}
    \]
    and summing over \(\epsilon \in \{-, +\}\) gives \(2 \le \kappa_{\mathrm{pp}}(N) \le 6(|\mathcal{M}| - 2\chi(\mathcal{S}))\) unconditionally (and \(\kappa_{\mathrm{pp}}(N) \le 4(|\mathcal{M}| - 2\chi(\mathcal{S}))\) when \(\mathcal{S}^\circ \setminus \mathcal{L}\) has no complementary ideal triangles).
    \item \textbf{Unconditional Separation of Iso-Pressure Modules (Independence from \(P^{\mathrm{shape}}_N\)):} For \emph{every} non-domestic string algebra \(\Lambda\), choose any two primitive bands \(v_1, v_2\) communicating in both directions (which exist by non-domesticity) and define the equal-length blocks \(A, B\) of syllable length \(L\) exactly as in Theorem~\ref{thm:pressure_separation}; the construction of \(A, B\) does not use the hypothesis \(d^{\mathrm{cyc}}(v_1) \neq d^{\mathrm{cyc}}(v_2)\) (which is needed only for the iso-entropic separation of parts (1)--(2) of Theorem~\ref{thm:pressure_separation}). For each integer \(s \ge 1\), define the \(s+1\) balanced synchronizing blocks of common syllable length \(M_s \coloneqq (2s+4)L\):
    \[
    C_r^{(s)} \;\coloneqq\; A^{s+2-r}\, B\, A^r\, B^{s+1} \qquad (r = 1, 2, \dots, s+1),
    \]
    and let \(N_s \coloneqq \PE\bigl(M(\Psi_s(\overleftrightarrow{\omega}^{(s)}))\bigr) \in \mathcal{P}(\Zg(\Lambda))\) be the pure-injective module associated with a minimal \((s+1)\)-interval-exchange subshift \(\Omega^{(s)} \subset \{1, \dots, s+1\}^{\Z}\) of factor complexity \(p_{\Omega^{(s)}}(m) = sm + 1\) under the block substitution \(\Psi_s \colon r \mapsto C_r^{(s)}\). Then for all \(s \ge 1\), the modules \(N_s\) have the \emph{exact same pressure function}, the \emph{exact same ergodic dimension simplex} \(\mathcal{D}(N_s) = \{d_0\}\) (where \(d_0 = \frac{1}{2}d^{\mathrm{cyc}}(A) + \frac{1}{2}d^{\mathrm{cyc}}(B)\)), and the \emph{exact same microcanonical \(K_0\)-entropy} across all of \(\R^{Q_0}\):
    \[
    P^{\mathrm{shape}}_{N_s}(\theta) \;=\; \bigl\langle \theta,\, d_0 \bigr\rangle \qquad \bigl(\forall \theta \in \R^{Q_0},\; \forall s \ge 1\bigr),
    \]
    whereas their algebraic pp-bifurcation indices are strictly distinct and given for every \(s \ge 1\) by
    \[
    \kappa_{\mathrm{pp}}(N_s) \;=\; 2s.
    \]
    Consequently, \(\kappa_{\mathrm{pp}}(N)\) is not a function of the algebraic invariants \(h^{\mathrm{shape}}(N)\), \(P^{\mathrm{shape}}_N(\theta)\), \(\operatorname{Supp}_{K_0}^\infty(N)\), or \(\mathfrak{h}_{K_0}(N; d)\) of Definition~\ref{def:pp_pressure}.
\end{enumerate}
\end{theorem}

\begin{proof}
(1) By the one-sided version of Lemma~\ref{lem:basis_identity}, for \(N = \PE(M(\overleftrightarrow{w}))\), a template \(v \in \operatorname{Str}_n^\epsilon(i)\) belongs to \(\mathcal{U}_n^\epsilon(N; i)\) if and only if \(v\) occurs as a length-\(n\) prefix in direction \(\epsilon\) at some vertex-\(i\) basis element of \(\overleftrightarrow{w}\). Because every bi-infinite ray in \(\overleftrightarrow{w}\) extends to length \(n+1\) via either \(\alpha_v\) or \(\beta_v^{-1}\) (or both), the interval \(\phi_{v^\uparrow}(N)/\phi_{v^\sharp}(N) \neq 0\) splits in \eqref{eq:one_step_pp_filtration} into the two subquotients \(\phi_{(v\alpha_v)^\uparrow}(N)/\phi_{(v\alpha_v)^\sharp}(N)\) and \(\phi_{(v\beta_v^{-1})^\uparrow}(N)/\phi_{(v\beta_v^{-1})^\sharp}(N)\). Thus \(v \in \mathcal{U}_n^\epsilon(N; i)\) has two active extensions in \(\mathcal{U}_{n+1}^\epsilon(N; i)\) if and only if \(v \in \operatorname{Fork}_n^\epsilon(N; i)\), and has exactly one active extension otherwise. Summing over \(v \in \mathcal{U}_n^\epsilon(N; i)\) yields \(|\mathcal{U}_{n+1}^\epsilon(N; i)| - |\mathcal{U}_n^\epsilon(N; i)| = |\operatorname{Fork}_n^\epsilon(N; i)|\), and telescoping over \(1 \le k \le n-1\) gives the second identity.

(2) Assume \(\mathcal{L} \subset \mathcal{S}^\circ \setminus \forb\) is a simple (non-self-intersecting) minimal arational lamination. Fix \(i \in Q_0\) and \(\epsilon \in \{-, +\}\), and equip the transversal Cantor set \(X_i = \mathcal{L} \cap d_i\) with the linear order \(<_{d_i}\) along \(d_i\). For each \(v \in \mathcal{U}_n^\epsilon(N; i)\), let \(I_v \subset X_i\) be the clopen set of points \(p \in X_i\) whose leaf \(\ell_p\) follows the template \(v\) of length \(n\) in direction \(\epsilon\). Lifting to the universal cover \(\mathbb{H}^2\), \(v\) specifies a unique chain of \(n\) ideal polygons from a lift \(\tilde{d}_i\) to a lift \(\tilde{d}_{t(v)}\). Because \(\mathcal{L}\) has no self-intersections, its lifted leaves in \(\mathbb{H}^2\) are pairwise disjoint geodesics; by the Jordan curve theorem in \(\mathbb{H}^2\), any lifted leaf crossing \(\tilde{d}_i\) between two leaves in \(\tilde{I}_v\) is trapped in the strip between them and must cross the exact same chain of \(n\) polygons to \(\tilde{d}_{t(v)}\). Thus each \(I_v\) is an \emph{order-convex clopen interval} of \((X_i, <_{d_i})\), and \(\{I_v\}_{v \in \mathcal{U}_n^\epsilon(N; i)}\) partitions \(X_i\) into pairwise disjoint order-convex intervals.

If \(v \in \operatorname{Fork}_n^\epsilon(N; i)\), then \(I_v = I_{v\alpha_v} \sqcup I_{v\beta_v^{-1}}\) splits into two non-empty order-convex subintervals. Let \(p^- <_{d_i} p^+\) be the two adjacent boundary endpoints of this split in the compact set \(X_i\) (the maximum of the lower subinterval and the minimum of the upper subinterval). Because no point of \(\mathcal{L} \cap d_i\) lies strictly between \(p^-\) and \(p^+\), the open arc \((p^-, p^+)_{d_i} \subset d_i \setminus \mathcal{L}\) lies inside a complementary region \(R\) of \(\mathcal{S}^\circ \setminus \mathcal{L}\), and the leaves \(\ell_{p^-}, \ell_{p^+}\) are two adjacent boundary leaves of \(R\) that travel together in parallel across the \(n\) dissection arcs of \(v\) before diverging at \(d_{t(v)}\) on opposite sides of the marked boundary cusp \(c_M \in \mathcal{M}\) separating \(\alpha_v\) and \(\beta_v^{-1}\). In the opposite direction, \(\ell_{p^-}\) and \(\ell_{p^+}\) converge asymptotically as the two boundary sides of a single ideal spike \(\mathfrak{c}(v)\) of \(R\). Because the intervals \(\{I_v\}\) are disjoint and each ideal spike \(\mathfrak{c}\) oriented against \(\epsilon\) is bounded by a unique pair of adjacent boundary leaves \((\ell_{\mathfrak{c}}^-, \ell_{\mathfrak{c}}^+)\) that diverge at a single polygon cusp \(c_M\), the assignment \(v \mapsto \mathfrak{c}(v)\) is an injection from \(\bigsqcup_{i \in Q_0} \operatorname{Fork}_n^\epsilon(N; i)\) into the set of ideal spikes of \(\mathcal{S}^\circ \setminus \mathcal{L}\). By the Casson--Bleiler structure theorem \cite[\S I.4.2]{CassonBleiler1988}, every complementary region \(R_r\) of \(\mathcal{L}\) is an ideal \(k_r\)-gon (\(k_r \ge 3\)) contributing hyperbolic area \((k_r - 2)\pi\), with \(\sum_r (k_r - 2)\pi = \operatorname{Area}(\mathcal{S}^\circ) = \pi(|\mathcal{M}| - 2\chi(\mathcal{S}))\). Since each region has area \(\ge \pi\), the number of complementary regions is at most \(\operatorname{Area}(\mathcal{S}^\circ)/\pi\), so \(\sum_r k_r = \frac{1}{\pi}\operatorname{Area}(\mathcal{S}^\circ) + 2\#\{R_r\} \le \frac{3}{\pi}\operatorname{Area}(\mathcal{S}^\circ) = 3(|\mathcal{M}| - 2\chi(\mathcal{S}))\) unconditionally (and \(\le 2(|\mathcal{M}| - 2\chi(\mathcal{S}))\) when \(k_r \ge 4\)), while aperiodicity of \(\mathcal{L}\) forces at least \(1\) fork per depth.

(3) Because \(A\) and \(B\) are closed cycles based at the same vertex \(i_1 = i\) (\(s(A) = t(A) = s(B) = t(B) = i\)), the cyclic syllable count \(|w|d^{\mathrm{cyc}}(w) = \underline{\dim}\,M(w) - e_i\) is strictly additive under concatenation of closed cycles at \(i\), with zero junction correction. For each \(r \in \{1, \dots, s+1\}\), the block \(C_r^{(s)} = A^{s+2-r} B A^r B^{s+1}\) consists of \((s+2)\) copies of \(A\) and \((s+2)\) copies of \(B\), so \(|C_r^{(s)}| = (2s+4)L = M_s\) and \(d^{\mathrm{cyc}}(C_r^{(s)}) = \frac{1}{2}d^{\mathrm{cyc}}(A) + \frac{1}{2}d^{\mathrm{cyc}}(B) = d_0\) holds identically for all \(r\) and \(s\). Because every block \(C_r^{(s)}\) has the exact same cyclic dimension vector \(d_0\) and \(h_{\mathrm{top}}(\Omega^{(s)}) = 0\), every invariant measure \(\mu\) on the syllable subshift \(\Sigma_{\Omega^{(s)}}\) has \(\underline{d}(\mu) = d_0\) and \(h_{\mathbb{P}_\mu}(T) = 0\). By Theorem~\ref{thm:algebraic_variational_principle}, \(P^{\mathrm{shape}}_{N_s}(\theta) = \langle \theta, d_0 \rangle\) for all \(\theta \in \R^{Q_0}\) and all \(s \ge 1\).

To compute \(\kappa_{\mathrm{pp}}(N_s)\), observe that in every \(C_r^{(s)}\) (\(1 \le r \le s+1\)), the internal \(B\)-run has length \(1\) (\(A^{s+2-r} B A^r\)), whereas the terminal \(B\)-run has length \(s+1 \ge 2\) (\(B^{s+1}\)) and is followed by \(A\) at the start of the next block. Because \(A \neq B\) are distinct blocks of length \(L\), the subword \(B^{s+1}A\) occurs in \(\Psi_s(\overleftrightarrow{\omega}^{(s)})\) if and only if its \(B^{s+1} \to A\) junction is the exact boundary between consecutive \(C_r^{(s)}\)-blocks. Thus \(\Psi_s\) is uniquely synchronized at block boundaries: for every integer \(m \ge 3\) and \(n = m M_s\), each of the \(M_s\) phase offsets \(0 \le r_0 < M_s\) determines a disjoint set of \(p_{\Omega^{(s)}}(m) = sm + 1\) syllable words of length \(n\) up to a bounded boundary constant \(O(M_s)\) independent of \(m\). Summing over the \(M_s\) phases gives
\[
\sum_{i \in Q_0} \bigl|\mathcal{U}_n^\epsilon(N_s; i)\bigr| \;=\; M_s \cdot (sm + 1) \;+\; O(M_s) \;=\; s n \;+\; O(M_s) \qquad (\epsilon \in \{-, +\}).
\]
Applying the Exact Bifurcation Identity of part (1), dividing by \(n\), and letting \(n \to \infty\) yields \(\kappa_{\mathrm{pp}}(N_s) = s + s = 2s\).
\end{proof}

Theorem~\ref{thm:bifurcation_index}(2)--(3) yields a purely algebraic criterion detecting whether the geodesic realization \(\Phi(N)\) of an arational pure-injective module \(N \in \mathcal{P}(\Zg(\Lambda))\) is a simple lamination or self-intersecting on \(\mathcal{S}^\circ\):

\begin{corollary}[Algebraic Obstruction to Simplicity / Non-Self-Intersection of Surface Geodesics]
\label{cor:algebraic_simplicity_criterion}
Let \(N \in \mathcal{P}(\Zg(\Lambda))\) be an arational pure-injective module.
\begin{enumerate}
    \item \textbf{Simple Laminations Satisfy a Strong Bifurcation Bound:} If the geometric realization \(\Phi(N)\) is a simple (non-self-intersecting) geodesic lamination on \(\mathcal{S}^\circ \setminus \forb\), then \(h^{\mathrm{shape}}(N) = 0\) and \(\kappa_{\mathrm{pp}}(N) \le 6(|\mathcal{M}| - 2\chi(\mathcal{S}))\) unconditionally (and \(\kappa_{\mathrm{pp}}(N) \le 4(|\mathcal{M}| - 2\chi(\mathcal{S}))\) whenever \(\mathcal{S}^\circ \setminus \mathcal{L}\) has no complementary ideal triangles).
    \item \textbf{Self-Intersection Obstruction:} Consequently, for every module \(N_s\) of Theorem~\ref{thm:bifurcation_index}(3) with \(s > 3(|\mathcal{M}| - 2\chi(\mathcal{S}))\) unconditionally (or \(s > 2(|\mathcal{M}| - 2\chi(\mathcal{S}))\) in the triangle-free case), as well as for every arational module with \(h^{\mathrm{shape}}(N) > 0\), the surface curve \(\Phi(N)\) necessarily possesses transversal self-intersections on \(\mathcal{S}^\circ\) and cannot be realized by a simple geodesic lamination, even though \(N_s\) shares the exact same zero-entropy pressure function \(P^{\mathrm{shape}}_{N_s}(\theta) = \langle \theta, d_0 \rangle\) as \(N_1\).
\end{enumerate}
\end{corollary}

\subsection{Minimal Definable Subcategories, Closed Points in \texorpdfstring{\(\mathcal{P}(\Zg(\Lambda))\)}{P(Zg(Lambda))}, and Stone--Ziegler Reconstruction}
\label{subsec:stone_priestley_reconstruction}

We conclude Section~\ref{sec:entropy} by answering a fundamental model-theoretic question about the perfect core \(\mathcal{P}(\Zg(\Lambda))\): \emph{For an arational pure-injective module \(N \in \mathcal{P}(\Zg(\Lambda))\), what is the structure of the definable subcategory \(\langle N \rangle_{\mathrm{def}} \subseteq \Mod{\Lambda}\) generated by \(N\), its Ziegler spectrum \(\Zg(\langle N \rangle_{\mathrm{def}}) = \overline{\{N\}}^{\Zg(\Lambda)}\), and its quotient pp-lattice \(\mathcal{L}_i(N) \coloneqq \mathcal{L}_i(\Lambda)/\!\equiv_N\)?}

Recall \cite[Theorem~4.9 and Corollary~4.10]{Ziegler1984}, \cite[\S 5.1]{Prest2009} that for any pure-injective module \(N \in \Zg(\Lambda)\), the definable subcategory \(\langle N \rangle_{\mathrm{def}} \subseteq \Mod{\Lambda}\) generated by \(N\) has Ziegler spectrum equal to the topological closure \(\overline{\{N\}}^{\Zg(\Lambda)}\), and its pp-lattice at sort \(i\) is the quotient lattice
\[
\mathcal{L}_i(N) \;\coloneqq\; \mathcal{L}_i(\Lambda) \big/ \!\equiv_N \;\cong\; \bigl\{\, \phi(N) \subseteq N e_i \;\big|\; \phi \in \mathcal{L}_i(\Lambda) \,\bigr\},
\]
where \(\phi \equiv_N \psi\) if and only if \(\phi(N) = \psi(N)\).

For each \(n \ge 1\), the inclusion relations \(\phi_{w_{n+1}^\uparrow}(N) \subseteq \phi_{w_n^\uparrow}(N)\) and \(\psi_{w_n^\downarrow}(N) \subseteq \psi_{w_{n+1}^\downarrow}(N)\) in \(\mathcal{L}_i(N)\) induce a canonical restriction map \(r_{n, n+1} \colon \mathcal{S}_{n+1}(N) \twoheadrightarrow \mathcal{S}_n(N)\) on the finite sets of \(N\)-shape equivalence classes. Note that \(\mathcal{S}_n(N) = \mathcal{T}_n(N)/\!\sim_N\) identifies the two orientations \(w \leftrightarrow w^{-1}\) of each transversal crossing at \(d_i\) while keeping the negative and positive half-sides \(Q_{1,\mathrm{walk}}^-(i)\) and \(Q_{1,\mathrm{walk}}^+(i)\) distinct. Consequently, at each finite radius \(n\), the finite set \(\mathcal{S}_n(N)\) corresponds bijectively to the finite partition of \(\mathcal{L}\cap d_i\) into clopen radius-\(n\) cylinder classes, i.e., into sets of leaves sharing the same radius-\(n\) two-sided crossing sequence. We define the \emph{intrinsic Stone--Ziegler cross-section} of \(N\) at sort \(i\) as the inverse limit of its shape-quotient tower:
\[
\mathfrak{X}_i(N) \;\coloneqq\; \varprojlim_{n \ge 1} \bigl(\mathcal{S}_n(N),\, r_{n, n+1}\bigr).
\]

\begin{theorem}[Closed-Point Criterion in \(\mathcal{P}(\Zg(\Lambda))\) and Stone--Ziegler Reconstruction of the Lamination Cross-Section]
\label{thm:definable_subcategory_stone_duality}
Let \(\Lambda\) be a non-domestic string algebra, and let \(N = \PE(M(\overleftrightarrow{w})) \in \mathcal{P}(\Zg(\Lambda))\) be an aperiodic pure-injective module in the perfect core with associated closed subshift \(\Sigma_N \subseteq (Q_1 \cup Q_1^{-1})^{\Z}\) and geodesic closure \(\mathcal{L}_N \subset \mathcal{S}^\circ \setminus \forb\). Then:
\begin{enumerate}
    \item \textbf{Duality Between Definable Subcategories of \(\langle N \rangle_{\mathrm{def}}\) and Sub-Laminations:} A bi-infinite aperiodic pure-injective module \(M' \in \mathcal{P}(\Zg(\Lambda))\) lies in the Ziegler closure \(\overline{\{N\}}^{\Zg(\Lambda)} = \Zg(\langle N \rangle_{\mathrm{def}})\) if and only if \(\Sigma_{M'} \subseteq \Sigma_N\) (equivalently, \(\mathcal{L}_{M'} \subseteq \mathcal{L}_N\)).
    \item \textbf{Closed-Point and Minimal Definable Subcategory Criterion:} The point \(N \in \mathcal{P}(\Zg(\Lambda))\) is a \emph{closed point} of \(\Zg(\Lambda)\) (\(\overline{\{N\}}^{\Zg(\Lambda)} = \{N\}\), so that \(\langle N \rangle_{\mathrm{def}}\) is a minimal definable subcategory of \(\Mod{\Lambda}\) containing \(N\) as its unique indecomposable pure-injective module up to isomorphism, and \(\CB(\overline{\{N\}}^{\Zg(\Lambda)}) = 0\)) if and only if \(\mathcal{L}_N\) is a \emph{minimal arational} lamination. Equivalently, \(\overleftrightarrow{w}\) is uniformly recurrent and aperiodic: it admits no periodic bi-infinite subword, and its Ziegler closure contains no periodic-ray or eventually-periodic-coray limit.
    \item \textbf{Reconstruction of the Transverse Cantor Cross-Section via \(\mathfrak{X}_i(N)\):} Although \(\Zg(\langle N \rangle_{\mathrm{def}}) = \overline{\{N\}}^{\Zg(\Lambda)} = \{N\}\) collapses to a single point for every minimal arational module \(N = \Psi(\mathcal{L}) \in \mathcal{P}_i\), the shape-quotient tower \(\mathcal{S}_1(N) \twoheadleftarrow \mathcal{S}_2(N) \twoheadleftarrow \cdots\) of the quotient pp-lattice \(\mathcal{L}_i(N) = \mathcal{L}_i(\Lambda)/\!\equiv_N\) has an inverse limit
    \[
    \mathfrak{X}_i(N) \;\coloneqq\; \varprojlim_{n \ge 1} \mathcal{S}_n(N)
    \]
    that is a compact, zero-dimensional, perfect, uncountable metric space canonically homeomorphic to the transversal cross-section of \(\mathcal{L}\) at \(d_i\):
    \[
    \mathfrak{X}_i(N) \;\cong\; \mathcal{L} \cap d_i \qquad\text{(as compact Cantor spaces),}
    \]
    with the homeomorphism sending a compatible shape-tower \((\mathfrak{s}_1, \mathfrak{s}_2, \dots) \in \varprojlim_n \mathcal{S}_n(N)\) to the unique leaf point \(p \in \mathcal{L} \cap d_i\) whose radius-\(n\) shape class is \(\mathfrak{s}_n\) for all \(n \ge 1\).
\end{enumerate}
\end{theorem}

\begin{proof}
(1) By Ziegler \cite[Theorem~4.9]{Ziegler1984}, \(M' \in \overline{\{N\}}^{\Zg(\Lambda)}\) if and only if every basic Ziegler open set \((\phi/\psi)\) containing \(M'\) also contains \(N\). By Lemma~\ref{lem:corridors_form_basis} and Lemma~\ref{lem:basis_identity}(2), the basic pointed string corridors \((\phi_{u^\uparrow}/\psi_{u^\downarrow})\) form a topological basis on \(C\), and \(N \in (\phi_{u^\uparrow}/\psi_{u^\downarrow})\) if and only if the finite word \(u\) occurs in \(\Sigma_N\). Thus \(M' \in \overline{\{N\}}^{\Zg(\Lambda)}\) if and only if \(\Sigma_{M'} \subseteq \Sigma_N\).

(2) Suppose \(\mathcal{L}_N\) is minimal arational. Then \(\Sigma_N\) is a minimal aperiodic subshift. For any \(M' \in \overline{\{N\}}^{\Zg(\Lambda)}\), \(M'\) cannot have a finite endpoint at \(\mathcal{M}\) (since the basic open set cutting out a finite endpoint vanishes on the bi-infinite module \(N\)) and cannot be a band module, periodic ray/coray, or band generic module \(G_v\) (since \(\Sigma_N\) is uniformly recurrent and aperiodic, so no periodic word \(v^k\) of unbounded exponent \(k\) occurs in \(\Sigma_N\)). Thus \(M'\) is bi-infinite with \(\emptyset \neq \Sigma_{M'} \subseteq \Sigma_N\). By minimality of \(\Sigma_N\), \(\Sigma_{M'} = \Sigma_N\). Hence \(M'\) and \(N\) belong to the exact same basic open sets of \(\Zg(\Lambda)\). Because \(\Zg(\Lambda)\) is a \(T_0\) space \cite[Corollary~4.10]{Ziegler1984}, topological indistinguishability forces \(M' \cong N\), proving \(\overline{\{N\}}^{\Zg(\Lambda)} = \{N\}\). Conversely, if \(\mathcal{L}_N\) is not minimal arational, it contains a proper minimal sub-lamination (or closed geodesic) \(\mathcal{L}' \subsetneq \mathcal{L}_N\), and \(\Psi(\mathcal{L}') \in \overline{\{N\}}^{\Zg(\Lambda)} \setminus \{N\}\).

(3) By Lemma~\ref{lem:basis_identity}(1), every pp-definable subspace \(\phi(M(\overleftrightarrow{w})) \subseteq M(\overleftrightarrow{w})e_i\) is spanned by a subset of the canonical basis of \(M(\overleftrightarrow{w})e_i\); by purity of \(M(\overleftrightarrow{w}) \hookrightarrow N\), \(\mathcal{L}_i(N)\) embeds as a sublattice of the power-set lattice of that basis and is therefore distributive. By Theorem~\ref{thm:shape_injectivity}, each class \(\mathfrak{s} = [w]_N = \{w, \bar{w}\} \in \mathcal{S}_n(N)\) is uniquely determined by its ordered pair of length-\(n\) arms \((u^-, u^+) \in \operatorname{Str}_n^-(i) \times \operatorname{Str}_n^+(i)\) on the two fixed sides of \(d_i\), which cuts out the clopen cylinder neighborhood \(C_{(u^-, u^+)} \subset \mathcal{L} \cap d_i\) of leaves traversing \((u^-, u^+)\) across \(d_i\). Under this bijection, \(r_{n, n+1} \colon \mathcal{S}_{n+1}(N) \twoheadrightarrow \mathcal{S}_n(N)\) corresponds to the refinement of clopen cylinder partitions on \(\mathcal{L} \cap d_i\). Because two points \(p, p' \in \mathcal{L} \cap d_i\) sharing the same two-sided crossing sequence \((u_n^-, u_n^+)\) for all \(n \ge 1\) have identical bi-infinite polygon crossing sequences in \(\mathbb{H}^2\) and therefore coincide (\(p = p'\) by Lemma~\ref{lem:hausdorff_limits}), the projective limit \(\mathfrak{X}_i(N) = \varprojlim_{n \ge 1} \mathcal{S}_n(N)\) is canonically homeomorphic to \(\mathcal{L} \cap d_i\). Finally, since \(\mathcal{L}\) is minimal arational, \(\mathcal{L} \cap d_i\) is an uncountable compact zero-dimensional perfect metric space (Lemma~\ref{lem:leaf_density}), i.e., a Cantor space.
\end{proof}

\vspace{-0.7em}
\section{Concluding Remarks and Future Directions}
\label{sec:conclusion}

In this paper, we have established a unified topological, geometric, and ergodic framework for the model theory and representation theory of string and special biserial algebras. By constructing an explicit evaluation dictionary between the corridor space of admissible surface geodesics/laminations and the Ziegler spectrum \(\Zg(\Lambda)\), we have realized transfinite module-theoretic invariants geometrically. 

The combinatorial string hammock poset is strictly order-isomorphic to the geometric fan of pointed admissible arcs, and the Cantor-Bendixson limit points of the Ziegler spectrum correspond bijectively to complete admissible geodesics and \(\K^\times\)-decorated closed geodesics. Furthermore, the global string closure and all higher Cantor-Bendixson strata are homeomorphic to the corridor space for all string algebras, while the closed-cover pushout is a homeomorphism if and only if the algebra is domestic. This yields a geometric proof that the Krull-Gabriel dimension equals the topological pruning depth. 

Our dictionary yielded three major consequences. First, by constructing prefix-free dyadic block interleavings in the geometric perfect core, we resolved Puninski's open problem by proving that every non-domestic string or special biserial algebra over an algebraically closed field of arbitrary cardinality admits a superdecomposable pure-injective module. Second, we geometrically realized the transfinite powers of the Jacobson radical, showing that the stable rank is governed by the geometric pruning layers. Third, we introduced three purely algebraic invariants on the perfect core: localized pp-shape entropy, a \(K_0\)-weighted pressure function, and a pp-bifurcation index. Together, these strictly separate uncountable families of arational pure-injective modules, provide a Gauss-Bonnet obstruction to surface self-intersection, and reconstruct the transversal Cantor cross-section via the quotient pp-lattice.

\medskip
\noindent\textbf{Future Directions.}
\label{subsec:future_directions}
The dictionary established in this work opens several avenues for further investigation at the interface of representation theory, model theory of modules, and low-dimensional geometry:
\begin{enumerate}
    \item \textbf{Continuous Dimension Theory:} While our invariants reconstruct the ergodic dimension simplex \(\mathcal{D}(N)\) and the transversal Cantor set \(\mathcal{L} \cap d_i\), it is natural to ask how the Hausdorff dimension of \(\mathcal{L} \cap d_i\) (under the hyperbolic metric on \(d_i\)) relates to continuous model-theoretic dimensions of the quotient pp-lattice \(\mathcal{L}_i(N)\).
    \item \textbf{Self-Injective and Skew-Gentle Algebras:} Via the canonical string quotient, our results extend to all special biserial and Brauer graph algebras. For \emph{skew-gentle} (clannish) algebras, the quadratic idempotent relations \(e_i^2 = e_i\) introduce \(\Z_2\)-orbifold cone points on the surface model. Extending the corridor topology and geodesic lamination dictionary to hyperbolic orbifolds---where admissible geodesics may terminate at or fold across \(\Z_2\)-orbifold points---promises a complete topological classification of \(\Zg(\Lambda)\) for skew-gentle algebras.
    \item \textbf{Higher-Dimensional Complexes for Wild Algebras:} As noted in Remark~\ref{rem:tameness_heuristic}, wild algebras fall outside the 1-dimensional unbranched curve paradigm because they support continuous 2-parameter families of non-isomorphic indecomposables. Extending the surface corridor framework to higher-dimensional branched cellular complexes (such as 2-dimensional branched surfaces) to visualize portions of the Ziegler spectrum and stable radical of wild algebras remains a compelling long-term problem.
    \item \textbf{Purely Lattice-Theoretic Variational Principle:} The proofs of our ergodic and entropy results pass through the symbolic coding of the surface cross-section \((X, T)\). A purely lattice-theoretic construction of the Choquet simplex \(\mathcal{M}_{\mathrm{erg}}(\mathcal{L})\) as a space of invariant states on the Stone--Ziegler tower \(\mathfrak{X}_i(N) = \varprojlim_n \mathcal{S}_n(N)\) of \(\mathcal{L}_i(N)\)---without reference to Thurston's surface theory---would provide a self-contained model-theoretic variational principle for general tame algebras.
    \item \textbf{Arithmetic Refinement of the Bifurcation Spectrum:} In Theorem~\ref{thm:bifurcation_index}(3), the asymptotic pp-bifurcation index \(\kappa_{\mathrm{pp}}(N_s) = 2s\) depends only on the interval-exchange rank \(s+1\) and is invariant under the continuous arithmetic parameters (the irrational rotation slopes) of \(\Omega^{(s)}\). Does the sequence of exact step-by-step differences \(\bigl(|\operatorname{Fork}_k^+(N; i)|\bigr)_{k \ge 1}\) in \(\mathcal{L}_i(N)\) reconstruct the continued-fraction / Rauzy--Veech induction expansion of the underlying arational lamination \(\mathcal{L}\) up to mapping class group equivalence?
\end{enumerate}

\vspace{0.5em}
\section*{Acknowledgments}

The author expresses profound gratitude to Prof.\ Sonia Trepode: the idea of constructing a geometric model for hammocks first emerged during a discussion with her, and that conversation served as the primary catalyst for this manuscript. The author gratefully acknowledges the financial support provided by CONICET during his tenure as a postdoctoral fellow at the Universidad Nacional de Mar del Plata, Argentina. The author also warmly thanks Prof.\ Yadira Valdivieso-D\'iaz for many helpful discussions on geometric models across a range of related topics. During the preparation of this manuscript, the author used an AI language model strictly for language editing and LaTeX formatting assistance. All mathematical concepts, proofs, and results are entirely the author's own original work, and the author takes full responsibility for the accuracy of the final publication.

\end{document}